\documentclass[a4paper,11pt]{amsart}

\usepackage{graphicx} 

\usepackage[margin=50pt]{geometry}
\usepackage{lipsum}

\usepackage{CJKutf8}

\usepackage{mathpazo}
\usepackage{textcomp}
\usepackage{nicefrac} 
\usepackage{xfrac}  
\usepackage{stmaryrd}

\usepackage{enumitem}
\usepackage{mathtools}
\usepackage[utf8]{inputenc}
\usepackage[all,matrix,arrow,curve]{xy}
\usepackage{amssymb}
\usepackage{amsfonts}
\usepackage{amsthm}
\usepackage{adjustbox}

\usepackage{float}
\usepackage{tabularx}

\usepackage{latexsym,pifont}
\usepackage{epsfig}
\usepackage{rotating}

\usepackage{ifpdf}
\ifpdf
\usepackage{epstopdf}
\usepackage{hyperref}
\else
\usepackage[hypertex]{hyperref}
\fi

\usepackage{amsmath,mathrsfs,extarrows,MnSymbol}
\usepackage{tikzsymbols}
\usepackage[all]{xy}

\xyoption{knot}

\usepackage{tikz}
\usetikzlibrary{matrix,arrows}

\def\la{\label}

\usepackage{bbding}

\usepackage{tikz-cd}
  \usepackage{tikz}
  \usepackage{multicol}
  \usepackage[linguistics]{forest}
        \usepackage{tikz-qtree}
        
        \forestset{sn edges/.style={for tree={parent anchor= north, child anchor=south}}}

\def\be{\begin{equation}}
\def\ee{\end{equation}}
\def\ba{\begin{eqnarray}}
\def\ea{\end{eqnarray}}
\def\id{\mathrm{id}}

\def\xcG{{\cal G}}

\def\g{\mathfrak{g}}

\newcommand{\p}{\mathrm p^{ \geq 2}}

\newtheorem{theorem}{Theorem}[section]
\newtheorem*{theorem*}{Theorem}
\newtheorem{proposition}[theorem]{Proposition}

\newtheorem{remark}[theorem]{Remark}
\newtheorem{example}[theorem]{Example}
\newtheorem{definition}[theorem]{Definition}
\newtheorem{corollary}[theorem]{Corollary}
\newtheorem{lemma}[theorem]{Lemma}

\newtheorem{propanition}[theorem]{Proposition/Definition}

\newtheorem{coronition}[theorem]{Corollary/Definition}

\newtheorem{conjecture}[theorem]{Conjecture}
\newtheorem{convention}[theorem]{Convention}

\let\a=\alpha \let\b=\beta
\let\G=\Gamma

\newcommand{\alxydim}[2]{\begin{aligned}\xymatrix#1{#2}\end{aligned}}

\newcommand{\brem}{\begin{remark}}
\newcommand{\erem}{\end{remark}}
\newcommand{\beg}{\begin{example}}
\newcommand{\eeg}{\end{example}}
\newcommand{\bedef}{\begin{definition}}
\newcommand{\exdef}{\end{definition}}
\newcommand{\berop}{\begin{proposition}}
\newcommand{\eerop}{\end{proposition}}
\newcommand{\belem}{\begin{lemma}}
\newcommand{\elem}{\end{lemma}}
\newcommand{\bethe}{\begin{theorem}}
\newcommand{\ethe}{\end{theorem}}
\newcommand{\becor}{\begin{corollary}}
\newcommand{\ecor}{\end{corollary}}
\newcommand{\beroof}{\noindent\begin{proof}}
\newcommand{\eroof}{\end{proof}}
\newcommand{\becon}{\begin{convention}}
\newcommand{\econ}{\end{convention}}
\newcommand{\efact}{\begin{flushright}$\checkmark$\end{flushright}\end{Fact}}
\newcommand{\becj}{\begin{conjecture}}
\newcommand{\ecj}{\begin{flushright}$\boxtimes$\end{flushright}\end{conjecture}}

\newcommand{\barr}{\begin{array}}
\newcommand{\earr}{\end{array}}
\newcommand{\ben}{\begin{enumerate}}
\newcommand{\een}{\end{enumerate}}
\newcommand{\bit}{\begin{itemize}}
\newcommand{\eit}{\end{itemize}}

\newcommand{\qq}{\begin{eqnarray}}
\newcommand{\qqq}{\end{eqnarray}}

\newcommand{\nn}{\nonumber}

\newcommand{\ovl}[1]{\overline{#1}}
\newcommand{\unl}[1]{\underline{#1}}

\newcommand{\Reqref}[1]{Eq.\,\eqref{#1}}

\newcommand{\Rxcite}[2]{Ref.\,\cite[#1]{#2}}

\newcommand\void[1]{}

\newcommand{\gt}[1]{\mathfrak{#1}}

\def\cA{\mathcal{A}}

\def\cE{\mathcal{E}}

\def\xcG{\mathcal{G}}

\def\cK{\mathcal{K}}

\def\cO{\mathcal{O}}

\def\cR{\mathcal{R}}

\def\cV{\mathcal{V}}

\def\xcX{\mathcal{X}}

\def\xcA{\mathscr{A}}

\def\xcD{\mathscr{D}}
\def\xcE{\mathscr{E}}
\def\xcF{\mathscr{F}}
\def\xcG{\mathscr{G}}

\def\xcK{\mathscr{K}}
\def\xcL{\mathscr{L}}
\def\xcM{\mathscr{M}}

\def\xcP{\mathscr{P}}

\def\xcR{\mathscr{R}}
\def\xcS{\mathscr{S}}

\def\xcU{\mathscr{U}}
\def\xcV{\mathscr{V}}

\def\xcX{\mathscr{X}}

\def\R{\mathbf{R}}

\def\t{\mathbf{t}}

\def\bB{{\mathbb{B}}}

\def\bE{{\mathbb{E}}}

\def\bN{{\mathbb{N}}}

\def\bR{{\mathbb{R}}}

\def\bV{{\mathbb{V}}}

\def\a{\alpha}
\def\b{\beta}
\def\g{\gamma}
\def\G{\Gamma}
\def\d{\delta}

\def\vep{\varepsilon}

\def\la{\lambda}

\def\om{\omega}
\def\Om{\Omega}

\def\si{\sigma}
\def\Si{\Sigma}

\def\t{\tau}

\def\z{\zeta}

\def\ggt{\gt{g}}

\newcommand{\sfd}{{\mathsf d}}

\newcommand{\sfI}{{\mathsf I}}

\newcommand{\sfJ}{{\mathsf J}}

\newcommand{\sfM}{{\mathsf M}}

\newcommand{\sfP}{{\mathsf P}}

\newcommand{\sfS}{{\mathsf S}}

\newcommand{\sfT}{{\mathsf T}}

\newcommand{\txA}{{\rm A}}

\newcommand{\txd}{{\rm d}}

\newcommand{\txg}{{\rm g}}
\newcommand{\txG}{{\rm G}}

\newcommand{\txH}{{\rm H}}

\newcommand{\txm}{{\rm m}}

\newcommand{\txP}{{\rm P}}
\newcommand{\breP}{{\rm \breve{P}}}

\newcommand{\txS}{{\rm S}}
\newcommand{\txT}{{\rm T}}

\newcommand{\txV}{{\rm V}}

\def\Cv{\v{C}}

\def\id{{\rm id}}
\newcommand{\pr}{{\rm pr}}

\def\too{\longrightarrow}
\def\ev{{\rm ev}}

\newcommand{\grpd}[2]{\hspace{-5pt}\alxydim{@C=.5cm@R=1.cm}{ #1 \ar@<.25ex>[r] \ar@<-.25ex>[r] & #2}\hspace{-3pt}}

\def\obj{{\rm Ob}}

\def\morf{{\rm Mor}}
\def\1morf{1{\rm -Mor}}
\def\2morf{2{\rm -Mor}}
\def\dim{{\rm dim}}

\def\ker{{\rm ker}}

\def\End{{\rm End}}

\newcommand{\Id}{{\rm Id}}
\def\Inv{{\rm Inv}}

\newcommand{\Gr}{{\rm {\bf Gr}}}
\newcommand{\BisGr}{{\rm Bisec({\bf Gr})}}

\def\p{\partial}

\newcommand{\Diff}{{\rm Diff}}

\def\emb{\hookrightarrow}

\def\bd1{{\boldsymbol{1}}}
\def\brd0{{\boldsymbol{0}}}

\def\rk{{\rm rk}}

\def\Ad{{\rm Ad}}

\newcommand\MCL{\theta_{\rm L}}
\newcommand\MCR{\theta_{\rm R}}

\def\x{\times}
\def\ox{\otimes}

\def\lx{{\hspace{-0.04cm}\ltimes\hspace{-0.05cm}}}
\def\rx{\rtimes}
\newcommand{\fibx}[2]{\hspace{1pt}{}_{#1}\hspace{-3pt}\x_{#2}\hspace{-4pt}}
\def\ract{\vartriangleleft}
\def\lact{\vartriangleright}
\def\mact{\blacktriangleleft}
\def\mlact{\blacktriangleright}
\def\ulact{\,\unl{\vartriangleright}\,}
\def\must{\stackrel{!}{=}}

\def\rstr{\big\vert}

\newcommand\colo{\hspace{-2pt}:}
\newenvironment{proproof}{
	\noindent{\em Proof of Proposition:}\ }{\hfill $\blacksquare$ \newline}
\newenvironment{lemproof}{
\noindent{\em Proof of Lemma:}\ }{\hfill $\blacksquare$ \newline}
\newenvironment{corproof}{
\noindent{\em Proof of Corollary:}\ }{\hfill $\blacksquare$ \newline}

\title{Principaloid bundles}

\author[T. Strobl]{Thomas Strobl}
\address{Thomas Strobl:\ Institut Camille Jordan,\
Universit\'e Claude Bernard Lyon 1,\ 69622 Villeurbanne,\ France  $\qquad \emph{and} \:$
The Erwin Schr\"odinger International Institute for Mathematics and Physics,\ University of Vienna,\ Boltzmanngasse 9A,\ 1090 Wien,\ Austria}
\email{stroblATmath.univ-lyon1.fr}

\author[R.R. Suszek]{Rafa\l~ R. Suszek}
\address{Rafa\l ~R.\ Suszek:\ Department of Mathematical Methods in Physics,\ Faculty of Physics, University of Warsaw,\ ul.\ Pasteura 5,\ 02-093 Warszawa,\ Poland}
\email{suszekATfuw.edu.pl}

\begin{document}
\tikzset{mystyle/.style={xshift = - 2.4ex, yshift = - 1.5ex, align=left}}

\tikzset{mystyle2/.style={xshift = - 1.7ex, yshift = - 2.ex, align=left}}

\tikzset{mystyle3/.style={xshift = + 2.4ex, yshift = - 1.5ex, align=left}}

\maketitle

\begin{abstract}
We present a novel generalisation of principal bundles---principaloid bundles:\ These are fibre bundles $\pi\hspace{-2pt}: \mathscr{P}\to B$ where the typical fibre is the arrow manifold $\mathscr{G}$ of a Lie groupoid $\hspace{-5pt}\alxydim{@C=.5cm@R=1.cm}{ \mathscr{G} \ar@<.25ex>[r] \ar@<-.25ex>[r] & M}\hspace{-3pt}$ and the structure group is reduced to the latter’s group of bisections.\ Each such bundle canonically comes with a bundle map $\mathscr{D} \hspace{-2pt}: \mathscr{P} \to \mathscr{F}$ to another fibre bundle $\mathscr{F}$ over the base $B$,\ with typical fibre $M$.\ Examples of principaloid bundles include ordinary principal ${\rm G}$-bundles,\ obtained for $\mathscr{G} := \hspace{-5pt}\alxydim{@C=.5cm@R=1.cm}{ {\rm G} \ar@<.25ex>[r] \ar@<-.25ex>[r] & \bullet}\hspace{-3pt}$,\ bundles associated to them,\ obtained for action groupoids $\mathscr{G} := {\rm G}\ltimes M$,\ and general fibre bundles if $\mathscr{G}$ is a pair groupoid. 
    
While $\pi$ is far from being a principal $\mathscr{G}$-bundle,\ we prove that $\mathscr{D}$ is one.\ Connections on the principaloid bundle $\pi$ are thus required to be $\mathscr{G}$-invariant Ehresmann connections.\ In the three examples mentioned above,\ this reproduces the usual types of connection for each of them.\ In a local description over a trivialising cover $\{O_i\}$ of $B$,\ the connection gives rise to Lie algebroid-valued objects living over bundle trivialisations $\{O_i \times M\}$ of $\mathscr{F}$.\ Their behaviour under bundle automorphisms,\ including gauge transformations,\ is studied in detail. 
    
Finally,\ we construct the Atiyah--Ehresmann groupoid $\,\hspace{-5pt}\alxydim{@C=.5cm@R=1.cm}{ \mathrm{At}(\mathscr{P}) \ar@<.25ex>[r] \ar@<-.25ex>[r] & \mathscr{F}}\hspace{-3pt}$ which governs symmetries of $\mathscr{P}$,\ this time mapping distinct $\mathscr{D}$-fibres to one another in general.\ It is a fibre-bundle object in the category of Lie groupoids,\ with typical fibre $\hspace{-5pt}\alxydim{@C=.5cm@R=1.cm}{ \mathscr{G} \ar@<.25ex>[r] \ar@<-.25ex>[r] & M}\hspace{-3pt}$ and base $\hspace{-5pt}\alxydim{@C=.5cm@R=1.cm}{ B\times B \ar@<.25ex>[r] \ar@<-.25ex>[r] & B}\hspace{-3pt}$.\ We show that those of its bisections which project to bisections of its base are in a one-to-one correspondence with automorphisms of $\pi$.
\end{abstract}

\tableofcontents

\section{Introduction}
\label{sec:introduction}

Principal bundles $\sfP$ and their associated bundles
as well as connections on them are of great importance in differential geometry and are also omnipresent in large parts of theoretical physics.\ In this paper,\ we propose and study a new generalisation of this framework where the structure Lie group of a principal bundle is replaced by a Lie groupoid in a particular way.

There are several equivalent ways of defining principal bundles.\ One is by saying that $\sfP$ is a fibre bundle over $\Sigma$ with typical fibre $\txG$,\ where $\txG$ is a (finite-dimensional) Lie group,\ for which the structure group $\Diff(\txG)$ is reduced to $L(\txG)$,\ the group of left translations of $\txG$.\ In other words,\ the fibre bundle $\pi_\sfP \colo \sfP \to \Sigma$ is equipped with an equivalence class of atlases with the property that all transition maps within one such atlas,\ when restricted to any point in the overlap of trivialisation charts in the base,\ take values in  $L(\txG)\subset \Diff(\txG)$.

In this description,\ it is evident that $\sfP$ carries a fibrewise right $\txG$-action.\ Since the Lie subgroup of $\Diff(\txG)$ consisting of diffeomorphisms of $\txG$ which commute with $L(\txG)$ is given by $R(\txG)$,\ the group of right translations of $\txG$,\ there is the following elegant alternative:\ A principal $\txG$-bundle $\sfP$ is given simply by a free and proper right $\txG$-action on a manifold $\sfP$. All the data above can be recovered from here.

This opens different ways for a generalisation to Lie groupoids,\ from which we here choose the first one as our starting point.\footnote{It deserves to be mentioned that the concept of a principaloid bundle, developed here in full detail, first arose in discussions of the first author with Oleksii Kotov in 2005.}

To set the notation,\ throughout this paper,\ a Lie groupoid is denoted by $\xcG$ or $\grpd{\xcG}{M}$,\ where $\xcG$ is the arrow manifold and $M$  the object manifold.\ Furthermore,\ $\bB \equiv{\rm Bisec}(\xcG)$ denotes the (generally infinite-dimensional) group of bisections of $\xcG$.\ Its elements $\b\in\bB$ are sections of the source map $s$ of $\xcG$,\ {\it i.e.},\ maps $\b\colo M\to\xcG$ such that $s\circ\b=\id_M$,\ with the additional property:\ $t\circ\b\in\Diff(M)$,\ written in terms of the target map $t$ of $\xcG$.\ $\bB$ can act on $\xcG$ from the left or from the right,\ essentially by left and right translations.\footnote{Whenever it cannot be avoided,\ maps with an infinite-dimensional manifold as domain or codomain---as for the $\bB$-actions above in the case where $M$ is not  a finite set---we assume the mapping to be smooth in the sense of Kriegl and Michor \cite{Kriegl:1997}.}\ We denote the corresponding two subgroups of $\Diff(\xcG)$ by $L(\bB)$ and $R(\bB)$.\ For example,\ given arbitrary $\beta \in \bB$ and $g\in \xcG$,\ the diffeomorphism $L_\b\in L(\bB)$ maps $g$ as follows: $L_\beta(g)= \beta(t(g)).g$,\ where the dot denotes the multiplication in $\xcG$. (For further details on this and other statements regarding Lie groupoids, their bisections, and related notions, see Appendices.)\smallskip

\noindent\textbf{Definition.} A \emph{principaloid bundle},\ or principaloid $\xcG$-bundle,\ is a fibre bundle $\pi_\xcP \colo \xcP \to \Si$ with typical fibre $\xcG$,\ whose structure group is reduced to $L(\bB)$ (or a subgroup thereof).
\smallskip

The commonly adopted generalisation of principal bundles, on the other hand,  
leads to what is called \emph{principal groupoid bundles},\ see,\ {\it e.g.}, \cite{MacKenzie:1987,Moerdijk:1991},\ or principal $\xcG$-bundles if the groupoid is specified.\ This notion is significantly different from the one considered here:\ In general,\ principal $\xcG$-bundles are not even fibre bundles.\ However, they also play an important role in the setting of principaloid bundles,\ as secondary or induced structures, for which reason we first recall this notion.

The definition of a principal $\xcG$-bundle is based on a particular choice of a $\xcG$-action on the total space of the bundle.\ This is already subtler than in the case of Lie groups.\ Indeed,\ for a general Lie groupoid $\grpd{\xcG}{M}$,\ multiplication is defined only for composable arrows,\ {\it i.e.},\ $g_1.g_2$ is defined only for elements $g_1,g_2 \in \xcG$ such that $t(g_2)=s(g_1)$.\ Consequently,\ for an action of a groupoid on a manifold $\breP$,\ one needs an additional map $\mu \colo \breP \to M$,\ called the moment map.\ An element $g\in \xcG$ then maps $\mu^{-1}(\{t(g)\})$ to $\mu^{-1}(\{s(g)\})$.\ Requiring this (right) action along $M$ to be---in an appropriate sense---proper and free turns $\breP$ into a principal $\xcG$-bundle,\ with a smooth base $\Sigma$ which is canonically identifiable with $\breP/\xcG$.\ The setting is summarised in the following diagram: 
\qq\label{diag:CechP}
\alxydim{@C=1cm@R=1.5cm}{ \breP \ar@{->>}[d]_{\pi_\breP} \ar[rd]^{\mu} & & \xcG \ar@{=>}[ld] \\ \Si\equiv\breP/\xcG & M & }\,.
\qqq
By the principality of the 
action,\ one finds that the fibre of $\pi_\breP$ at a given point $\sigma \in \Si$ is diffeomorphic to the $t$-fibre $t^{-1}(\{\mu(\pi_\breP^{-1}(\{\sigma\}))\})$ in $\xcG$.\ Since $t \colo \xcG \to M$ is only required to be a surjective submersion,\ $\breP$ does not have a typical fibre in general. 

The groupoid $\xcG$ comes with two surjective submersions:\ $s$ (the source map) and $t$ (the target map) to its object manifold $M$,\ and so it should not come as a surprise that every principaloid bundle is canonically endowed with two important maps:\ One modelled on $s$ as the moment map $\mu \colo \xcP \to M$,\ and one modelled on $t$ as the base projection of a surjective submersion $\xcD \colo \xcP \to \xcF$,\ 
a bundle map to another fibre bundle $\xcF$ over $\Sigma$,\ with typical fibre $M$,\ which is canonically induced by $\xcP$.  We shall show in the present paper that these data combine in the following manner:\medskip 

\noindent\textbf{Theorem A.} {\it Every principaloid bundle $\xcP \to \Sigma$ is a principal $\xcG$-bundle over $\xcF$ in a canonical way,
\qq \label{diag:diagram2}
\alxydim{@C=1cm@R=1.5cm}{ \xcG \ar@{->>}[d]_{t} \ar@{^{(}.>}[r] & \xcP \ar@{->>}[d]_{\xcD} \ar[rd]^{\mu} & & \xcG \ar@{=>}[ld] \\ M \ar@{^{(}.>}[r] & \xcF\equiv\xcP/\xcG \ar@{->>}[d]_{\pi_\xcF} & M &  \\ & \Si & }\,.
\qqq
The $\xcG$-action preserves the $\pi_\xcP$-fibres,\ where $\pi_\xcP=\pi_\xcF\circ\xcD$.}\medskip

The present setting is not only natural,\  but also subsumes other known constructions as special cases:\ On the one hand, for a Lie groupoid over a point, $\grpd{\txG}{\bullet}$, we recover a principal $\txG$-bundle $\sfP \to \Si$. Note that here the group of bisections $\bB\cong \txG$ so that the structure group becomes the usual one,\ $\txG$,\ and $\xcP\cong \txP$. Furthermore, the $\xcG$-action in \eqref{diag:diagram2} reduces to the standard right $\txG$-action on $\txP$.

On the other hand,\ if we choose $\xcG = M \times M$,\ the pair groupoid of a connected manifold $M$,\ then $\bB \cong \mathrm{Diff}(M)$ and $\xcP \cong \xcF \times M$,\ where $\xcF$ is a general fibre bundle over $\Si$ with typical fibre $M$,\ and an unrestricted structure group $\Diff(M)$.\ Note that here the right $\xcG$-action of \eqref{diag:diagram2}  corresponds to diffeomorphisms of the second factor in $\xcF \times M$,\ leaving $\xcF$ untouched.

Finally,\ let the structure Lie groupoid be an action Lie groupoid, $\xcG := \txG \ltimes M$,\ defined for an arbitrary $\txG$-action $\lambda \colo \txG \times M \to M$ on a manifold $M$.\ Consider a principaloid $\xcG$-bundle in which the structure group $\bB$ is reduced to $\bB_0=\{\ (g,\cdot) \;\vert\; g\in\txG \ \} \cong \txG$.\ Then,\ $\xcP \cong \sfP \x M$,\ where $\sfP$ is an ordinary principal $\txG$-bundle over $\Si$.\ Moreover,\ the bundle $\xcF$ becomes none other than the bundle associated to $\sfP$ by the action $\la$.\ In fact,\ the map $\xcD \colo \txP \times M\cong\xcP \to \xcF$ turns out to coincide with the corresponding quotient map in the usual construction of the associated bundle.

Thus,\ principaloid bundles unify principal $\txG$-bundles,\ associated bundles,\ as well as general fibre bundles under a common roof---while simultaneously offering a rich class of additional examples that invite further exploration.

There are also other ways of characterising a principaloid bundle worth mentioning here.\ 
We first observe that a right-$\xcG$-equivariant diffeomorphism of $\xcG$,\ for the right action given by right translations,\ always comes---as shall be shown---in the form of a left-multiplication by a bisection,\ $\Diff_{r(\xcG)}(\xcG)=L(\bB)$.\ Furthermore,\ every $\xcG$-action on a manifold $N$ gives rise to a $\bB$-action on $N$.\ While the former action is, for a given element, defined only on submanifolds of $N$,\ mapping $\mu$-fibres to $\mu$-fibres,\ the latter one is an honest group action,\ defined over all of $N$.\ We then have the following results:

\medskip 
\noindent\textbf{Theorem B.} There are the following \emph{alternative} ways of characterising a principaloid bundle (in the second case with the additional assumption that every $g\in \xcG$ admits a $\beta \in \bB$ such that $\beta(s(g))=g$): 
\bit
\item A fibre bundle in the category of right $\xcG$-modules with typical fibre given by the canonical right $\xcG$-module $\xcG$ and with the trivial $\xcG$-action on the base. 
\item A fibre bundle in the category of $\bB$-spaces with typical fibre given by the canonical right $\bB$-space $\xcG$,\ with the trivial $\bB$-action on the base,\ 
together with an equivalence class of $\bB$-equivariant atlases.
\eit
\medskip 
Such definitions are closer to the definition of an ordinary principal $\txG$-bundle which uses a principal action of the structure group $\txG$.\ Note that due to the identity $\bB \cong \txG$,\ which holds true for the groupoid  $\grpd{\txG}{\bullet}$,\ there arise two generalisations of such a $\txG$-action for general groupoids:\ by a (partially defined) $\xcG$-action,\ or by a (globally defined) $\bB$-action.\smallskip

Connections are central objects in the study of fibre bundles.\ We define a connection $\Theta$ on a principaloid $\xcG$-bundle as a $\xcG$-invariant Ehresmann connection $\Theta$ on $\xcP \to \Si$. Since $\xcG$ acts on $\xcP$ by (diffeomorphically) mapping $\mu$-fibres to $\mu$-fibres,\ see Diagram  \eqref{diag:diagram2},\ the above condition on $\Theta$ is meaningful only if the horizontal subbundle $\txH\xcP \subset \txT \xcP$ lies inside $\ker \txT \mu$,\ which we therefore implicitly assume when demanding that $\Theta$ be $\xcG$-invariant.\ The $\xcG$-invariance further implies that $\Theta$ induces an Ehresmann connection $\Theta_\xcF$ on $\xcF \cong \xcP/\xcG$.

Since every principaloid $\xcG$-bundle $\xcP$ is also a $\bB$-space,\ a (compatible) connection on $\xcP$ is also $\bB$-invariant with respect to the globally defined right $\bB$-action $\xcR$ on $\xcP$.\ Under the assumption about $\xcG$ mentioned in Theorem B---that every arrow $g$ has a global bisection $\beta$ passing through it---the $\xcR(\bB)$-invariance of $\Theta$ implies its $\varrho(\xcG)$-invariance, where $\varrho$ denotes the right $\xcG$-action of Diagram \eqref{diag:diagram2}.\ (The fact that the horizontal distribution must lie in the kernel of $\txT \mu$ is a consequence of large stabiliser groups of the $\bB$-action here.)

In the three special examples of principaloid bundles discussed above,\ such a connection is in a one-to-one correspondence with (i) for $\xcG = \txG$:\ an ordinary principal connection,\ (ii) for $\xcG=M\x M$:\ an unconstrained Ehresmann connection on the general fibre bundle $\xcF$,\ and (iii) for $\xcG=\txG \ltimes M$ and structure group $\bB_0$:\ an ordinary principal connection together with the connection on the associated bundle $\xcF$ canonically induced by it. 
 
For ordinary principal $\txG$-bundles $\sfP$,\ $\txG$-invariant Ehresmann connections can be described by Lie algebra-valued 1-forms on the total space of the bundle.\ This does not generalise to principaloid bundles as Lie algebroid-valued 1-forms.\ The reason is simple:\ An Ehresmann connection on a bundle $\xcP \to \Sigma$ can be understood as a projector $\Theta\colo\txT\xcP\to\txV\xcP$ onto the vertical subbundle $\txV\xcP\subset \txT\xcP$,\ modelled on $\txT\xcG$.\ Now,\ for a group,\ $\xcG \cong \txG$,\ the tangent bundle trivialises as $\txT \txG \cong \txG \times \ggt$,\ where $\ggt = \mathrm{Lie}(\txG)$ is the Lie algebra of $\txG$,\ which,\ as a vector space,\ we can identify with the tangent space of $G$ at the identity.\ The Lie algebroid $E \to M$ of a Lie groupoid $\grpd{\xcG}{M}$,\ on the other hand,\ contains only those tangent vectors at the identity bisection of $\xcG$ which are tangent to the $s$-fibres, $E\cong \Id^*\ker\,\txT s$.
 
Nevertheless,\ we shall show that a connection 1-form $\Theta$ locally gives rise to 1-forms valued in the tangent Lie algebroid of the structure groupoid $\xcG$,\ generalising the local Lie algebra-valued 1-forms $A \in \Omega^1(O,\ggt)$ in a trivialisation of $\sfP$ over $O \subset \Si$.\ The situation is,\ however,\ more intricate than in the case of ordinary principal bundles,\ essentially due to the appearance of the bundle $\xcF\to\Si$ with typical fibre $M$ of non-zero dimension in Diagram \eqref{diag:diagram2},\ whereas for principal $\txG$-bundles $\xcF$ reduces to $\Si$.

In a local trivialisation of $\xcP$ over $O \subset \Si$,\ the connection 1-form $\Theta$ on $\xcP$ is in a one-to-one correspondence with a section $\txA \in \Gamma(\pr_1^*T^*O \otimes \pr_2^*E)$ over the corresponding trivialisation $O\times M$ of $\xcF$.\ Hence,\ at a point $(\sigma,m)\in O\times M$,\ one has, after natural idendtifications,
\qq\label{Aintro}   \txA(\sigma,m) \; \in \;T_\sigma^*O \otimes E_m \,,
\qqq
where,\ as above,\ $E \to M$ denotes the Lie algebroid of $\grpd{\xcG}{M}$ and $E_m$ its fibre over $m\in M$.

The Lie algebroid-valued $\txT O$-foliated 1-forms $A$ are glued together according to an affine transformation law on overlaps of trivialisation domains,\ which generalises the familiar transformation law for local connection 1-forms on ordinary principal bundles in a natural way.\ Specific formulas are technically involved,\ so we refer the reader to Definition \ref{def:principoidle-conn-cech} and Theorem \ref{thm:loc-data-conn} in the main text below for further details.\smallskip 

Automorphisms of a principaloid bundle $\xcP$ are bundle automorphisms which commute with the right $\xcG$-action on $\xcP$.\ There is a groupoid behind these automorphisms,\ the Atiyah groupoid $\mathrm{At}(\xcP)$.\ The object manifold of this groupoid is composed of the $\xcD$-fibres in $\xcP$ and is thus parametrised by points in the bundle $\xcF$.\ In total,\ we get two commuting  principal groupoid actions on $\xcP$, as summarised in the following diagram:
\qq\label{diag:Trident}
\alxydim{@C=.75cm@R=1cm}{ & & \xcG \ar@{^{(}.>}[d] & & \\ {\rm At}(\xcP) \ar@{=>}[rd] & & \xcP \ar[dl]_{\xcD} \ar[rd]^{\mu} \ar[dd]^{\pi_\xcP} & & \xcG \ar@{=>}[ld] \\ & \xcF & & M & \\ & & \Si & &}\,,
\qqq
where we dropped maps to $\Si$ or $\Si \times \Si$ on the left for aesthetic reasons, so as to arrive at a diagram resembling a trident. Note that while $\xcP/\xcG = \xcF$, the action of ${\rm At}(\xcP)$ projects nontrivially 
to the base $\Sigma$ of $\xcP$ and one obtains $\xcP/\mathrm{At(\xcP)} = M$.

As a bundle over $\Si$,\ ${\rm At}(\xcP)$ has typical fibre $\xcP$,\ where the \Cv ech 1-cocycle of $\xcP$ is used for the gluing on overlaps by employing the \emph{right}  $\bB$-action on $\xcP$.\ It can also be viewed as a fibre bundle over $\Si \times \Si$,\ in which case the typical fibre is $\xcG$---and gluings by means of the \Cv ech 1-cocycle are performed using the left-multiplication over the first copy of $\Si$ and the right-multiplication over the second copy of $\Si$.\ In fact,\ one can view ${\rm At}(\xcP)$ as a  fibre-bundle object in the category of Lie groupoids as follows:
\qq\label{diag:category}
\alxydim{@C=.5cm@R=1.cm}{\xcG \: \ar@{^{(}.>}[r] & {\rm At}(\xcP) \ar[d]^{\pi} \\ & {\rm Pair}(\Si)}\,.
\qqq
Adding the kernel of the groupoid morphism $\pi$ into the picture,\ one obtains an exact sequence of Lie groupoids ${\rm Ad}(\xcP)\emb{\rm At}(\xcP)\twoheadrightarrow{\rm Pair}(\Si)$. Here,\ by definition,\ ${\rm Ad}(\xcP)$ is the preimage of the identity bisection $\Id(\Si)\subset\Si\x\Si$ along the epimorphism $\pi$ of Lie groupoids;\ it is a fibre bundle over $\Sigma$ with typical fibre $\xcG$,\ where the gluing  is performed by conjugation or the adjoint action of the transition 1-cocyle.\ Differentiation of this leads to the corresponding exact sequence of Lie algebroids,\ studied as particular $Q$- or $\mathrm{dg}$-bundles in \cite{Kotov:2007nr}.\ Objects of $\mathrm{At}(\xcP)$ are $\xcD$-fibres of $\xcP$ sitting over points in $\xcF$,\ and,\ as such,\ also over points in $\Sigma$.\ A groupoid element of $\mathrm{At}(\xcP)$,\ characterised by a triple
$(\sigma_2,g_1,\sigma_1)$, $g_1\in \xcG$, $\si_1,\si_2 \in \Sigma$ in a local chart,\ then maps a point $p \in \xcP$,\ given by $(\sigma_1,g_2)$ for some $g_2 \in t^{-1}(\{s(g_1)\}) \subset \xcG$ in a local chart,\ to the point $(\sigma_2,g_1.g_2)$ in this description. In this manner,\ the action of ${\rm At}(\xcP)$ on $\xcP$ covers the canonical left action of ${\rm Pair}(\Si)$ on its base $\Si$.

Every groupoid action gives rise to an action of the group of bisections,\ where bisections can be also viewed as  submanifolds in the groupoid mapping to $M$ diffeomorphically by means of the source and target maps.\ Applying this logic to the right wing of Diagram \eqref{diag:Trident},\ {\it i.e.},\ to the $\xcG$-action on $\xcP$,\ one recovers the (right) $\bB$-action on $\xcP$.\ For the left wing,\ {\it i.e.},\ the $\mathrm{At(\xcP)}$-action on $\xcP$ along $\xcF$,\ we obtain a characterisation or description of the automorphisms.\ Find below an expanded version of the trident diagram 
\qq
\alxydim{@C=1.1cm@R=1.5cm}{ \Ad(\xcP) \ar@{^{(}~>}[r] \ar[dd] & {\rm At}(\xcP) \ar@{=>}[rrd] \ar@{~>>}[dd]_\pi & \b_\Phi(\xcF) \ar@{_{(}->}[l] & & \xcP \ar[dl]_{\xcD} \ar[rd]^{\mu} \ar[dd]^{\pi_\xcP} & & \b(M) \ar@{^{(}->}[r] & \xcG \ar@{=>}[lld] \\ & & & \xcF \ar[dd]_(.4){\pi_\xcF} \ar[ul]_{\b_\Phi} & & M \ar[ur]^{\b} & & \\ \Si\cong\Id(\Si) \ar@{^{(}->}[r] & \Si\x\Si \ar@{=>}[rrd] & (f,\id_\Si)(\Si) \ar@{_{(}->}[l] & & \Si \ar[dl]^{\id_\Si} & & & \\ & & & \Si \ar[ul]_{(f,\id_\Si)} & & & & }\,.\cr \label{diag:Beast}
\qqq 
where the exact sequence of Lie groupoids mentioned before is  highlighted by the wiggly arrows.\ Now recall that objects of $\mathrm{At(\xcP)}$ are fibres of $\xcD \colo \xcP \to \xcF$.\ A $\pi_\xcF$-fibre $\pi_\xcF^{-1}(\{\sigma_1\})$ assembles all $\xcD$-fibres sitting over the point $\sigma_1 \in \Sigma$ into the $\pi_\xcP$-fibre over that point,\ which is isomorphic to the groupoid $\xcG$.\ If one takes an arbitrary bisection of the groupoid $\grpd{\mathrm{At(\xcP)}}{\xcF}$,\ it will map two $\xcD$-fibres over $\sigma_1$ to $\xcD$-fibres sitting over,\ in general,\ two different points in $\Sigma$.\ This is,\ however,\ in contradiction with the fact that an automorphism of the bundle $\pi_\xcP \colo \xcP \to \Sigma$ must map $\pi_\xcP$-fibres to $\pi_\xcP$-fibres.\ This can be cured if one considers only bisections $\beta_\Phi$ of the groupoid $\grpd{\mathrm{At(\xcP)}}{\xcF}$ which are $\pi$-projectable,\ {\it i.e.},\ those bisections which induce diffeomorphisms $f$ on the base $\Sigma$ of $\xcP$. 

Note that such complications disappear in the case of an ordinary principal bundle $\sfP$, resulting from ours for the groupoid $\grpd{\txG}{\bullet}$, in which case the trident degenerates to a W-shaped diagram, where on the left $\xcF$ is replaced by the base $\Sigma$. In such a case,\ there is a one-to-one correspondence between arbitrary bisections of the Atiyah groupoid $\mathrm{At(\sfP)}$ and automorphisms of  $\sfP$.\medskip
\noindent\textbf{Theorem C.} {\it Every principaloid bundle $\xcP \to \Sigma$ gives rise to a fibre-bundle object \qq\nn
\left({\rm At}(\xcP) \Rightarrow \xcF\right) \stackrel{\pi}{\too} \left({\rm Pair}(\Sigma) \Rightarrow \Si \right)
\qqq
in the category of Lie groupoids.\ The $\pi$-projectable bisections of its total space ${\rm At}(\xcP) \Rightarrow \xcF$ are in a one-to-one correspondence with automorphisms of $\xcP$,\ and those projecting to ${\rm Id}(\Si)$ with its vertical autormorphisms,\ {\it i.e.},\ with gauge transformations of $\xcP$.\ These gauge transformations act on local connection data  $A_i$ over $O_i \times M$---see \eqref{Aintro}---subordinate to a cover $\{O_i\}_{i \in I}$ of $\Si$ as in \eqref{gaugetrafo} in the main text. 

Every principaloid bundle $\xcP \to \Sigma$ is a principal ${\rm At}(\xcP)$-bundle over $M$ in a canonical way.\ Together with the structure of a principal $\xcG$-bundle of Theorem B, this combines into a principal$({\rm At}(\xcP),\xcG)$-bibundle object in ${\rm {\bf Bun}}(\Si)$,\ described by the Trident Diagram \eqref{diag:Trident}.} 

\medskip 

We foresee numerous applications of the present formalism.\ One of these concerns gauge theories lying outside the traditional framework of principal bundles,\ such as the Poisson sigma model \cite{Ikeda:1994,Schaller:1994es} or the generalisations of Yang--Mills theories provided in \cite{S04, KS15}.\ Another class would be the study of conditions under which higher-geometric and -cohomological structures defined on the manifold $M$ descend to the bundle $\xcF$ upon pullback to $\xcP$ along the moment map $\mu$---{\it i.e.},\ the modelling of descent of objects such as $n$-gerbes and field theories determined by the associated differential characters to characteristic foliations of Lie groupoids through construction of $\xcG$-equivariant structures on these objects,\ along the lines of \cite{Gawedzki:2010rn,Gawedzki:2012fu}.

\bigskip 

\noindent{\bf Acknowledgements:} Both authors gratefully acknowledge the kind support and hospitality of Institut Camille Jordan at Universit\'e Claude Bernard in Lyon,\ the Department of Mathematical Methods in Physics at the University of Warsaw,\ and the Erwin Schr\"odinger Institute of Mathematics and Physics in Vienna,\ where various parts of this work were carried out.\ In the final stage of work,\ the second author was also supported by an STSM mobility grant from the CA21109 COST Action CaLISTA.

The authors are indebted to Peter W.\ Michor for reading the manuscript and kindly sharing his remarks with them.\ They also thank Alexander Schmeding for his valuable comments on the structure of the group of bisections of a Lie groupoid. T.S.\ thanks Oleksii Kotov for important discussions at an early stage of this project and Anton Yu.\ Alekseev for remarks on the Introduction.

\section{Principaloid bundles}
\label{sec:definition}

Let $\xcG \sim \grpd{\xcG}{M} \sim \Gr$ denote a Lie groupoid and $\bB := \BisGr$ its group of bisections.\ There is a canonical right $\xcG$-action on $\xcG$ by right translations $r$,\ and three canonical $\bB$-actions:\ the left-multiplication $L$,\ the right-multiplication $R$,\ and the conjugation $C$.\ By means of the target map $t\colo \xcG \to M$,\ $L$ induces also an action of $\bB$ on $M$, which we denote by $t_*$ and call the shadow action.\ We refer to Appendix \ref{app:Liegrpd} for further details,\ notation,\ and conventions.

\subsection{Three definitions,\ and relations between them}

Let us consider three generalisations of the definition of a principal $\txG$-bundle,\ with $\txG$ a Lie group,\ and so simultaneously a Lie groupoid over a point $\grpd{\txG}{\bullet}$ wih ${\rm Bisec}(\txG)\cong\txG$.\ We start with one which emphasises the existence of a right action of the structure groupoid.
\bedef\label{def:principaloid-G-bndl}
A {\bf principaloid $\xcG$-bundle} is a fibre-bundle object in the category of right $\xcG$-modules with a fibrewise right action of $\xcG$ and typical fibre given by the canonical right $\xcG$-module $\xcG$.
\exdef
\brem\label{rem:princ-G-expl}
Let us unwrap this first definition:\ It describes a fibre bundle $(\xcP,\Si,\xcG,\pi_\xcP)$ with a base given by a manifold $\Si$,\ to which the total space $\xcP$ maps by a surjective submersion $\pi_\xcP\colo \xcP\to\Si$ (the base projection).\ The total space $\xcP$ carries the structure of a right $\xcG$-module $(\xcP,\mu,\varrho)$,\ see Def.\,\ref{def:gr-mod},\ with an action 
\qq\nn
\varrho\colo\xcP{}_{\mu}\hspace{-3pt}\x_t\hspace{-1pt}\xcG\too\xcP,\ (p,g)\longmapsto\varrho(p,g)\equiv\varrho_g(p)\equiv p\mact g
\qqq
which preserves $\pi_\xcP$-fibres.\ The base admits an open cover $\{O_i\}_{i\in I}$ with local trivialisations 
\qq\label{eq:PG-loc-triv}
\xcP\t_i\colo \pi_\xcP^{-1}(O_i)\xrightarrow{\ \cong\ } O_i\x\xcG\,,\qquad i\in I\,,
\qqq
which are $\xcG$-equivariant diffeomorphisms in the sense of Def.\,\ref{def:gr-mod}.\ In other words,\ for all $\si\in O_i$ and $g,h\in\xcG$,\ we have the equivalence
\qq\nn 
(g,h)\in\xcG{}_s\hspace{-3pt}\x_t\hspace{-1pt}\xcG\qquad\Longleftrightarrow\qquad(\xcP\t_i^{-1}(\si,g)),h)\in\xcP{}_{\mu}\hspace{-3pt}\x_t\hspace{-1pt}\xcG\,,
\qqq 
and the identity 
\qq\nn
\varrho_h\bigl(\xcP\t^{-1}(\si,g)\bigr)=\xcP\t^{-1}\bigl(\si,r_h(g)\bigr)\equiv\xcP\t^{-1}(\si,g.h)
\qqq
holds true whenever either condition is satisfied.
\erem
There is another definition,\ which takes as its base the choice of the structure group and its realisation on the typical fibre.\ This definition is,\ then,\ more directly related to the abstract notion of principal bundle as a fibre $\txG$-bundle with the specific choice of the realisation of the structure group on the typical fibre,\ namely:\ $L(\txG)\subset\Diff(\xcG)$,\ see \cite{Kolar:1993}.
\bedef\label{def:principoidle-cech}
A {\bf preprincipaloid bundle} is a fibre bundle with typical fibre $\xcG$ which admits a reduction of the structure group to the realisation $L(\bB)\subset\Diff(\xcG)$ of the group $\bB$ of bisections of $\xcG$ by left-multiplication.\ A {\bf principaloid bundle} is a preprincipaloid bundle together with a choice of such a reduction. 
\exdef
\brem\label{rem:defprincipoidleunwrap}
A preprincipaloid bundle is a fibre bundle $(\widetilde\xcP,\Si,\xcG,\pi_\xcP)$ which admits a choice of local trivialisations 
\qq\label{eq:P-loc-triv}
\widetilde\xcP\t_i\colo \pi_{\widetilde\xcP}^{-1}( O_i)\xrightarrow{\ \cong\ } O_i\x\xcG\,,\qquad i\in I\,,
\qqq
associated with an open cover $\{ O_i\}_{i\in I}$ of its base $\Si$,\ for which the corresponding transition maps take values in the group of bisections $\bB$ of $\xcG$,\ {\it i.e.},\ for every $(\si,g)\in O_{ij}\x\xcG$,
\qq\nn
\widetilde\xcP\t_i\circ\widetilde\xcP\t_j^{-1}(\si,g)=\bigl(\si,\b_{ij}(\si)\lact g\bigr)\,,
\qqq
where
\qq\nn
\b_{ij}\in C^\infty(O_{ij},\bB)\,.
\qqq
Accordingly,\ there exists a bundle isomorphism 
\qq\label{eq:clutch-xcP}
\widetilde\xcP\cong\bigsqcup_{i\in I}\,( O_i\x\xcG)/\sim_{ L{}_{\b_{\cdot\cdot}}}\,,
\qqq
which puts the bundle on the right-hand side in the role of a model of $\widetilde\xcP$.

Fixing such a trivialisation,\ up to equivalence,\ constitutes a choice---of which there could be many---turning $\widetilde\xcP$ into a principaloid bundle.\ Here,\ two trivialisations are considered equivalent if their joint refinement has transition functions in $L(\bB)$.
\erem
Given that left-multiplication commutes with right-multiplication,\ it is clear that once a reduction of the structure group to $L(\bB)$ has been fixed,\ $\bB$ acts on $\widetilde\xcP$ from the \emph{right} by bundle automorphisms which cover the identity on $\Si$ (termed vertical):\ Such an action is induced from right-mutiplications $R$ in local trivialisations as they glue smoothly over intersections of the trivialisation charts.\ Note that different reductions,\ {\it i.e.},\ different equivalence classes of trivialisations,\ in general give rise to different $\bB$-actions.

The $\bB$-action induced on the principaloid bundle $\widetilde\xcP$ can be presented as a group homomorphism 
\qq\label{eq:R-act-Poid}
\xcR\colo \bB\too{\rm Aut}_{{\rm {\bf Bun}}(\Si)}(\widetilde\xcP)_{\rm vert}\,.
\qqq
We have,\ for all $(\si,g)\in O_i\x\xcG$ and $i\in I$,
\qq\nn
\xcR_\b\circ\widetilde\xcP\t_i^{-1}(\si,g):=\widetilde\xcP\t_i^{-1}\bigl(\si,R_\b(g)\bigr)\equiv\widetilde\xcP\t_i^{-1}\bigl(\si,g\ract\b\bigr)\,,
\qqq
and so trivialisations within an equivalence class are (tautologically) $\bB$-equivariant.

\brem
Principaloid bundles can be viewed as bundles (non-canonically) associated to principal $\bB$-bundles by the action $L\colo\bB\to\Diff(\xcG)$.\ Such principal $\bB$-bundles are represented---in virtue of the classic Clutching Theorem---by (classes of) the corresponding transition 1-cocycles $\{\b_{ij}\}_{(i,j)\in I^{\x 2}_\cO}$.\ It is not hard to see that principaloid bundles,\ while offering us the practical advantage of finite dimensionality,\ actually carry the same \v Cech-cohomological information as the underlying infinite-dimensional objects.\ Indeed,\ the 1-cocycle $\{\b_{ij}\}_{(i,j)\in I^{\x 2}_\cO}$ can readily be recovered from its realisation $\,L_{\b_{ij}}\,$ on the typical fibre $\xcG$ through evaluation on the identity bisection $\Id(M)\subset\xcG$,
\qq\nn
L_{\b_{ij}(\si)}(\Id_m)\equiv\bigl(\b_{ij}(\si)\bigr)\bigl(t\bigl(\Id_m\bigr)\bigr).\Id_m=\bigl(\b_{ij}(\si)\bigr)(m)\,.
\qqq
\erem
Finally,\ one may ask to what extent the existence of the $\bB$-action is a defining property of a principaloid bundle.\ In order to study this question,\ let us introduce another notion,\ adapted to the right action of $\bB$.
\bedef\label{def:Bequiv-Gbndl}
A {\bf $\xcG$-fibred $\bB$-bundle} $\widehat\xcP$ over $\Si$ is a fibre-bundle object in the category of right $\bB$-spaces,\ with a fibrewise right action of $\bB$ and typical fibre given by the canonical right $\bB$-space $\xcG$,\ on which $\bB$ acts by right-multiplication. 
\exdef

We shall now examine relations between the three definitions given above.\ In the first step,\ we establish an equivalence between the first two definitions of principaloid bundles.
\bethe\label{thm:principaloid-GB-equiv}
Every principaloid $\xcG$-bundle carries a canonical structure of a principaloid bundle,\ and {\it vice versa}.
\ethe

\beroof
Consider a principaloid $\xcG$-bundle $\xcP$ with $\xcG$-equivariant trivialisations $\xcP\t_i,\ i\in I$ as in Remark \ref{rem:princ-G-expl}.\ As their smooth inverses are also (automatically) $\xcG$-equivariant,\ we obtain,\ for arbitrary $\si\in O_i\cap O_j\equiv O_{ij}$ and $ g\in\xcG$,
\qq\nn
\xcP\t_i\circ\xcP\t_j^{-1}(\si, g)=\bigl(\si,t_{ij}(\si,g)\bigr)\,,
\qqq
for some smooth maps $t_{ij}\colo O_{ij}\to\Diff_{r(\xcG)}(\xcG),\ \si\mapsto t_{ij}(\si,\cdot)$ with values in
\qq\nn
\Diff_{r(\xcG)}(\xcG):=C^\infty_\xcG(\xcG,\xcG)\subset \Diff(\xcG)\,,
\qqq
the {\bf group of $r(\xcG)$-equivariant diffeomorphisms of} $\xcG$,\ see Def.\,\ref{def:gr-mod}.\ This subgroup of $\Diff(\xcG)$ is identified in the following proposition.
\berop\label{prop:rGequiv-LB}
For every Lie groupoid $\xcG$,
\qq\nn
\Diff_{r(\xcG)}(\xcG)=L(\bB)\,.
\qqq
\eerop
\begin{proproof}
The relation $L(\bB)\subset\Diff_{r(\xcG)}(\xcG)$ follows straightforwardly from the fact that $L(\bB)$ is implemented by left translations on $\xcG$,\ which commute with right translations.

For the reverse relation,\ note that for every $\Phi\in\Diff_{r(\xcG)}(\xcG)$ and $g\in\xcG$,\ we have the identity
\qq\nn
\Phi(g)\equiv\Phi\bigl(\Id_{t(g)}.g\bigr)=\Phi\bigl(\Id_{t(g)}\bigr).g\,.
\qqq
Define a smooth map
\qq\nn
\varphi\colo M\too\xcG,\ m\longmapsto\Phi(\Id_m)
\qqq
to rewrite the above as
\qq\label{eq:Phinalus}
\Phi(g)=\varphi\bigl(t(g)\bigr).g\,.
\qqq
This further implies 
\qq\nn
(s\circ\varphi)\bigl(t(g)\bigr)=s\bigl(\Phi(g).g^{-1}\bigr)=s\bigl(g^{-1}\bigr)=t(g)\equiv\id_M\bigl(  t( g)\bigr)\,,
\qqq
which yields the identity
\qq\nn
s\circ\varphi=\id_M
\qqq
in virtue of the surjectivity of $t$.

Up to now,\ we have not fully taken into account the assumption of invertibility of $\Phi$,\ which further constrains $\varphi$.\ By a standard argument,\ the $\xcG$-equivariance of $\Phi$ implies the same property of its inverse,\ and so---in the light of the reasoning presented above---there exists a smooth map $\widetilde\varphi\colo M\to\xcG$ such that,\ for every $g\in\xcG$,
\qq\nn
\Phi^{-1}(g)=\widetilde\varphi\bigl(t(g)\bigr).g\,.
\qqq
We then obtain equalities
\qq\nn
g&=&\bigl(\Phi^{-1}\circ\Phi\bigr)(g)=\widetilde\varphi\bigl(t\bigl(\varphi\bigl(t(g)\bigr).g\bigr)\bigr).\varphi\bigl(t(g)\bigr).g=\widetilde\varphi\bigl((t\circ\varphi)\bigl(t(g)\bigr)\bigr).\varphi\bigl(t(g)\bigr).g\,,\cr\cr
g&=&\bigl(\Phi\circ\Phi^{-1}\bigr)(g)=\varphi\bigl((t\circ\widetilde\varphi)\bigl(t(g)\bigr)\bigr).\widetilde\varphi\bigl(t(g)\bigr).g\,,
\qqq
valid for every $g\in\xcG$,\ or,\ equivalently,
\qq\nn
&\widetilde\varphi\bigl((t\circ\varphi)(m)\bigr).\varphi(m)=\Id_m\,,&\cr\cr
&\varphi\bigl((t\circ\widetilde\varphi)(m)\bigr).\widetilde\varphi(m)=\Id_m\,,&
\qqq
valid for every $m\in M$.\ Upon evaluating the target map on both sides of the two identities,\ we arrive at the conclusion that $t\circ\widetilde\varphi$ is a smooth inverse of $t\circ\varphi$,\ and so,\ in particular,
\qq\nn
t\circ\varphi\in\Diff(M)\,.
\qqq

From our hitherto analysis,\ the mapping $\varphi$ emerges as a bisection of $\xcG$,\ and \Reqref{eq:Phinalus} yields
\qq\nn
\Phi\equiv L_{\varphi}\,,
\qqq
which is exactly the postulated result.
\end{proproof}

In virtue of the above proposition,\ the transition maps $t_{ij}$ of the principaloid $\xcG$-bundle $\xcP$ take the special form  
\qq\nn
t_{ij}(\si)=L_{\b_{ij}(\si)}\,,
\qqq
and so right-$\xcG$-equivariant trivialisations determine a reduction of the structure group of the bundle as in Def.\,\ref{def:principoidle-cech}.

For the converse claim,\ consider a principaloid bundle $(\widetilde\xcP,\Si,\xcG,\pi_{\widetilde\xcP})$ with a choice of reduction as in Remark \ref{rem:defprincipoidleunwrap}.\ For the structure of a right $\xcG$-module on $\widetilde\xcP$,\ define smooth maps 
\qq\nn
\mu_i\colo\pi_{\widetilde\xcP}^{-1}(O_i)\too M, \,\widetilde\xcP\t_i^{-1}(\si,g)\longmapsto s(g)\,.
\qqq
On overlaps $O_{ij}\x\xcG\ni(\si,g)$ of trivialisation charts,\ we find 
\qq\nn
\mu_j\bigl(\widetilde\xcP\t_i^{-1}(\si,g)\bigr)=\mu\bigl(\widetilde\xcP\t_j^{-1}\bigl(\si,\b_{ji}(\si)\lact g\bigr)\bigr)\equiv s\bigl(\b_{ji}(\si)\lact g\bigr)=s(g)\equiv\mu_i\bigl(\widetilde\xcP\t_i^{-1}(\si,g)\bigr)\,.
\qqq
Hence,\ the $\mu_i$ are restrictions $\mu_i=\mu\rstr_{\pi_{\widetilde\xcP}^{-1}(O_i)}
$ of a globally smooth map
\qq\nn
\mu\colo\widetilde\xcP\too M\,.
\qqq

Next,\ consider smooth maps
\qq\nn
\varrho_i\colo\pi_{\widetilde\xcP}^{-1}(O_i){}_{\mu_i}\hspace{-3pt}\x_t\hspace{-1pt}\xcG\too\pi_{\widetilde\xcP}^{-1}(O_i),\ \bigl(\widetilde\xcP\t_i^{-1}(\si,g),h\bigr)\longmapsto\widetilde\xcP\t_i^{-1}(\si,g.h)\,.
\qqq
Again,\ the above glue to a globally smooth map 
\qq\nn
\varrho\colo\widetilde\xcP{}_{\mu}\hspace{-3pt}\x_t\hspace{-1pt}\xcG\too\widetilde\xcP
\qqq
since,\ for arbitrary $(\si,g)\in O_{ij}\x\xcG$ and $h\in t^{-1}(\{s(g)\})$,
\qq\nn
\varrho_j\bigl(\widetilde\xcP\t_i^{-1}(\si,g),h\bigr)&=&\varrho_j\bigl(\widetilde\xcP\t_j^{-1}\bigl(\si,\b_{ji}(\si)\lact g\bigr),h\bigr)\equiv\widetilde\xcP\t_j^{-1}\bigl(\si,\bigl(\b_{ji}(\si)\lact g\bigr).h\bigr)=\widetilde\xcP\t_j^{-1}\bigl(\si,\b_{ji}(\si)\lact(g.h)\bigr)\cr\cr
&=&\widetilde\xcP\t_i^{-1}\bigl(\si,g.h\bigr)\equiv\varrho_i\bigl(\widetilde\xcP\t_i^{-1}(\si,g),h\bigr)\,,
\qqq
where in the third equality we used \Reqref{eq:BisAct-vs-str-v}.\ The triple $(\widetilde\xcP,\mu,\varrho)$ satisfies axioms (GrM1)--(GrM3) of a right $\xcG$-module as it is locally modelled on the canonical right-$\xcG$-module structure of \Reqref{eq:can-R-Gr-mod},\ in conformity with Def.\,\ref{def:principaloid-G-bndl}.
\eroof

\brem Given the equivalence between the notions of principaloid bundle and principaloid $\xcG$-bundle,\ we shall use them interchangeably,\ choosing the latter,\ in particular,\ whenever we want to emphasise the underlying structure groupoid.
\erem 

Having established the equivalence of the first two generalisations of the notion of principal $\txG$-bundle to Lie-groupoidal symmetries,\ we may,\ next,\ investigate a relation between principaloid bundles and $\xcG$-fibred $\bB$-bundles.\ To this end,\ observe that a $\xcG$-fibred $\bB$-bundle is a fibre bundle $(\widehat\xcP,\Si,\xcG,\pi_{\widehat\xcP})$ on whose total space $\widehat\xcP$ there is \emph{given} an action 
\qq\label{eq:hatR-act-hatP}
\widehat\xcR\colo \bB\too{\rm Aut}_{{\rm {\bf Bun}}(\Si)}(\widehat\xcP)_{\rm vert}
\qqq
such that for \emph{every} choice of local trivialisations
\qq\label{eq:hatP-loc-triv}
\widehat\xcP\t_i\colo \pi_{\widehat\xcP}^{-1}(O_i)\xrightarrow{\ \cong\ } O_i\x\xcG\,,\qquad i\in I
\qqq
we have $\bB$-equivariance
\qq\nn
\widehat\xcR{}_\b\circ\widehat\xcP\t_i^{-1}(\si,g)=\widehat\xcP\t_i^{-1}\bigl(\si,R_\b(g)\bigr)\,,\qquad(\si,g)\in O_i\x\xcG\,.
\qqq
We then obtain,\ for arbitrary $\si'\in O_i\cap O_j\equiv O_{ij}$ and $ g\in\xcG$,
\qq\nn
\widehat\xcP\t_i\circ\widehat\xcP\t_j^{-1}(\si', g)=\bigl(\si',t_{ij}(\si', g)\bigr)\,,
\qqq
where the $t_{ij}\colo O_{ij}\to\Diff_{\bB}(\xcG),\ \si\mapsto t_{ij}(\si,\cdot)$ are smooth maps with values in
\qq\label{eq:DiffBG}
\Diff_{\bB}(\xcG)=\bigl\{\ \Psi\in\Diff(\xcG) \:\vert\;\forall \b\in\bB \colo \:\Psi\circ R_{\b}= R_{\b}\circ\Psi \ \bigr\}\,,
\qqq
the {\bf group of $ R(\bB)$-equivariant diffeomorphisms of $\xcG$}.\ Thus,\ whenever,\ for some choice of local trivialisations,\ the transition maps $t_{ij}$ take the special form:
\qq\label{eq:tij-Lbij}
t_{ij}(\si)=L_{\b_{ij}(\si)}\,,
\qqq
expressed in terms of a \Cv ech 1-cocycle $\b_{ij}\colo O_{ij}\to\bB,\ (i,j)\in I^{\x 2}_\cO$,\ the $\xcG$-fibred $\bB$-bundle $\widehat\xcP$ is a principaloid bundle.\ This,\ however,\ is not true in general,\ and so \emph{not} every $\xcG$-fibred $\bB$-bundle is a principaloid bundle.\ This is readily seen by means of the following example\footnote{The authors thank  Damian Kayzer and Jakub Filipek for suggesting and working out (a special case of) this example.}.
\beg\label{eg:SWKF}
Consider,\ after \Rxcite{Remark 2.18\,b)}{Schmeding:2016},\ the pair groupoid of a manifold $M = M_1\sqcup M_2$ defined by two non-diffeomorphic connected manifolds $M_1$ and $M_2$ of the same dimension.\ Since $\bB\cong\Diff(M)$,\ there is clearly no global bisection through a point $(m_2,m_1)\in M_2\x M_1\subset M\x M$ as there is no (global) diffeomorphism which could map $m_1$ to $m_2$.\ We have the following decomposition of the arrow manifold of ${\rm Pair}(M)$ into disjoint connected components
\qq\nn
M\x M=M_1\x M_1\sqcup M_1\x M_2\sqcup M_2\x M_1\sqcup M_2\x M_2\,.
\qqq
Take an arbitrary diffeomorphism of $M \x M $ with the corresponding decomposition
\qq\label{eq:bad-diffeo}
\Phi=\id_{M_1\x M_1}\sqcup\id_{M_1\x M_2}\sqcup (f_2\x\id_{M_1})\sqcup\id_{M_2\x M_2}\,,
\qqq
written for some $f_2\in\Diff(M_2)\setminus\{\id_{M_2}\}$.\ As the diffeomorphism has a trivial restriction to the second argument in its domain,\ it manifestly commutes with 
\qq\nn
R(\bB)=\bigl\{\ \id_{M }\x(f_1\sqcup f_2) \ \vert\ (f_1,f_2)\in\Diff(M_1)\x\Diff(M_2)\ \bigr\}\subset\Diff(M )\,.
\qqq
On the other hand,\ we have
\qq\nn
L(\bB)=\bigl\{\ (f_1\sqcup f_2)\x\id_{M } \ \vert\ (f_1,f_2)\in\Diff(M_1)\x\Diff(M_2)\ \bigr\}\subset\Diff(M )\,,
\qqq
and so,\ clearly,\ $\Phi\notin L(\bB)$.
\eeg
\noindent Restricting the class of structure groupoids as below leads to an equivalence of the two notions,\ the pricipaloid ($\xcG$-)bundle and the $\xcG$-fibred $\bB$-bundle.
\bedef\label{def:Id-red}
A Lie groupoid $\xcG$ shall be called {\bf $\Id$-reducibile} if through every point $g\in\xcG$,\ there is a global bisection $\beta_g\in \bB$ such that 
\qq\label{eq:Id-red}
g=\beta_g(s(g))\equiv R_{\beta_g}\bigl(\Id_{t(g)}\bigr)\,.
\qqq
\exdef
\beg 
There is a large class of $\Id$-reducibile Lie groupoids including all source-connected ones,\ see \cite[Thm.\,3.1]{Zhong:2009},\ but evidently also every action groupoid. 
\eeg
\bethe\label{thm:Bequiv-as-princoid}
For an $\Id$-reducibile Lie groupoid $\xcG$,\ there is a one-to-one correspondence between principaloid ($\xcG$-)bundles and $\xcG$-fibred $\bB$-bundles.
\ethe

\beroof
That every principaloid bundle is a $\xcG$-fibred $\bB$-bundle is implied directly by the existence of the right action $\xcR$ of $\bB$,\ see \eqref{eq:R-act-Poid}.

Conversely,\ transition maps $t_{ij}$ of a $\xcG$-fibred $\bB$-bundle associated with arbitrary local trivialisations take values in $\Diff_\bB(\xcG)$,\ see \eqref{eq:DiffBG},\ and so they take the desired form \eqref{eq:tij-Lbij} in direct consequence of a proposition which we prove next.
\berop\label{prop:RB-comm-LB}
For every $\Id$-reducibile Lie groupoid $\xcG$,\ 
the commutant or centraliser $C_{\Diff(\xcG)}(R(\bB))\equiv R(\bB)'$ within $ \,\Diff(\xcG)$ of the subgroup $R(\bB)$,\ the group of right-multiplications by bisections,\ is $L(\bB)$,\ the group of left-multiplications by bisections:
\qq\label{eq:comm-RB-in-DiffG}
\Diff_{\bB}(\xcG)\equiv  R(\bB)'=L(\bB)\,.
\qqq
\eerop
\begin{proproof}
Consider a diffeomorphism $\Phi\in\Diff(\xcG)$ with the property that,\ for all $(g,\b)\in\xcG\x\bB$,\ one has 
\qq\label{eq:Rinv-diffG}
\Phi(g\ract\b)=\Phi(g)\ract\b\,,
\qqq
for every point $ g\in\xcG$.\ We shall demonstrate that this implies the identity
\qq\label{eq:sPhis}
s\circ\Phi=s\,.
\qqq
Intuitively,\ this identity is clear:\ If it were \emph{not} satisfied at,\ say,\ $g\in\xcG$,\ we could deform a given $\b$,\ for which \eqref{eq:Rinv-diffG} holds true at $g$,\ by replacing $\b$ with its product $\b\cdot\d$ with a bisection $\d\in\bB$ differing from $\Id$ only in a neighbourhood of $s(\Phi(g))$ disjoint\footnote{The existence of such a neighbourhood is guaranteed by the Hausdorff topology of $M$.} from that of $s(g)$.\ When plugged back into \Reqref{eq:Rinv-diffG},\ the deformed bisection would give us the same left-hand side,\ but a changed right-hand side,\ and so a contradiction.\ Below,\ we make this reasoning precise by passing to the tangent of $\bB$.

Making use of Prop.\,\ref{prop:BisAct-vs-str}(v),\ we may rewrite this identity as 
\qq\nn
\Phi\bigl(l_g\circ R_\b(\Id_{s(g)})\bigr)=\Phi\bigl(R_\b\circ l_g(\Id_{s(g)})\bigr)\equiv\Phi\bigl( g\ract\b\bigr)=\Phi\bigl( g\bigr)\ract\b\equiv R_\b\circ l_{\Phi(g)}\bigl(\Id_{s\circ\Phi(g)}\bigr)=l_{\Phi(g)}\circ R_\b\bigl(\Id_{s\circ\Phi(g)}\bigr)\,.
\qqq
This yields
\qq\label{eq:Rinv-diffG-cast}
\Id_{ s\circ\Phi( g)}\ract\b= l_{\Phi( g)^{-1}}\circ\Phi\circ l_{ g}\bigl(\Id_{ s( g)}\ract\b\bigr)\,.
\qqq
Invoking Thm.\,\ref{thm:LieBis} and taking $\b\equiv F^t_\xcX\circ\Id$ to be the value of the flow $F^\cdot_\xcX$ of a vector field $\xcX\in\G_{\rm c}(\ker\,\txT t\rstr_{\Id(M)})$,\ attained at a `time' $t\approx 0$ away from the initial condition $\Id$,\ 
we take the derivative of relation \eqref{eq:Rinv-diffG-cast} in the parameter $t$ at $t=0$,\ whereby we obtain
\qq\nn
\xcX\bigl(\Id_{ s\circ\Phi( g)}\bigr)=  \txT_{\Id_{ s( g)}}\bigl( l{}_{\Phi( g)^{-1}}\circ\Phi\circ l{}_{ g}\bigr)\bigl(\xcX\bigl(\Id_{ s( g)}\bigr)\bigr)\,.
\qqq
Next,\ we apply the following argument to eliminate the possibility of $ s\circ\Phi( g)\neq s( g)$:\ Assume $\xcX$ to satisfy $\xcX(\Id_{ s\circ\Phi( g)})\neq 0$,\ and so also $\xcX(\Id_{ s( g)})\neq 0$.\ If the inequality were true,\ we could take separating neighbourhoods $ O_1\ni\Id_{ s( g)}$ and $ O_2\ni\Id_{ s\circ\Phi( g)}$ (with $O_1\cap O_2=\emptyset$) and replace the vector field $\xcX$ with $\xcX^\chi\equiv\chi\lact\xcX$ for some $\chi\in C^\infty(\Id(M),\bR)$ with ${\rm supp}\,\chi\subset O_2$ and $\chi(\Id_{ s\circ\Phi( g)})=1$,\ whereupon we would arrive at the contradiction
\qq\nn
0\neq\xcX\bigl(\Id_{ s\circ\Phi( g)}\bigr)\equiv\xcX^\chi\bigl(\Id_{ s\circ\Phi( g)}\bigr)=\txT_{\Id_{ s( g)}}\bigl( l{}_{\Phi( g)^{-1}}\circ\Phi\circ l{}_{ g}\bigr)\bigl(\xcX^\chi\bigl(\Id_{ s( g)}\bigr)\bigr)=0\,.
\qqq

With this result in hand,\ we may use the assumption of $\Id$-reducibity to find,\ for every $g\in\xcG$,
\qq\nn
\Phi(g)&=&\Phi\bigl(\Id_{t(g)}\ract \beta_g\bigr)= \Phi\bigl(\Id_{t(g)}\bigr)\ract \beta_g\equiv \Phi\bigl(\Id_{t(g)}\bigr).\beta_g\bigl((t_*\beta_g)^{-1}\bigl((s\circ\Phi)\bigl(\Id_{t(g)}\bigr)\bigr)\bigr)\cr\cr
&=&\Phi\bigl(\Id_{t(g)}\bigr).\beta_g\bigl((t_*\beta_g)^{-1}\bigl(s\bigl(\Id_{t(g)}\bigr)\bigr)\bigr)=\Phi\bigl(\Id_{t(g)}\bigr).\beta_g\bigl((t_*\beta_g)^{-1}\bigl(t(g)\bigr)\bigr)\,,
\qqq
and so,\ by Prop.\,\ref{prop:rasoR}(i),
\qq\nn
\Phi(g)=\Phi\bigl(\Id_{t(g)}\bigr).g\,,
\qqq
see \Reqref{eq:BgtBgt}.\ From this point onwards,\ the proof continues exactly as that of Prop.\,\ref{prop:rGequiv-LB}.
\end{proproof}

The above reasoning implies that \emph{every} trivialisation determines a reduction of the structure group of $\widehat\xcP$ as in Def.\,\ref{def:principoidle-cech}.
\eroof

For {\em non}-bisection-complete structure groupoids $\xcG$,\ instead,\ principaloid bundles form a subset of the set of $\xcG$-fibred $\bB$-bundles.\ In general,\ this subset is proper,\ which can be established directly through examination of $\bB$-actions on the total space of a (simple) preprincipaloid bundle induced by reductions from different $\bB$-equivariance classes,\ see Example \ref{eg:acts-from-reds} below.\medskip

We close the present section with a corollary of the reasoning employed in the last proof,\ which shall prove useful when we come to discuss connections on pricipaloid bundles.
\berop\label{prop:RinV-in-kerTs}
For every Lie groupoid $\xcG$,\ vector fields from the subspace
\qq\nn
\G(\txT\xcG)^{R(\bB)}:=\{\ \xcX\in\G(\txT\xcG) \quad\vert\quad \forall_{\b\in\bB}\colo R_{\b\,*}(\xcX)=\xcX \ \}\,,
\qqq
{\it i.e.},\ invariant under right-multiplication by bisections,\ are tangent to $s$-fibres,\
\qq\nn
\G(\txT\xcG)^{R(\bB)}\subset\G(\ker\,\txT s)\,.
\qqq
\eerop
\beroof
By a standard argument,\ the flow $\Phi_\xcX$ of an $R(\bB)$-invariant vector field $\xcX\in\G(\txT\xcG)^{R(\bB)}$ is $R(\bB)$-equivariant.\ Now,\ recall identity \eqref{eq:sPhis},\ which uses merely the $R(\bB)$-equivariance of $\Phi$,\ and does \emph{not} require that the latter be a global diffeomorphism.\ It applies to $\Phi_\xcX$ verbatim,\ yielding
\qq\nn
s\circ\Phi_\xcX=s\,.
\qqq
Thus,\ integral lines of $\xcX$ are contained in $s$-fibres,\ and so their tangents lie in $\ker\,\txT s$.\ However,\ the tangents are,\ by definition,\ just values attained by $\xcX$ along the corresponding integral lines,\ and the latter go through every point in $\xcG$,\ which implies the statement of the proposition. 
\eroof

\subsection{The induced principal groupoid bundle}

In the literature,\ there is yet another generalisation of principal $\txG$-bundles with a structure group $\txG$ replaced by a structure groupoid $\xcG$,\ which goes under the name of principal groupoid bundles \cite{MacKenzie:1987,Moerdijk:1991}.\ In this case,\ the main requirement is that there be an action of the groupoid $\xcG$ on the total space $\breP$ along a moment map $\mu \colo \breP \to M$.\ The action is required to be transitive and free along the fibres in $\breP$,\ which,\ in general,\ is not a fibre bundle,\ but a surjective submersion over the base $\Sigma$. (See Appendix \ref{app:princ-grpd-bndl} for details).

The notion of a principaloid bundle $\xcP$,\ introduced in the present paper,\ is different from this generalisation.\ However,\ every principaloid bundle $\xcP \to \Sigma$ is a principal $\xcG$-bundle in a canonical way,\ but over a different base,\ which itself is a fibre bundle over $\Sigma$.\ More specifically:
\bethe 
\label{thm:duck-as-prince}
Every principaloid bundle $\xcP$ canonically induces a fibre bundle $\xcF$ with model
\qq\label{eq:clutch-xcF}
\xcF\cong\bigsqcup_{i\in I}\,\bigl( O_i\x M\bigr)/\sim_{  t_*\b_{\cdot\cdot}}\,,
\qqq
written in terms of the transition 1-cocycle $\{ \b_{ij}\}_{(i,j)\in I^{\x 2}_\cO}$ of $\xcP$,\ which we realise by the shadow action of $\bB$ on $M$.\ It comes with a bundle map
\qq\label{diag:sitting-duck}
\alxydim{@C=.5cm@R=1.cm}{ \xcP \ar[rr]^{\xcD} \ar[dr]_{\pi_\xcP} & & \xcF \ar[dl]^{\pi_\xcF} \\ & \Si & }\quad,
\qqq
locally modelled on $  t\colo \xcG\to M$.
The quintuple $(\xcP,\xcF,\xcD,\mu,\varrho)$,\ in which $\mu\colo\xcP\to M$ is the moment map for the defining $\xcG$-action $\varrho\colo \xcP{}_{\mu}\hspace{-3pt}\x_t\hspace{-1pt}\xcG\to\xcP$ on $\xcP$,\ carries a canonical structure of a principal-$\xcG$-bundle object in the category of fibre bundles over $\Si$,\ in the sense captured by the commutative diagram
\qq\label{diag:duck-as-prince}
\alxydim{@C=.75cm@R=1.5cm}{ \Si \ar[d]_{\id_\Si} & & \xcP \ar[ll]_{\pi_\xcP} \ar@{->>}[d]_{\xcD} \ar[rd]^{\mu} & & \xcG \ar@{=>}[ld] \\ \Si & & \xcF \ar[ll]^{\pi_\xcF} & M & }\,.
\qqq
\ethe
\beroof
Formula \eqref{eq:clutch-xcF} implies that the induced bundle admits local trivialisations 
\qq\label{eq:F-loc-triv}
\xcF\t_i\colo \pi_\xcF^{-1}( O_i)\xrightarrow{\ \cong\ } O_i\x M
\qqq 
such that the corresponding transition mappings are
\qq\nn
\xcF\t_i\circ\xcF\t_j^{-1}\colo  O_{ij}\x M\too O_{ij}\x M,\ (\si,m)\longmapsto\bigl(\si,  t_*\bigl(\b_{ij}(\si)\bigr)(m)\bigr)\,.
\qqq
Now,\ the local mappings
\qq\nn
\xcD_i\colo \pi_\xcP^{-1}(O_i)\too\pi_\xcF^{-1}(O_i),\ \xcP\t_i^{-1}\bigl(\si, g\bigr)\longmapsto\xcF\t_i^{-1}\bigl(\si,  t(g)\bigr)
\qqq
glue up smoothly at $\si'\in O_{ij}$ as
\qq\nn
&&\xcD_j\bigl(\xcP\t_i^{-1}\bigl(\si', g\bigr)\bigr)=\xcD_j\bigl(\xcP\t_j^{-1}\bigl(\si', L_{\b_{ji}(\si')}\bigl( g\bigr)\bigr)\bigr)\equiv\xcF\t_j^{-1}\bigl(\si',  t\bigl(\b_{ji}(\si')\bigl(  t( g)\bigr).g\bigr)\bigr)\cr\cr
&=&\xcF\t_j^{-1}\bigl(\si',  t_*\b_{ji}(\si')\bigl(  t( g)\bigr)\bigr)=\xcF\t_i^{-1}\bigl(\si',  t_*\b_{ij}(\si')\circ  t_*\b_{ji}(\si')\bigl(  t( g)\bigr)\bigr)\cr\cr
&=&\xcF\t_i^{-1}\bigl(\si',  t_*\bigl(\b_{ij}(\si')\cdot\b_{ji}(\si')\bigr)\bigl(  t( g)\bigr)\bigr)=\xcF\t_i^{-1}\bigl(\si',  t( g)\bigr)\equiv\xcD_i\bigl(\xcP\t_i^{-1}\bigl(\si', g\bigr)\bigr)\,.
\qqq

Taking the above into account,\ alongside the structure of a right-$\xcG$-module space $(\xcP,\mu,\varrho)$ on $\xcP$ detailed in the second part of the proof of Thm.\,\ref{thm:principaloid-GB-equiv},\ we demonstrate the existence,\ on the fibred square $\xcP\x_\xcF\xcP\equiv\xcP\,{}_{\xcD}\hspace{-2pt}\x_{\xcD}\hspace{-1pt}\xcP$,\ of a unique map
\qq\label{eq:PhiP}
\phi_\xcP\colo\xcP\x_\xcF\xcP\too\xcG
\qqq
with the property $(\pr_1,\phi_\xcP)=(\pr_1,\varrho)^{-1}$.\ The latter property can be expressed more explicitly as
\qq\label{eq:div-maP-props}
t\circ\phi_\xcP(p_1,p_2)=\mu(p_1)\,,\qquad\qquad p_2=\varrho\bigl(p_1,\phi_\xcP(p_1,p_2)\bigr)\,,
\qqq
with arbitrary $(p_1,p_2)\in\xcP\x_\xcF\xcP$.\ Consider smooth maps
\qq\nn
\phi_i\colo\pi_\xcP^{-1}(O_i)\x_\xcF\pi_{O_i}^{-1}(O_i)\too\xcG,\ \bigl(\xcP\t_i^{-1}(\si,g_1),\xcP\t_i^{-1}(\si,g_2)\bigr)\longmapsto g_1^{-1}.g_2\,,
\qqq
whose well-definedness hinges on the identity $s(g_1^{-1})=t(g_1)=t(g_2)$,\ derived from
\qq\nn
\xcF\t_i^{-1}\bigl(\si,t(g_1)\bigr)\equiv\xcD\circ\xcP\t_i^{-1}(\si,g_1)=\xcD\circ\xcP\t_i^{-1}(\si,g_2)\equiv\xcF\t_i^{-1}\bigl(\si,t(g_2)\bigr)\,.
\qqq
These maps satisfy the equality,\ written for an arbitrary $\si\in O_{ij}$,
\qq\nn
&&\phi{}_j\bigl(\xcP\t_i^{-1}(\si,g_1),\xcP\t_i^{-1}(\si,g_2)\bigr)=\phi{}_j\bigl(\xcP\t_j^{-1}\bigl(\si,\b_{ji}(\si)\lact g_1\bigr),\xcP\t_j^{-1}\bigl(\si,\b_{ji}(\si)\lact g_2\bigr)\bigr)\cr\cr
&\equiv&\bigl(\b_{ji}(\si)\lact g_1\bigr)^{-1}.\bigl(\b_{ji}(\si)\lact g_2\bigr)=\bigl(g_1^{-1}\ract \b_{ij}(\si)\bigr).\bigl(\b_{ji}(\si)\lact g_2\bigr)=g_1^{-1}.\bigl(\b_{ij}(\si)\lact\bigl(\b_{ji}(\si)\lact g_2\bigr)\bigr)\cr\cr
&=&g_1^{-1}.\bigl(\bigl(\b_{ij}(\si)\cdot \b_{ji}(\si)\bigr)\lact g_2\bigr)=g_1^{-1}.g_2\equiv\phi{}_i\bigl(\xcP\t_i^{-1}(\si,g_1),\xcP\t_i^{-1}(\si,g_2)\bigr)\,,
\qqq
in whose derivation we have invoked Eqs.\,\eqref{eq:BisAct-vs-str-iv} and \eqref{eq:BisAct-vs-str-vi}.\ Hence,\ they glue up to a globally smooth map
\qq\nn
\phi_\xcP\colo \xcP\x_\xcF\xcP\too\xcG\,,\qquad\phi_\xcP\rstr_{\pi_\xcP^{-1}(O_i)\x_\xcF\pi_\xcP^{-1}(O_i)}\equiv\phi_i\,.
\qqq
The desired properties \eqref{eq:div-maP-props} are readily verified in the local presentation.
\eroof

\bedef 
We shall call $\xcF$ the {\bf shadow} ({\bf bundle}) of the principaloid bundle $\xcP$ and $\xcD \colo \xcP \to \xcF$ the {\bf sitting-duck map}.\footnote{Recall that a sitting duck is an easy target.} We call the map \eqref{eq:PhiP} the {\bf division map} of $\xcP$. 
\exdef

\brem
\label{rem:quotient}
In the light of Thm.\,\ref{thm:Godement},\ the above theorem implies that the shadow bundle $\xcF$ is diffeomorphic to $\xcP/\xcG$ and the sitting-duck map is the quotient map $\xcP\to\xcP/\xcG$.\ Note,\ however,\ that in contrast to ordinary Lie-group actions,\ freeness and properness of a $\xcG$-action does not,\ in general,\ give a principaloid bundle.
\erem
\brem
The division map $\phi_\xcP$ of $\xcP$ defined in the proof above is a straightforward generalisation of the standard division map $\phi_\txP\colo\txP\x_\Si\txP\to\txG$ of a principal $\txG$-bundle.\ Indeed,\ if we take into account Example \ref{eg:princ-as-princ},\ we obtain $\xcP\x_\xcF\xcP\equiv\txP\x_\Si\txP$ for $\xcG=\txG$,\ and the mapping is seen to assign to any pair of points in the $t$-fibre over a given point in the base (in a local trivialisation) the unique arrow that connects one of them to the other. 
\erem

We close the present section,\ dedicated to a more detailed study of the defining action of the structure groupoid on a principaloid bundle,\ with a reinterpretation of the $\bB$-action \eqref{eq:R-act-Poid} as an induced one.\ The induction is a consequence of the following general result,\ which shall become of particular relevance in our discussion of automorphisms of principaloid bundles.
\berop\label{prop:G-in-B_act}
Let $(X,\mu,\varrho)$ be a right $\xcG$-module.\ The action $\varrho\colo X{}_{\mu}\hspace{-3pt}\x_{t}\hspace{-1pt}\xcG\to X$ canonically induces a right action $\varrho_\bB\colo X\x\bB\too X$ of $\bB$ on $X$,\ as determined by the following commutative diagram:
\qq\nn
\alxydim{@C=3.cm@R=1.5cm}{ X{}_{\mu}\hspace{-3pt}\x_{\id_M}\hspace{-1pt}M\x\bB \ar[r]^{\quad\id_X\x\Inv\circ\ev} & X{}_{\mu}\hspace{-3pt}\x_t\hspace{-1pt}\xcG \ar[d]^{\varrho} \\ X\x\bB \ar[u]^{(\id_X,\mu)\x\Inv} \ar[r]_{\varrho_\bB} & X }\,,
\qqq
where $\ev\colo M\x \bB\too\xcG,\ (m,\b)\longmapsto\b(m)$ is the {\bf evaluation map}.\ An analogous statement holds true for left $\xcG$-modules.
\eerop
\beroof
The diagram yields the map
\qq\nn
\varrho_\bB\colo X\x\bB\too X\colo (x,\b)\longmapsto x\mact\b^{-1}\bigl(\mu(x)\bigr)^{-1}\,.
\qqq
Here,\ $\b^{-1}$ is the inverse of $\b$ in the group $\bB$,\ whereas the inverse applied after evaluation of this $\b^{-1}$ on $\mu(x)$ is the one in $\xcG$.\ The well-definedness of the map follows from the identity
\qq\nn
t\bigl(\b^{-1}\bigl(\mu(x)\bigr)^{-1}\bigr)\equiv t\bigl(\b\bigl((t_*\b)^{-1}\bigl(\mu(x)\bigr)\bigr)\bigr)=\mu(x)\,.
\qqq
It now suffices to check the axioms of a group action,\ invoking those of the groupoid module along the way.\ We have,\ by (GrM2) from Def.\,\ref{def:gr-mod} in Appendix \ref{app:princ-grpd-bndl},\ for every $x\in X$,
\qq\nn
\varrho_\bB(x,\Id)\equiv x\mact\bigl(\Id_{\mu(x)}^{-1}\bigr)^{-1}=x\mact\Id_{\mu(x)}=x \,.
\qqq
Similarly,\ by (GrM1),\ \eqref{eq:BisGr-act-obj},\ and (GrM3),\ for every $\b_1,\b_2\in\bB$,
\qq\nn
&&\varrho_\bB\bigl(\varrho_\bB(x,\b_1),\b_2\bigr)\equiv\bigl(x\mact\b_1^{-1}\bigl(\mu(x)\bigr)^{-1}\bigr)\mact\b_2^{-1}\bigl(\mu\bigl(x\mact\b_1^{-1}\bigl(\mu(x)\bigr)^{-1}\bigr)\bigr)^{-1}\cr\cr
&=&\bigl(x\mact\b_1^{-1}\bigl(\mu(x)\bigr)^{-1}\bigr)\mact\b_2^{-1}\bigl(s\bigl(\b_1^{-1}\bigl(\mu(x)\bigr)^{-1}\bigr)\bigr)^{-1}=\bigl(x\mact\b_1^{-1}\bigl(\mu(x)\bigr)^{-1}\bigr)\mact\b_2^{-1}\bigl(\bigl(t_*\b_1^{-1}\bigr)\bigl(\mu(x)\bigr)\bigr)\bigr)^{-1}\cr\cr
&=&x\mact\bigl(\b_2^{-1}\bigl(\bigl(t_*\b_1^{-1}\bigr)\bigl(\mu(x)\bigr)\bigr).\b_1^{-1}\bigl(\mu(x)\bigr)\bigr)^{-1}\equiv x\mact\bigl(\b_2^{-1}\cdot\b_1^{-1}\bigr)\bigl(\mu(x)\bigr)^{-1}=x\mact(\b_1\cdot\b_2)^{-1}\bigl(\mu(x)\bigr)^{-1}\cr\cr
&\equiv&\varrho_\bB\bigl(x,\b_1\cdot\b_2\bigr)\,.
\qqq
\eroof
\brem\label{rem:Bisec-geom-act}
One may also view the statement of Prop.\,\ref{prop:G-in-B_act} in a more geometric manner:\ The induced action $\varrho_\bB$ can be realised through restriction of $\varrho$ to submanifolds in $\xcG$ obtained by evaluation of bisections on $M$.\ Specialising this to the principaloid bundle $\xcP$,\ we obtain the following diagrammatic representation of the induced action:
\qq\nn
\alxydim{@C=.75cm@R=1.5cm}{ \Si \ar@{=}[d] & & \xcP \ar[ll]_{\pi_\xcP} \ar@{->>}[d]_{\xcD} \ar[rd]^{\mu} & \b(M) \ar@{^{(}->}[r] & \xcG \ar@{=>}[ld] \\ \Si & & \xcF \ar[ll]^{\pi_\xcF} & M \ar[u]^{\b} & }\,.
\qqq
\erem

For later purposes,\ it will be useful to,\ {\it vice versa},\ implement the $\xcG$-action by the $\bB$-action on $\xcP$.\ This is not always possible,\ since,\ in general,\ not every $g\in\xcG$ admits a global bisection passing through it.\ Thus,\ we stumble upon the condition of $\Id$-reducibility once again.
\becor\label{cor:rho-vs-R}
For every principaloid bundle $\xcP$ with a bisection-complete structure groupoid,\ the following identity holds true for arbitrary $(p,g)\in\xcP{}_{\mu}\hspace{-3pt}\x_t\hspace{-1pt}\xcG$
\qq\nn
\varrho(p,g)=p\ract \beta^{-1}_{g^{-1}}\,,
\qqq
where $\beta_g\in\bB$ is any bisection through $g$,\ as specified in Def.\,\ref{def:Id-red}.
\ecor

\subsection{Examples}

\beg\label{eg:princ-as-princ}
If $\xcG=\txG$,\ a Lie group,\ then $\bB \cong \txG$ and a principaloid bundle reduces to a standard principal $\txG$-bundle.\ Similarly,\ as $\xcF \cong \Si$,\ so does the induced principal groupoid bundle \eqref{diag:duck-as-prince}.\ Hence,\ in this special case,\ the two notions coincide. 
\eeg
\beg\label{eg:princ-as-gen}
If $\xcG={\rm Pair}(M)$,\ then $\bB\cong\Diff(M)$ and a principaloid bundle $\xcP\cong\xcF\x M$,\ where $\xcF$ is a generic fibre bundle with typical fibre $M$ and no restriction on the structure group.\ Here,\ the sitting-duck map $\xcD =\pr_1$ and the moment map $\mu = \pr_2$.\ The groupoid ${\rm Pair}(M)$ acts just on the second factor $M$ and does not affect $\xcF$.\ So,\ in this case,\ we recover general fibre bundles $\xcF$.\ 
\eeg

Principaloid bundles are particular fibre bundles,\ where the fibre is $\xcG$ and the structure group is reduced to $\bB = {\rm Bisec}(\xcG)$.\ At the same time,\ and as the previous example shows,\ general fibre bundles $\xcF$,\ with fibre $M$ and no restriction on the structure group $\Diff(M)$,\ can be described as particular principaloid bundles.\ This may be compared to the relation between Poisson manifolds $(M,\Pi$) and Lie algebroids $E\to M$:\ While the former induce particular Lie algebroids on $E=\txT^*M$,\ the latter are particular Poisson manifolds,\ where $E^*$ carries a canonical fibre-linear Poisson structure.

\beg\label{eg:duck-for-action}
Let $\xcG=\txG\x M$,\ the action groupoid associated with an action $\lambda$ of a Lie group $\txG$ on a manifold $M$.\ In this case,\ there exists an embedding of the group $G$ into the group $\bB$ of bisections of $\xcG$,\ given by $\iota \colo \txG\ \to \bB, \, g \mapsto \beta_g$,\ where $\beta_g(m)=(g,m)$.\ Let $\bB_0 := \iota(G)$ and $\xcP_0$ be a principaloid bundle in which the structure group is restricted/reduced to $\bB_0 \subset \bB$.

Then,\ $\xcP_0\cong\sfP\x M$,\  where $\sfP$ is a principal $\txG$-bundle with the transition mappings $\b_{ij} \colo O_{ij}\to\txG$ induced from the transition 1-cocycle with values in $\bB_0\cong \txG$.\ The induced shadow bundle $\xcF_0\cong\sfP\x_\la M$ is just the associated bundle obtained from $P$ by the action $\lambda$,\ with quotient map $\xcD \colo \sfP\x M \to \sfP\x_\la M$.\ Note that here still \emph{all} of $\bB$ acts on $\xcP_0=\sfP\x M$ and this action does not factorise into separate actions on $\sfP$ and $M$,\ except for the subgroup $\bB_0$.\ The moment map $\mu_0 = \pr_2$ intertwines the $\bB$-action on $\xcP_0$ and $M$. 

So, in the above fashion, we recover standard principal $G$-bundles $P\to \Sigma$ together with associated ones in the form of couples $(\xcP_0,\xcF_0)$. 
\eeg

\beg
For a symplectic groupoid of Example \ref{eg:symplgrpd},\ with a symplectic arrow manifold $(\xcG,\om)$ and a Poisson object manifold $(M,\Pi)$,\ it is customary to \emph{restrict} the choice of bisections to those $\b\in{\rm Bisec}(\xcG)$---termed {\bf lagrangean}---for which the corresponding submanifolds $\b(M)\subset\xcG$ are lagrangean in $(\xcG,\om)$,\ see \cite[Chap.\,II]{Coste:1987}.\ These form a Fr\'echet--Lie subgroup ${\rm Bisec}_{\rm lagr}(\xcG,\om)$ within ${\rm Bisec}(\xcG)$,\ which is regular in the sense of Kriegl and Michor,\ see \cite{Rybicki:2001}.\ The corresponding left- and right-multiplications of $(\xcG,\om)$ by ${\rm Bisec}_{\rm lagr}(\xcG,\om)\subset{\rm Bisec}(\xcG)$ are symplectic,\ see \cite[Prop.\,1.2,\ Chap.\,II]{Coste:1987}, and project to Poisson actions on $(M,\Pi)$ along both:\ the source and target maps,\ \cite[p.\,22]{Coste:1987}.\ We shall refer to principaloid bundles with the structure group reduced to ${\rm Bisec}_{\rm lagr}(\xcG,\om)$ as {\bf symplectic principaloid $(\xcG,\om)$-bundles} in future studies.\ In some cases,\ it may be natural to further restrict the structure group,\ so as to obtain,\ {\it e.g.},\ hamiltonian $\bB$-actions on $(\xcG,\omega)$ or on $(M,\Pi)$.
\eeg   

We conclude the present section by studying an example of two manifestly inequivalent $\xcG$-structures induced on a preprincipaloid bundle with a bisection-complete structure groupoid by representatives of two different equivariance classes of reductions.\ By Prop.\,\ref{prop:G-in-B_act},\ these give rise to two different $\bB$-structures on that bundle.
\beg\label{eg:acts-from-reds}
Consider a trivial preprincipaloid bundle $\widetilde\xcP=\Si\x\xcG$ over base $\Si$ with typical fibre $\xcG={\rm SO}(3)\x\bR^3$ given by the arrow manifold of the action groupoid $\grpd{{\rm SO}(3)\x\bR^3}{\bR^3}$ associated with the defining representation of ${\rm SO}(3)$.\ This trivial bundle admits a family of trivialisations
\qq\nn
\t_{\vec d}\colo\widetilde\xcP\too\Si\x\xcG,\ \bigl(\si,(R,\vec r)\bigr)\longmapsto\bigl(\si,(R,\vec r+\vec d)\bigr)
\qqq
indexed by vectors $\vec d\in\bR^3$.\ Each of these trivialisations comes with a $\xcG$-action---and so also a $\bB$-action---on $\widetilde\xcP$ of the following form:\ For every $(\si,(R,\vec r))\in\widetilde\xcP$ and $M\in{\rm SO}(3)$,
\qq\nn
&&\varrho_{\vec d}\bigl(\bigl(\si,(R,\vec r)\bigr),\bigl(M,M^{-1}\bigl(\vec r+\vec d\bigr)\bigr)\bigr):=\t_{\vec d}^{-1}\circ(\id_\Si\x r_{(M,M^{-1}(\vec r+\vec d))})\circ\t_{\vec d}\bigl(\si,(R,\vec r)\bigr)\cr\cr
&=&\t_{\vec d}^{-1}\bigl(\si,\bigl(R,\vec r+\vec d\bigr).\bigl(M,M^{-1}\bigl(\vec r+\vec d\bigr)\bigr)\bigr)=\t_{\vec d}^{-1}\bigl(\si,\bigl(R\cdot M,M^{-1}\bigl(\vec r+\vec d\bigr)\bigr)\bigr)=\bigl(\si,\bigl(R\cdot M,M^{-1}\bigl(\vec r+\vec d\bigr)-\vec d\bigr)\bigr)\,.
\qqq
The $\xcG$-orbits of $(\si,(\bd1,\vec 0))\in\widetilde\xcP$ under $\varrho_{\vec d}\,$ take the form: 
\qq\nn
\varrho_{\vec d}\bigl(\bigl(\si,(\bd1,\vec 0)\bigr),\xcG\bigr)=\bigl\{\ \bigl(\si,\bigl(M,M^{-1}\vec d-\vec d\bigr)\bigr) \ \vert\ M\in{\rm SO}(3) \ \bigr\}\,.
\qqq
Thus,\ for any two distinct choices of $\vec d$,\ we obtain different $\xcG$-orbits.\ Consider,\ in particular,\ the case $\vec d=0$.

Each such choice gives rise to a different structure of a principaloid bundle $\xcP_{\vec d}\,$ on $\widetilde\xcP$.\ Indeed,\ for $\xcP_{\vec d{}_1}$ to be equivalent to $\xcP_{\vec d{}_2}\,$ for arbitrary $\vec d{}_2\neq\vec d{}_1$,\ we would need to be able to have a common refinement of the respective trivialisations such that the transition maps between them would lie in $L(\bB)$.\ In the case at hand,\ bisections are of the form
\qq\nn
\b\colo \bR^3\too{\rm SO}(3)\x\bR^3,\ \vec r\longmapsto\bigl(\widetilde M\bigl(\vec r\bigr),\vec r\bigr)
\qqq
for some $\widetilde M\in C^\infty(\bR^3,{\rm SO}(3))$.\ A simple calculation yields the following expression for the left-multiplication of $(R,\vec r)\in\xcG\cong\widetilde\xcP{}_\si$ with $\si\in\Si$ by the above $\b\in\bB$:
\qq\label{eq:LB-on-R3}
L_\b\bigl(R,\vec r\bigr)=\bigl(\widetilde M\bigl(R\vec r\bigr)\cdot R,\vec r\bigr)\,.
\qqq
On the other hand,\ the transition map between $\t_{\vec d{}_1}$ and $\t_{\vec d{}_2}\,$  is given by the formula
\qq\nn
\bigl(\t_{\vec d{}_1}\circ\t_{\vec d{}_2}^{-1}\bigr)\bigl(\si,(R,\vec r)\bigr)=\bigl(\si,(R,\vec r+\vec d{}_1-\vec d{}_2)\bigr)\neq\bigl(\si,(R,\vec r\bigr)\bigr)\,,
\qqq
which is manifestly not of the form \eqref{eq:LB-on-R3}.
\eeg

\section{Connections}\label{sec:connections}

\subsection{Connections on principaloid bundles} 
\label{sub:connections1}
\bedef\label{def:principaloid-conn}
A (compatible) {\bf connection} on a principaloid bundle $\xcP$ is an Ehresmann connection $\txT\xcP=\txV\xcP\oplus\txH\xcP$ with a $\xcG$-invariant horizontal distribution $\txH\xcP\subset\ker\,\txT\mu$.\ In other words,\ for every $g\in\xcG$ and $p\in\mu^{-1}(\{t(g)\})$,\ the horizontal distribution satisfies
\qq \label{eq:conn} 
\txT\varrho_g (\txH_p\xcP)=\txH_{\varrho_g(p)}\xcP\,,
\qqq
where $\varrho_g(p)\equiv\varrho(p,g)\equiv p\mact g$.
\exdef 
\brem
The restriction $\txH\xcP\subset\ker\,\txT\mu$ is a prerequisite for the imposition of the condition of $\xcG$-invariance as the $\xcG$-action is $\mu$-fibred:\ An element $g\in\xcG$ acts only on $\mu^{-1}(\{t(g))\}$.\ Furthermore,\ as the momentum $\mu$ is modelled on $s$ in a local trivialisation $\xcP\t_i$,\ the kernel $\ker\,\txT\mu$ becomes $\pr_1^*\txT O_i\oplus\pr_2^*\ker\,\txT s$ in that trivialisation,\ which makes the restriction meaningful from the point of view of the definition of a horizontal distribution.  
\erem
\noindent As usual,\ we may replace the notion of a compatible (Ehresmann) connection with that of a smooth distribution of projectors over $\xcP$.
\bedef\label{def:Ginv-principal-conn} 
A (compatible) {\bf connection 1-form} on $\xcP$ is a $\txV\xcP$-valued 1-form $\Theta\in\Om^1(\xcP,\txV\xcP)\equiv\G(\txT^*\xcP\ox\txV\xcP)$
such that,\ when viewed as an element in $\End{(\txT\xcP)}$:
\bit
\item[(GC1)] $\Theta\rstr_{\txV\xcP}=\id_{\txV\xcP}$\, (a projection onto $\txV\xcP$);
\item[(GC2)] $\txT\mu\circ\Theta=\txT\mu$\, (kernel inclusion);
\item[(GC3)] $\forall\ g\in\xcG\,,\ p\in\mu^{-1}(\{t(g)\})\colo \ \txT_p\varrho_g(\ker\,\Theta_p)=\ker\,\Theta_{\varrho_g(p)}$\, ($\xcG$-invariance).
\eit
\exdef
\noindent The equivalence of the two definitions is established in the following proposition.
\berop\label{prop:conn-1-Econn}
Connections on a principaloid bundle $\xcP$ are in one-to-one correspondence with connection 1-forms on $\xcP$.
\eerop
\beroof
The choice of a decomposition $\txT\xcP=\txV\xcP\oplus\txH\xcP$ determines a 1-form $\Theta$ by declaring the latter form to vanish on the horizontal 
subbundle and to be the identity on the vertical one,\ which in particular implies (GC1) by construction.\ Conversely,\ given  $\Theta$,\ property (GC1) determines a decomposition with $\txH\xcP:=\ker\,\Theta$.\ With this identification,\ condition (GC2) becomes equivalent to the relation $\txH\xcP\subset\ker\,\txT\mu$,\ and (GC3) is a restatement of the $\xcG$-invariance of the horizontal distribution.
\eroof
\beg
If $\xcG=\txG$,\ then a connection on the principal $\txG$-bundle $\xcP\equiv\txP$ is a standard principal $\txG$-connection.\ In addition,\ since canonically $\txV\txP\cong\txP\x\ggt$,\ the corresponding connection 1-form $\Theta\in\G(\txT^*\txP\ox \txV \txP)=
\Om^1(\txP)\ox\ggt$ can be represented as a Lie algebra-valued 1-form in this case. 
\eeg
\beg
If $\xcG={\rm Pair}(M)$,\ then the principaloid bundle $\xcP$ can be identified with $\xcF\x_\Sigma \xcM$, where $\xcM \equiv \Sigma\x M$.\ A connection on $\xcP\cong\xcF\x_\Sigma \xcM$ corresponds uniquely to a pair composed of a connection on $\xcF$ and one on $\xcM$.\ As $\mu$ is the canonical projection to the factor $M$ in $\xcM$,\ the horizontal subbundle is constrained to lie within $\pr_1^*\txT\xcF\oplus\pr_2^*\widetilde\pr{}_1^*\txT\Si$,\ where $\pr_1\colo\xcF\x_\Si\xcM\to\xcF$,\  $\pr_2\colo\xcF\x_\Si\xcM\to\xcM$,\ and $\widetilde\pr{}_1\colo \Si\x M\to\Si$ denote the respective canonical projections.\ Note that the right action of $\xcG=M\x M$ is trivial on $\xcF$,\ and thus the connection there is unrestricted by the condition of $\xcG$-invariance.\ On $\xcM=\Sigma\x M$,\ $\xcG$ acts non-trivially only on the second factor $M$---by arbitrary replacements of points.\ Hence,\ the entire subbundle $\txH_{\mathrm{can}}\xcM := \widetilde\pr{}_1^*\txT\Sigma\oplus 0$ is preserved by the $\xcG$-action,\ and this corresponds to the canonical flat connection on $\xcM$. 

In the equivalent identification $\xcP \cong \xcF \x M$,\ one has $\txV\xcP\cong\pr_1^*\txV\xcF\oplus\pr_2^*\txT M$,\ and,\ according to what we established above,\ the connection 1-form splits as $\Theta=\pr_1^*\cA+\id_{\pr_2^*\txT M}$,\ with the trivial component on the factor $M$ and an arbitrary connection 1-form $\cA\in\Om^1(\xcF,\txV\xcF)$ on the factor $\xcF$.

We thus find,\ when choosing $\xcG$ to be the pair groupoid,\ that there is a one-to-one correspondence between principaloid bundles with connection $(\xcP,\txH\xcP)$ as in Def.\,\ref{def:principaloid-conn} and ordinary fibre bundles with Ehresmann connection $(\xcF,\txH\xcF)$.
\eeg

Since the shadow bundle $\xcF$ can be viewed as the quotient of $\xcP$ by the $\xcG$-action,\ see Remark \ref{rem:quotient},\ a connection on $\xcP$,\ being $\xcG$-invariant,\ descends to one on $\xcF$ along the quotient map $\xcD$.
\begin{propanition}
A connection on a principaloid bundle $\xcP$ canonically induces an Ehresmann connection on $\xcF$ along the sitting-duck map $\xcD$ as
\qq\label{eq:Whitney-shadow}
\txT\xcF=\txV\xcF\oplus\txH\xcF\,,\qquad\qquad\txV\xcF=\txT\xcD(\txV\xcP)\,,\qquad\txH\xcF=\txT\xcD(\txH\xcP)\,.
\qqq
We shall call it the {\bf shadow connection} on $\xcF$.
\end{propanition}
\beroof
Let $\,\txH\xcP\,$ be the horizontal distribution over $\xcP$.\ Using the sitting-duck map $\xcD\colo \xcP\to\xcF$,\ we define,\ over an arbitrary point $\xcD(p)\in\xcF$,\ a subspace
\qq\nn
\txH_{\xcD(p)}\xcF:=\txT_p\xcD\bigl(\txH_p\xcP\bigr)\subset\txT_{\xcD(p)}\xcF\,.
\qqq
We readily check that the right-hand side does not depend on the choice of the representative $p$ in the $\xcD$-fibre over $\xcD(p)$.\ Indeed,\ in virtue of Thm.\,\ref{thm:duck-as-prince},\ $\xcG$ acts transitively on the $\xcD$-fibre,\ which ensures that for any $q\in\xcD^{-1}(\{\xcD(p)\})$,\ we always find (a unique) $g\in t^{-1}(\{\mu(p)\})$ such that $q=\varrho_g(p)$,\ and so---by the assumption of $\xcG$-invariance of the horizontal distribution--- 
\qq\nn
\txT_q\xcD\bigl(\txH_q\xcP\bigr)\equiv\txT_{\varrho_g(p)}\xcD\bigl(\txH_{\varrho_g(p)}\xcP\bigr)=\bigl(\txT_{\varrho_g(p)}\xcD\circ\txT_p\varrho_g\bigr)(\txH_p\xcP)=\txT_p(\xcD\circ\varrho_g)(\txH_p\xcP)=\txT_p\xcD(\txH_p\xcP)\,.
\qqq
Here,\ the last equality is implied by the identity $\xcD\circ\varrho_g\rstr_{\mu^{-1}(\{t(g)\})}=\xcD\rstr_{\mu^{-1}(\{t(g)\})}$,\ true for every $g\in\xcG$,\ which,\ in turn,\ follows from the corresponding model identity $t\circ r_g\rstr_{s^{-1}(\{t(g)\})}=t\rstr_{s^{-1}(\{t(g)\})}$.

At this stage,\ it remains to verify that at every point $f=\xcD(p)\in\xcF$ as above,\ we have 
\qq\label{eq:FEhres-ind-PEhres}
\txT_f\xcF=\txV_f\xcF\oplus\txH_f\xcF\,.
\qqq
Being modelled on the surjective submersion $t$,\ both the sitting-duck map $\xcD \colon \xcP \to \xcF$ and its restriction $\xcD\rstr\colo \xcP_\sigma \to \xcF_\sigma$ to a fibre over every point $\sigma \in \Si$ are also surjective submersions.\ Consequently,  $\txT_p\xcD(\txT_p\xcP)=\txT_f\xcF$ and $\txT_p\xcD(\txV_p\xcP)= \txV_f\xcF$.\ Thus,\ by the Rank Theorem,\ \eqref{eq:FEhres-ind-PEhres} follows.\ Altogether,\ then,\ by the Constant-Rank Theorem,\ we obtain an induced Ehresmann connection
$
\txT\xcF=\txV\xcF\oplus\txH\xcF$
with the vector subbundle
\qq\label{eq:shadow-horizon}
\txH\xcF=\txT\xcD(\txH\xcP)\,.
\qqq
\eroof

\bedef
Let $\xcP$ be a principaloid bundle.\ Given a connection $\txT\xcP=\txV\xcP\oplus\txH\xcP$ on $\xcP$ and the corresponding shadow connection $\txT\xcF=\txV\xcF\oplus\txT\xcD(\txH\xcP)$ on $\xcF$,\ the {\bf shadow connection 1-form} on $\xcF$ is a $\txV\xcF$-valued 1-form $\Theta_\xcF\in\Om^1(\xcF,\txV\xcF)\equiv\G(\txT^*\xcF\ox\txV\xcF)$ with the property $\txT\xcD(\txH\xcP)=\ker\,\Theta_\xcF$.
\exdef

\berop\label{prop:shadconn-idef}
The shadow connection 1-form is uniquely determined by the identity
\qq\label{eq:shadconn-idef}
\Theta_\xcF\circ\txT\xcD=\txT\xcD\circ\Theta\,.
\qqq
\eerop
\beroof
The identity follows straightforwardly from \eqref{eq:Whitney-shadow}.\ Uniqueness of the 1-form thus defined is ensured by the surjectivity of $\txT\xcD$,\ itself a consequence of the surjective submersivity of $\xcD$.
\eroof

\beg A connection on $\xcP_0$ of Example \ref{eg:duck-for-action} is in a one-to-one correspondence with a connection on the principal $G$-bundle $\sfP$.\ The induced connection on $\xcF_0\cong\sfP\x_\la M$ is the same as the one obtained from the standard procedure for associated bundles,\ see,\ {\it e.g.},\ \cite{Kobayashi:1963}.
\eeg

\brem
In standard principal $\txG$-bundles $\sfP$,\ it is customary to describe connections by Lie algebra-valued connection 1-forms,\ $\omega \in \Omega^1(\sfP,\ggt)$.\ This is possible due to the freeness of the $\txG$-action on $\sfP$,\ leading to the global identification $\txV\xcP\cong\sfP\x\ggt$.\ For a general principaloid bundle,\ an analogous identification with the Lie algebroid $E=\mathrm{Lie}(\xcG)$ is not available:\ For $\dim\,M >0$,\ the fundamental distribution of the $\xcG$-action on $\xcP$ is a proper subbundle of $ \, \txV\xcP$,\ whose rank equals $\dim \,\xcG - \dim \, M$.

However,\ we recover $\Gamma(E)$-valuedness for certain local 1-forms on the base,\ see,\ in particular,\ formula \eqref{eq:GammaE} below.
\erem

\subsection{Local description}

There is also a definition of compatible connection 1-forms on principaloid bundles using local data:
\bedef\label{def:principoidle-conn-cech}
{\bf Local connection data} on a principaloid bundle $\xcP$ over $\Si$ are an assignment of a $\txT O_i$-foliated $\pr_2^*E$-valued 1-form $\txA_i$ on $O_i\x M$,\  $\txA_i\in\G(\pr_1^*\txT^*O_i\ox\pr_2^*E)$ to the trivialising cover $\cO\equiv\{O_i\}_{i\in I}$ of $\,\Si$.\ The 1-forms are subject to the following {\bf gluing law} over double intersections $O_{ij}\equiv O_i\cap O_j\ni\si$:
\qq\label{eq:glue-Aij}
\txA_i\bigl(\si,\b_{ij}(\si)\ulact m\bigr)=\txT_{\Id_m}C_{\b_{ij}(\si)}\circ\txA_j(\si,m)-\theta_{\rm R}\circ\txT_\si(\ev_m\circ \b_{ij})\,.
\qqq
Here,\ $\{ \b_{ij}\}_{(i,j)\in I^{\x 2}_\cO}$ is the transition 1-cocycle of $\xcP$ associated with $\cO$,\ and the pullback of the right-invariant Maurer--Cartan 1-form $\theta_{\rm R}$ on $\xcG$ is written in terms of the {\bf evaluation map} $\ev_m\colo\bB\to\xcG,\ \b\mapsto\b(m)$ at $m\in M$.
\exdef
\brem
 For standard principal $\txG$-bundles,\ with $\txG$ a matrix Lie group,\ the relation between two local connection 1-forms $A$ and $A'$ restricted to a common domain $O$,\   $A, A' \in \Omega^1(O,\ggt)$,\  takes the well-known form 
 \qq \label{eq:A'}
 A' = g \, A \, g^{-1} -  \, \txd g \, g^{-1} \, \equiv \, \Ad_g A - g^* \theta_{\rm R} \, , 
 \qqq 
 where $g \colo O \to \txG\,$ and $\theta_{\rm R}$ is the right-invariant Maurer-Cartan form on $\txG$.\ To recover a similar formula from \Reqref{eq:glue-Aij},\ let us define $A_\sigma :=\txA_j(\sigma,\cdot) \in \pr_1^*\txT_\sigma^* O \otimes_\R \Gamma \left((\pr_2^*E)_{(\sigma,\cdot)}\right) \cong  \txT^*_\sigma O \otimes_\R \Gamma (E)$,\  
 where $O:=O_{ij}$ and $\pr_1$ and $\pr_2$ refer to the factors in $O \times M$;\ and,\ likewise,\ $A'_\si := \txA_i(\sigma,\cdot) \in  \txT^*_\sigma O \otimes_\R \Gamma (E)$.\ Note that for $M$ compact,\ $\Gamma (E)$ can be identified with the Lie algebra of $\bB$ \cite{Schmeding:2019}.\ Then,\ for $\beta(\sigma) := \beta_{ij}(\si) \in \bB$,\ \Reqref{eq:glue-Aij} takes the form
 \qq \label{eq:Asigma}
A_\sigma' = \Ad_{\beta(\sigma)}A_\sigma - \txT_ {\beta(\sigma)} r_{\beta(\sigma)^{-1}} \circ \txT_\si \beta \, , 
 \qqq
where the last term can be symbolically rewritten as ``$\txd \beta \, \beta^{-1}$'',\ the pullback by $\b$ of the right-invariant (${\rm Lie}(\bB)$-valued) Maurer-Cartan form on the infinite dimensional Lie group $\bB$,\ evaluated at $\sigma$.\ In fact,\ if we contract the latter formula with a vector $v \in \txT_\sigma O$ and write 
\qq \label{eq:GammaE} A(v) := \langle A_\sigma , v \rangle \; \in \; \Gamma(E), 
\qqq 
then,\ both formulas,\ \eqref{eq:A'} and $\eqref{eq:Asigma}$,\ take the \emph{same} form 
\qq A'(v) = \Ad_{\beta(\si)} \: A(v)- \txT_ {\beta(\sigma)} r_{\beta(\sigma)^{-1}} \left( \txT_\si \beta(v) \right) \, ,
\qqq 
since for $\xcG =\txG$,\ one has $\beta(\si) \equiv g(\si)$ and $\Gamma(E) \cong \ggt$.\ Note that,\ for a general Lie groupoid,\ $A(v)$ is still to be evaluated at $m \in M$,\ whereupon $A(v)\vert_m \in E_m$.

While the above rewriting brings \Reqref{eq:glue-Aij} into a form more familiar from the standard setting of local connection 1-forms on principal $G$-bundles, the original version avoids all infinite dimensionality of $\bB$ and the Lie algebra of that Lie group altogether,\ and is thus preferable in our context. 
\erem
\brem
The appearance of the right-invariant Maurer-Cartan form on $\xcG$ in the gluing law \eqref{eq:glue-Aij} calls for some care.\ Since the differential form $\theta_{\rm R}$ is foliated,\ see Example \ref{eg:MC} in the appendices,\ it only takes arguments from the distribution $\ker\,\txT s\subset\txT\xcG$.\ Hence,\ it is necessary that the linear operator $\txT_\si(\ev_m\circ \b_{ij})$ map an arbitrary vector $v\in\txT_\si\Si$ into that distribution.\ That this is the case follows from the fact that the evaluation map $\ev_m$ fixes the source of all the arrows in $\{(\ev_m\circ\b_{ij})(\si)\}_{\si\in O_{ij}}$ to be $m$,\ so that,\ accordingly,\ $\txT_\si(\ev_m\circ \b_{ij})(v)\in(\ker\,\txT s)_{(\b_{ij}(\si))(m)}$.
\erem
\noindent Equivalence of Definitions \ref {def:Ginv-principal-conn} and \ref{def:principoidle-conn-cech} is stated in
\bethe\label{thm:loc-data-conn}
For every principaloid bundle $\xcP$,\ a connection 1-form determines local connection data,\ and {\it vice versa}.\ The relation between the two is given by the formula
\qq\nn
\bigl(\xcP\t_i^{-1\,*}\Theta\bigr)(\si,g)=\id_{\txT\xcG}\rstr_g+\txT_{\Id_{t(g)}}r_g\circ\txA_i\bigl(\si,t(g)\bigr)\,,\qquad (\si,g)\in O_i\x\xcG\,.
\qqq
\ethe
\beroof
Consider local trivialisations $\xcP\t_i\colo \pi_\xcP^{-1}(O_i)\xrightarrow{\ \cong\ }O_i\x\xcG$ over an open cover $\{ O_i\}_{i\in I}$ of the base $\Si$ of $\xcP$,\ and the associated transition 1-cocycle $ \{\b_{ij}\in C^\infty(O_{ij},\bB)\}$.\ Upon pullback\footnote{Recall that pullbacks along diffeomorphisms can be defined for arbitrary tensors,\ and not just differential forms,\ see,\ {\it e.g.},\ \cite{Thirring:1978}.\ For example,\ for a 1-1 tensor $\t$ and $\phi \in \mathrm{Diff}(M)$,\ one has: $(\phi^*\t)(v,\a)=\t(\txT_m\phi(v),\a\circ\txT_{\phi(m)}\phi^{-1})$ for every $v \in \txT_m M$ and $\alpha \in \txT_m^* M\equiv(\txT_m M)^*$.} of $\Theta$ along $\xcP\t_i$,\ we obtain $\xcP\t_i^{-1\,*}\Theta\in\G(\txT^*(O_i\x\xcG)\ox\pr_2^*\txT\xcG)$.\ Taking into account axiom (GC1),\ we can now write the pullback of $\Theta$ along $\xcP\t_i$ as
\qq\label{eq:loc-Theta-gen}
\bigl(\xcP\t_i^{-1\,*}\Theta\bigr)(\si,g)=\id_{\txT\xcG}\rstr_{\txT_g\xcG}+\G_i(\si,g)
\qqq
in terms of a vector-valued 1-form $\G_i\colo\pr_1^*\txT O_i\to\pr_2^*\txT\xcG$,\ where we use the identifications $\txT(O_i\x\xcG) \cong \pr_1^*\,\txT O_i\oplus \pr_2^*\,\txT\xcG$ and $\pr_2^*\,\txT\xcG\vert_{(\sigma,g)} \cong \txT_g\xcG$.\ Likewise,\ the notation for the two terms in the sum \eqref{eq:loc-Theta-gen} is adapted to the identification $\txT^*(O_i\x\xcG) \cong \pr_1^*\,\txT^*O_i\oplus \pr_2^*\,\txT^*\xcG$.\ Using local coordinates $\{\si^\mu\}^{\mu\in\ovl{1,\dim\,\Si}}$ on a neighbourhood of $\si\in O_i$,\ we write
\qq\nn
\G_i(\si,\cdot)=\txd\si^\mu(\si)\ox\G^\si_{i\,\mu}(\cdot)\,,
\qqq
where $\G^\si_{i\,\mu}\in\G(\txT\xcG)$. 
\belem\label{lem:smallth-in-kerTs}
The vector fields $\G^\si_{i\,\mu}$ lie in the distribution $\ker\,\txT s$.
\elem
\begin{lemproof}
In virtue of Prop.\,\ref{prop:Bequiv-conns-as-forms},\ the pullback \eqref{eq:loc-Theta-gen} determines the local presentation of the horizontal subspace over $O_i\x\xcG$ as the kernel of the 1-form $\id_{\txT\xcG}+\G_i$.\ This implies that the component of every horizontal vector field along the fibre is determined by the unconstrained component along the base and lies in the image of $\G_i$.\ But the field has to belong to $\ker\,\txT\mu$,\ which becomes $\pr_1^*\txT O_i\oplus\pr_2^*\ker\,\txT s$ in the local model as the momentum $\mu$ is modelled on the source map $s$.
\end{lemproof}

\belem
The local presentation \eqref{eq:loc-Theta-gen} can be written as
\qq\label{eq:princoid-conn-loc}
\bigl(\xcP\t_i^{-1\,*}\Theta\bigr)(\si,g)=\id_{\txT\xcG}\rstr_g+\txT_{\Id_{t(g)}}r_g\circ\txA_i\bigl(\si,t(g)\bigr)
\qqq
in terms of a locally smooth $E$-valued 1-form
\qq\label{eq:loc-gauge-potoid}
\txA_i(\si,m)\in\txT^*_\si O_i\ox_\bR E_m\,.
\qqq
\elem
\begin{lemproof}
In the light of the reasoning given in the proof of the previous lemma,\ the $\xcG$-invariance of the horizontal distribution translates into the same property of the $\G^\si_{i\,\mu}\in\G(\ker\,\txT s)$,\ that is $\G^\si_{i\,\mu}\in\G(\txT\xcG)_{\rm R}$,\ see Example \ref{eg:tanLiealgbrd}.\ Hence,\ for any $g\in\xcG$,
\qq\nn
\G_{i\,\mu}^\si(g)=\G_{i\,\mu}^\si\bigl(r_g\bigl(\Id_{t(g)}\bigr)\bigr)=T_{\Id_{t(g)}}r_g\bigl(\G_{i\,\mu}^\si\bigl(\Id_{t(g)}\bigr)\bigr)\,,
\qqq
and so
\qq\nn
\G_i(\si,g)=\txd\si^\mu(\si)\ox\G^\si_{i\,\mu}(g)=T_{\Id_{t(g)}}r_g\circ\bigl(\txd\si^\mu(\si)\ox\G_{i\,\mu}^\si\bigl(\Id_{t(g)}\bigr)\bigr)\,.
\qqq
Altogether,\ then,\ we obtain
\qq\nn
\txA_i(\si,m)=\txd\si^\mu(\si)\ox\G_{i\,\mu}^\si\bigl(\Id_m\bigr)\,.
\qqq
\end{lemproof}

We shall next investigate the gluing laws for the local 1-forms $\txA_i$ over the intersections $O_{ij}$.\ Clearly,\ these are determined by the behaviour of $\id_{\txT\xcG}$ --- viewed,\ as before,\ as an endomorphism of $\txT(O_{ij}\x\xcG)$ --- under the fibre identifications 
\qq\nn
\g{}_{ij}\colo O_{ij}\x\xcG\too O_{ij}\x\xcG,\ (\si,g,j)\longmapsto\bigl(\si,\check\g{}_{ij}(\si,g),i\bigr)\,,
\qqq
written in terms of the smooth maps
\qq\label{eq:gamcheckij}
\check\g{}_{ij}\colo O_{ij}\x\xcG\too\xcG,\ (\si,g)\longmapsto \b_{ij}(\si)\lact g\equiv r_g\bigl(\b_{ij}(\si)\bigl(t(g)\bigr)\bigr)\,.
\qqq
These are the identifications that underlie the model \eqref{eq:clutch-xcP} of the principaloid bundle $\xcP$,\ 
\qq\nn
\xcP\t_i^{-1}\circ\g_{ij}\rstr_{O_{ij}\x\xcG}=\xcP\t_j^{-1}\rstr_{O_{ij}\x\xcG}\,.
\qqq
Let us fix coordinates $\{\z^\a\}^{\a\in\ovl{1,\dim\,\xcG}}$ on a neighbourhood $\xcU$ of $ g $ (with $\p_\a\equiv\frac{\p\ }{\p\z^\a}$) and coordinates $\{\widetilde\z{}^\a\}^{\a\in\ovl{1,\dim\,\xcG}}$ on a neighbourhood $\widetilde\xcU$ of $\check\g{}_{ij}(\si,g)$ in $\xcG$ (with $\widetilde\p{}_\a\equiv\frac{\p\ }{\p\widetilde\z{}^\a}$),\ and compute,\ for $(v,V)\in\txT_\si O_{ij}\oplus\txT_g\xcG\cong\txT_{(\si,g)}(O_{ij}\x\xcG)$ arbitrary,\ how the local presentation \eqref{eq:princoid-conn-loc} pulls back along the identifications.\ We obtain
\qq\nn
&&\bigl(\bigl(\g_{ij}^*\xcP\t_i^{-1\,*}\Theta\bigr)(\si,g)\bigr)(v,V)=\cr\cr
&=&\txT_{\check\g{}_{ij}(\si,g)}L_{\b_{ij}(\si)^{-1}}\circ\bigl(\check\g{}_{ij}^*\txd\widetilde\z{}^\a(\si,g)\ox\widetilde\p{}_\a\rstr_{\check\g{}_{ij}(\si,g)}+\txT_{\Id_{t(\check\g{}_{ij}(\si,g))}}r_{\check\g{}_{ij}(\si,g)}\circ\txA_i\bigl(\si,t\bigl(\check\g{}_{ij}(\si,g)\bigr)\bigr)\bigr)(v,V)\cr\cr
&=&\txT_{\check\g{}_{ij}(\si,g)}L_{\b_{ij}(\si)^{-1}}\bigl(\txT_g L{}_{\b_{ij}(\si)}(V)+\txT_\si\check\g{}_{ij}(\cdot,g)(v)+\txT_{\Id_{t(\check\g{}_{ij}(\si,g))}}r_{\check\g{}_{ij}(\si,g)}\circ\txA_i\bigl(\si,t\bigl(\check\g{}_{ij}(\si,g)\bigr)\bigr)(v)\bigr)\cr\cr
&=&V+\txT_{\check\g{}_{ij}(\si,g)}L_{\b_{ij}(\si)^{-1}}\circ\bigl(\txT_\si\check\g{}_{ij}(\cdot,g)+\txT_{\Id_{t(\check\g{}_{ij}(\si,g))}}r_{\check\g{}_{ij}(\si,g)}\circ\txA_i\bigl(\si,t\bigl(\check\g{}_{ij}(\si,g)\bigr)\bigr)\bigr)(v)\,.
\qqq
For the right-hand side of \eqref{eq:princoid-conn-loc} to come from a globally smooth object on $\xcP$,\ the above must be equal to the result of evaluating $(\xcP\t_j^{-1\,*}\Theta)(\si,g)$ on the same pair $(v,V)$,\ 
\qq\nn
\bigl(\bigl(\xcP\t_j^{-1\,*}\Theta\bigr)(\si,g)\bigr)(v,V)=V+\txT_{\Id_{t(g)}}r_g\circ\txA_j\bigl(\si,t(g)\bigr)(v)\,,
\qqq
which yields the gluing law
\qq\nn
\txT_{\Id_{t(g)}}r_g\circ\txA_j\bigl(\si,t(g)\bigr)\must\txT_{\check\g{}_{ij}(\si,g)}L_{\b_{ij}(\si)^{-1}}\circ\bigl(\txT_\si\check\g{}_{ij}(\cdot,g)+\txT_{\Id_{t(\check\g{}_{ij}(\si,g))}}r_{\check\g{}_{ij}(\si,g)}\circ\txA_i\bigl(\si,t\bigl(\check\g{}_{ij}(\si,g)\bigr)\bigr)\bigr)\,.
\qqq
Now,\ differentiating formula \eqref{eq:gamcheckij} with respect to $\si$ (for $g$ fixed),\ one finds
\qq\nn
\txT_\si\check\g{}_{ij}(\cdot,g)\equiv\txT_\si\bigl(r_g\bigl(\b_{ij}(\cdot)\bigl(t(g)\bigr)\bigr)\bigr)=\txT_{\b_{ij}(\si)(t(g))}r_g\circ\txT_\si\bigl(\ev_{t(g)}\circ \b_{ij}\bigr)\,,
\qqq
with the evaluation map $\ev_{t(g)}$ smooth in the relevant sense,\ see (the proof of) \cite[Lemma 10.15]{Michor:1980}.\ Its presence in the operator $\txT_\si(\ev_{t(g)}\circ \b_{ij})$ ensures that the latter maps an arbitrary vector from $\txT_\si O_{ij}$ to the tangent $(\ker\,\txT s)_{\b_{ij}(\si)(t(g))}\subset\txT_{\b_{ij}(\si)(t(g))}\xcG$ of the $s$-fibre through $\b_{ij}(\si)(t(g))$ with $s(\b_{ij}(\si)(t(g)))=t(g)$.\ Hence,\ $\txT_{\b_{ij}(\si)(t(g))}r_{\b_{ij}(\si)(t(g))^{-1}}\circ\txT_\si\bigl(\ev_{t(g)}\circ \b_{ij}\bigr)$ maps $\txT_\si O_{ij}$ to $\txT_{\b_{ij}(\si)(t(g))}r_{\b_{ij}(\si)(t(g))^{-1}}((\ker\,\txT s)_{\b_{ij}(\si)(t(g))})=(\ker\,\txT s)_{\Id_{t(g)}}$,\ {\it i.e.},\ to the fibre of the Lie algebroid $E\equiv\Id^*\ker\,\txT s=M{}_{\Id}\hspace{-2pt}\x_{\pi_{\txT\xcG}}\hspace{-1pt}\ker\,\txT s$ over $t(g)$,\ in keeping with \cite[Prop.\,1.8.]{Schmeding:2019}.\ Taking into account the chain rule
\qq\label{eq:chain-r}
\txT_{\b_{ij}(\si)(t(g))}r_{\b_{ij}(\si)(t(g))^{-1}}=\txT_{\check\g{}_{ij}(\si,g)}r_{(\check\g{}_{ij}(\si,g))^{-1}}\circ\txT_{\b_{ij}(\si)(t(g))}r_g\,,
\qqq
we thus arrive at the reworked form of the gluing law
\qq\nn
&&\txT_{\Id_{t(g)}}r_g\circ\txA_j\bigl(\si,t(g)\bigr)\cr\cr
&\must&\txT_{\check\g{}_{ij}(\si,g)}L_{\b_{ij}(\si)^{-1}}\circ\txT_{\Id_{t(\check\g{}_{ij}(\si,g))}}r_{\check\g{}_{ij}(\si,g)}\circ\bigl(\txT_{\b_{ij}(\si)(t(g))}r_{\b_{ij}(\si)(t(g))^{-1}}\circ\txT_\si\bigl(\ev_{t(g)}\circ \b_{ij}\bigr)+\txA_i\bigl(\si,t\bigl(\check\g{}_{ij}(\si,g)\bigr)\bigr)\bigr)\,,
\qqq 
A straightforward calculation,\ using the inverse of \eqref{eq:chain-r} in the second equality,
\qq\nn
&&\txT_g r_{g^{-1}}\circ\txT_{\check\g{}_{ij}(\si,g)}L_{\b_{ij}(\si)^{-1}}\circ\txT_{\Id_{t(\check\g{}_{ij}(\si,g))}}r_{\check\g{}_{ij}(\si,g)}=\txT_{\b_{ij}(\si)(t(g))}L_{\b_{ij}(\si)^{-1}}\circ\txT_{\check\g{}_{ij}(\si,g)}r_{g^{-1}}\circ\txT_{\Id_{t(\check\g{}_{ij}(\si,g))}}r_{\check\g{}_{ij}(\si,g)}\cr\cr
&=&\txT_{\b_{ij}(\si)(t(g))}L_{\b_{ij}(\si)^{-1}}\circ\txT_{\Id_{t(\check\g{}_{ij}(\si,g))}}r_{\b_{ij}(\si)(t(g))}\,,
\qqq
subsequently enables us to rewrite the gluing law as
\qq\nn
&&\txA_j\bigl(\si,t(g)\bigr)\cr\cr
&\must&\txT_{\b_{ij}(\si)(t(g))}L_{\b_{ij}(\si)^{-1}}\circ\txT_{\Id_{t(\check\g{}_{ij}(\si,g))}}r_{\b_{ij}(\si)(t(g))}\circ\bigl(\txT_{\b_{ij}(\si)(t(g))}r_{\b_{ij}(\si)(t(g))^{-1}}\circ\txT_\si\bigl(\ev_{t(g)}\circ \b_{ij}\bigr)+\txA_i\bigl(\si,t\bigl(\check\g{}_{ij}(\si,g)\bigr)\bigr)\bigr)\,.
\qqq 
Finally,\ using identity \eqref{eq:Trahg-TRahb} in the form
\qq\nn
&&\txT_{\b_{ij}(\si)(t(g))}L_{\b_{ij}(\si)^{-1}}\circ\txT_{\Id_{t(\check\g{}_{ij}(\si,g))}}r_{\b_{ij}(\si)(t(g))}=\txT_{\b_{ij}(\si)(t(g))}L_{\b_{ij}(\si)^{-1}}\circ\txT_{\Id_{t_*\g{}_{ij}(\si)(t(g))}}R_{\b_{ij}(\si)}\rstr_{(\ker\,\txT s)_{\Id_{t(\check\g{}_{ij}(\si,g))}}}\cr\cr
&\equiv&\txT_{\Id_{t_*\g{}_{ij}(\si)(t(g))}}C_{\b_{ij}(\si)^{-1}}\rstr_{\ker\,\txT s_{\Id_{t(\check\g{}_{ij}(\si,g))}}}\,,
\qqq
taking into account the definition of the Maurer-Cartan form, see Example \ref{eg:MC},\ and replacing $t(g)$ by $m$,\ we ultimately arrive at the anticipated gluing law
\qq\nn
\txT_{\Id_m}C_{\b_{ij}(\si)}\circ\txA_j(\si,m)&=&\txA_i\bigl(\si,t_*\bigl(\b_{ij}(\si)\bigr)(m)\bigr)+\txT_{\b_{ij}(\si)(m)}r_{\b_{ij}(\si)(m)^{-1}}\circ\txT_\si\bigl(\ev_m\circ \b_{ij}\bigr)\cr\cr
&\equiv&\txA_i\bigl(\si,t_*\bigl(\b_{ij}(\si)\bigr)(m)\bigr)+\bigl(\ev_m\circ \b_{ij}\bigr)^*\theta_{\rm R}(\si)\,.
\qqq

For the converse,\ note that the gluing law \eqref{eq:glue-Aij} ensures descent of local 1-forms given by the right-hand side of \eqref{eq:princoid-conn-loc} to the model of $\xcP$ given in \eqref{eq:clutch-xcP}.
\eroof

\bedef
Let $\{\txA_i\}_{i\in I}$ be local data of a connection 1-form $\Theta$ on a principaloid bundle $\xcP$.\ The $\pr_1^*\txT O_i$-foliated $\ker\,\txT s$-valued 1-forms given,\ for any $(\si,g)\in O_i\x\xcG,\ i\in I$,\ by the formul\ae
\qq\nn
\G_i(\si,g)=\txT_{\Id_{t(g)}}r_g\circ\txA_i\bigl(\si,t(g)\bigr)
\qqq
are called the {\bf Christoffel forms} of $\Theta$.
\exdef

\berop\label{prop:gauge-coupl}
Let $\xcP$ be a principaloid bundle over $\Si$ with local trivialisations $\xcP\t_i\colo \pi_\xcP^{-1}(O_i)\xrightarrow{\ \cong\ }O_i\x\xcG$ over an open cover $\{O_i\}_{i\in I}$ of $\Si$.\ Let $\Theta\in\Om^1(\xcP,\txV\xcP)$ be a connection 1-form on $\xcP$ with local connection data $\{\txA_i\}_{i\in I}$.\ The shadow connection 1-form induced by $\Theta$ on the shadow $\xcF$ of $\xcP$ is the $\txV\xcF$-valued 1-form $\Theta_\xcF\in\Om^1(\xcF,\txV\xcF)$ with local presentations
\qq\nn
\bigl(\xcF\t_i^{-1\,*}\Theta_\xcF\bigr)(\si,m)=\id_{\txT M}+\rho(\Id_m)\circ\txA_i(\si,m)\,,\qquad\qquad (\si,m)\in O_i\x M\,,\qquad i\in I\,,
\qqq
written in terms of the local trivialisations $\xcF\t_i\colo\pi_\xcF^{-1}(O_i)\xrightarrow{\ \cong\ }O_i\x M$ corresponding to the $\xcP\t_i$,\ and of the anchor map $\rho$ of the tangent Lie algebroid $E$ of $\xcG$.
\eerop
\beroof 
The shadow connection 1-form $\Theta_\xcF$ is a smooth family of projectors onto $\txV\xcF$ and the latter subbundle within $\txT\xcF$ has $\txT M$ as a model.\ This explains the first (canonical) term $\id_{\txT M}$ in the local presentation of $\Theta_\xcF$.\ The structure of the non-canonical second term follows from the identity between the anchor $\rho\equiv\txT t$ and the local model of the tangent $\txT\xcD$ of the sitting-duck map inducing $\txH\xcF$ as in \eqref{eq:shadow-horizon}.
\eroof

\subsection{Existence and completeness}

\bethe[Existence of connections on a principaloid bundle]
On every principaloid bundle,\ there exists a connection.
\ethe
\beroof
We invoke Prop.\,\ref{prop:conn-1-Econn} and argue for the existence of a connection 1-form.\ To this end,\ consider a trivialising cover $\cO\equiv\{O_i\}_{i\in I}$ (assumed locally finite,\ see \cite{Nagata:1965}) of the base $\Si$ of the principaloid bundle $\xcP$,\ together with the associated trivialisations $\xcP\t_i$.\ Define the $I$-indexed family of smooth vector-valued 1-forms
\qq\nn
\cA_i:=\txT\xcP\t_i^{-1}\circ\jmath_2\circ\pi_2\circ\txT\xcP\t_i\,,\qquad i\in I
\qqq
in terms of the canonical maps:\ the projection $\pi_2\colo \txT(O_i\x\xcG)\cong\pr_1^*\txT O_i\oplus\pr_2^*\txT\xcG\to\pr_2^*\txT\xcG$ and the injection $\jmath_2\colo\pr_2^*\txT\xcG\emb\pr_1^*\txT O_i\oplus\pr_2^*\txT\xcG\cong\txT(O_i\x\xcG)$,\ where $\pr_1$ and $\pr_2$ are the canonical projections onto the first and the second factor in $O_i\x\xcG$,\ respectively.\ Use a smooth partition of unity $\{h_i\}_{i\in I}$ associated with $\cO$ to define the globally smooth vector-valued 1-forms
\qq\nn
\cA:=\sum_{i\in I}\,\pi_\xcP^*h_i\,\cA_i\,.
\qqq
These satisfy the axioms (GC1)--(GC3) by construction.
\eroof

Since the set of principaloid bundles include general fibre bundles,\ see Example \ref{eg:princ-as-gen},\ we do not expect connections to be complete always.\ However,\ completeness of the shadow connection on $\xcF$ is sufficient for the completeness of the connection on $\xcP$ under some additional assumptions:
\bethe
Let $\xcP$ be a principaloid bundle with structure groupoid $\grpd{\xcG}{M}$ and connection 1-form $\Theta$.\ If $\xcG$ admits a complete riemannian metric $\txg$ with respect to which right-invariant vector fields are bounded,\ and the connection $\Theta_{\xcF}$ on the shadow $\xcF$ of $\xcP$ induced from $\Theta$ is complete,\ then $\Theta$ is complete.\ In particular,\ if $\xcG$ admits a complete $\xcG$-invariant metric (in the sense of \cite{delHoyo:2015}),\ and $M$ is compact,\ then every connection $\Theta$ on $\xcP$ is complete.
\ethe
\beroof 
The problem of proving completeness of a connection on a fibre bundle $\pi_\xcE\colo\xcE\to B$ with base $\Si$ and typical fibre $F$ can,\ in general,\ be reduced to the following tasks (see \cite[Problem VII.12]{Greub:1972}):
\bit
\item[$1^\circ$] finding a horizontal lift $\widetilde\g{}_i$ to $\xcE$ for any curve $\g_i\colo]a_i,b_i[\to O_i$,\ defined for some $-\infty<a_i<b_i<\infty$,\ and contained entirely in an open set $O_i\subset B$ belonging to a given open cover $\{O_i\}_{i\in I}$ over which $E$ trivialises;
\item[$2^\circ$] checking that whenever there is a curve $\g_{ij}\colo]a_i,b_j[\to O_i\cup O_j,\ (i,j)\in I^{\x 2}_\cO$,\ defined for some $-\infty<a_i<a_j<b_i<b_j<\infty$ and glued smoothly from two smooth pieces $\g_i\colo]a_i,b_i[\to O_i$ and $\g_j\colo]a_j,b_j[\to O_j$ over $]a_j,b_i[\subset\g_{ij}^{-1}(O_i\cap O_j)$,\ the corresponding lifts $\widetilde\g{}_i$ and $\widetilde\g{}_j$ glue smoothly in $\xcE$.
\eit
In the former task,\ the key point is to ensure extendibility of every local lift from its domain given by a semi-closed interval $[c_i,d_i[\subset]a_i,b_i[$ resp.\ $]d_i,c_i]\subset]a_i,b_i[$,\ bound on one side by an initial condition in $\xcE_{\g_i(c_i)}$ at $c_i$,\ to its closure $[c_i,d_i]$ resp.\ $[d_i,c_i]$,\ see \cite[Sec.\,17.9]{Michor:2008}.\ This goal is attained through a limiting procedure which uses a complete metric structure on $F$.

Specialising the above to the case of a principaloid bundle $\xcE=\xcP$ with connection 1-form $\Theta$ for a proof of the first statement of the proposition,\ we start by considering the very formulation of the initial problem in $\pi_\xcP^{-1}(O_i)\cong O_i\x\xcG$.\ The assumed completeness of $\Theta_\xcF$ enables us to lift $\g_i$ horizontally to $\xcF$.\ In this manner,\ we obtain a smooth curve $\unl\g{}_i\colo]a_i,b_i[\to O_i\x M$ through an arbitrary point $m_i\equiv(\pr_2\circ\unl\g{}_i)(c_i)\in M$ for some $c_i\in]a_i,b_i[$,\ with a base projection $\pr_1\circ\unl\g{}_i=\g_i$ and a velocity field $\dot{\unl\g}_i\in\ker\,(\xcF\t_i^{-1\,*}\Theta_\xcF)$.\ In the next step,\ we restrict the local datum $\txA_i$ of $\Theta$ to $\unl\g{}_i$,\ whereby we obtain a time-dependent vector field
\qq\nn
\xcX_i\colo\, ]a_i,b_i[{}_{\pr_2\circ\unl\g{}_i}\hspace{-3pt}\x_{t}\hspace{-1pt}\xcG\too\ker\,\txT s\subset\txT\xcG,\ (t,g)\longmapsto\txT_{\Id_{\pr_2\circ\unl\g{}_i}}r_g\bigl(\txA_i\bigl(\dot{\unl\g}_i\bigr)\bigr)\equiv\G_i\bigl(\g_i(t),g\bigr)\bigl(\dot{\unl\g}_i\bigr)\,.
\qqq
By assumption,\ this field is bounded with respect to the complete metric $\txg$ on $\xcG$ over every compact interval within $]a_i,b_i[$.\ Let us,\ next,\ consider a lift of $\g_i$ to $O_i\x\xcG$ through the initial point $g_i\in t^{-1}(\{m_i\})$ at $c_i$.\ Such a lift takes the form $\widetilde\g{}_i=(\g_i,\widehat\g{}_i)$,\ where $\widehat\g{}_i$ is an integral curve of $\xcX_i$ through $\widehat\g{}_i(c_i)=g_i$.\ {\it A priori},\ the latter exists and is unique locally,\ that is over some interval,\ say,\ $[c_i,d_i[\subset]a_i,b_i[$.\ In order to extend the local solution all the way to $b_i$---and likewise to $a_i$ on the other side of $c_i$---we need to verify that the local solution reaches a point $\widetilde g{}_i\in\xcG$ in the fibre over $d_i$ as we take $c_i\to d_i$.\ However,\ using the aforementioned boundedness of $\xcX_i$ over $[c_i,d_i]$,
\qq\nn
\forall\,(t,g)\in[c_i,d_i]{}_{\pr_2\circ\unl\g{}_k}\hspace{-3pt}\x_{t}\hspace{-1pt}\xcG\colo \Vert\xcX_i(t,g)\Vert_\txg\leq C_i<\infty\,,
\qqq
we may estimate,\ for any $n_i\in\bN^\x$,\ the metric distance between the consecutive end-points $\widehat\g{}_i(d_i-\frac{1}{n_i})$ and $\widehat\g{}_i(d_i-\frac{1}{n_i+1})$ as
\qq\nn
0\leq d_\txg\bigl(\widehat\g{}_i(d_i-\tfrac{1}{n_i}),\widehat\g{}_i(d_i-\tfrac{1}{n_i+1})\bigr)\leq\int_{d_i-\frac{1}{n_i}}^{d_i-\frac{1}{n_i+1}}\,\sfd t\,\Vert\dot{\widehat\g}{}_i(t)\Vert_\txg=\int_{d_i-\frac{1}{n_i}}^{d_i-\frac{1}{n_i+1}}\,\sfd t\,\Vert\xcX_i\bigl(t,\widehat\g{}_i(t)\bigr)\Vert_\txg\leq \tfrac{C_i}{n_i(n_i+1)}\xrightarrow[\ n_i\to\infty\ ]{}0\,.
\qqq
The assumed completeness of the metric space $(\xcG,d_\txg)$ enables us to infer from this estimate that $\widehat\g{}_i$,\ and so also $\widetilde\g{}_i$ do,\ indeed,\ extend to $d_i$,\ with 
\qq\nn
\widetilde\g{}_i(d_i)=(\g_i(d_i),\lim_{n_i\to\infty}\,\widehat\g{}_i(d_i-\tfrac{1}{n_i}))\,.
\qqq
Taking the endpoint $\widetilde\g{}_i(d_i)$ of this extension as a new initial condition,\ we may,\ next,\ proceed with the extension to the right until we reach $b_i$,\ see the proof of \cite[Prop.\,23.9]{Michor:2008}.\ Thus,\ we obtain a $\Theta$-horizontal lift $\widetilde\g{}_i$ of $\g_i$ for any interval $]a_i,b_i[$ through an arbitrary point $g_i\equiv(\pr_2\circ\widetilde\g{}_i)(c_i)\in t^{-1}(\{m_i\})$ at $c_i\in]a_i,b_i[$.\ A smooth gluing of these lifts over intersections $O_{ij}$ is ensured by the gluing law \eqref{eq:glue-Aij} for the corresponding local data of $\Theta$.

For a proof of the second statement of the proposition,\ note,\ first,\ that every connection on a fibre bundle with a compact typical fibre is complete---see,\ {\it e.g.},\ \cite[Sec.\,17.9]{Michor:2008}---and so the shadow connection $\Theta_\xcF$ is complete for every $\Theta$.\ Moreover,\ the norm squared of a right-invariant vector field with respect to a right-invariant metric on $\xcG$,\ viewed as a function on $\xcG$ with values in $\bR_{\geq 0}$,\ is a pullback along $t$ of a smooth function on the identity bisection,\ which is diffeomorphic with $M$,\ and hence compact by assumption.\ As a continuous function on a compact set,\ the norm is therefore bounded,\ and so both assumption from the first statement of the proposition are reproduced.
\eroof

\brem
Conditions under which there exists a riemannian metric on a Lie groupoid,\ and consequences of its existence,\ were discussed at length in \cite{delHoyo:2015}.
\erem

\subsection{$\bB$-invariance}

In the case of a bisection-complete strucure groupoid,\ an alternative definition of the principaloid bundle in terms of the $\bB$-action is available.\ The underlying characterisation of the geometric object through its symmetry extends to connections in a natural way.\ By construction,\ a connection on $\xcP$,\ {\it i.e.},\ the choice of the horizontal distribution $\txH\xcP$,\ is $\xcG$-invariant,\ and hence also $\bB$-invariant in the light of Prop.\,\ref{prop:G-in-B_act}.\ That the latter property determines the connection uniquely,\ for the distinguished class of strucure groupoids introduced in Def.\,\ref{def:Id-red},\ follows from the analysis presented below.

\bedef\label{def:Bisec-principal-Ehresmann} 
A (compatible) {\bf connection} on a $\xcG$-fibred $\bB$-bundle $\widehat\xcP$ is a $\bB$-invariant Ehresmann connection $\txT\widehat\xcP=\txV\widehat\xcP\oplus\txH\widehat\xcP$,\ {\it i.e.},\ a connection such that,\ for every $p\in\widehat\xcP$ and $\b \in \bB$,\ the horizontal distribution $\txH\widehat\xcP$ satisfies
\qq \label{eq:B-conn} 
\txT_p\widehat\xcR{}_\b (\txH_p\widehat\xcP) = \txH_{\widehat\xcR{}_\b(p)}\widehat\xcP\,.
\qqq
\exdef 
\brem
For a general groupoid $\xcG$ with an object manifold $M$ of non-zero dimension,\ the action of the group of bisections $\bB$ on $\xcP$ is far from free.\ Choosing a bisection $\beta$ in the stabiliser of a point $p \in \xcP$,\ $\xcR_\b(p)=p$,\ the map $\txT\xcR_\b$ may not be the identity. Correspondingly,\  condition \eqref{eq:B-conn} then poses severe restrictions on the horizontal subspace $\txH_p\xcP$, namely that it lie inside $\ker \txT \mu$.
\erem
\noindent As before,\ we may equivalently work with a smooth distribution of projectors over $\widehat\xcP$.
\bedef\label{def:Bisec-principal-conn} 
Let $\widehat\xcP$ be a $\xcG$-fibred $\bB$-bundle.\ A (compatible) {\bf connection 1-form} on $\widehat\xcP$ is a $\txV\widehat\xcP$-valued 1-form $\widehat\Theta\in\Om^1(\widehat\xcP,\txV\widehat\xcP)\equiv\G(\txT^*\widehat\xcP\ox\txV\widehat\xcP)$
such that,\ when viewed as an element in $\End{(\txT\widehat\xcP)}$:
\bit
\item[(BC1)] $\widehat\Theta\rstr_{\txV\widehat\xcP}=\id_{\txV\widehat\xcP}$\, (a projection onto $\txV\widehat\xcP$);
\item[(BC2)] $\forall\ \b\in\bB\colo \quad \widehat\Theta\circ\txT\xcR_\b=\txT\xcR_\b\circ\widehat\Theta$\, ($\bB$-equivariance).
\eit
\exdef
\noindent And,\ again,\ we readily establish the equivalence of the two definitions.
\berop\label{prop:Bequiv-conns-as-forms}
Connections on a $\xcG$-fibred $\bB$-bundle $\widehat\xcP$ are in a one-to-one correspondence with connection 1-forms on $\widehat\xcP$.
\eerop
\beroof
The choice of a decomposition $\txT\widehat\xcP=\txV\widehat\xcP\oplus\txH\widehat\xcP$ determines a 1-form $\widehat\Theta$ by declaring the latter to vanish on the horizontal subbundle and to be the identity on the vertical one,\ which in particular implies (BC1) by construction.\ Conversely,\ given $\widehat\Theta$,\ property (BC1) determines a decomposition with $\txH\widehat\xcP:=\ker\,\widehat\Theta$. 

Note that (BC2) is trivially satisfied on vertical vectors due to (BC1).\ When applied to a horizontal vector,\ {\it i.e.},\ one in the kernel of $\widehat\Theta$,\ (BC2) reduces to 
\qq\nn
\widehat\Theta\bigl(\xcR_\b(p)\bigr) \circ \txT_p\xcR_\b \vert_{\txH_p\widehat\xcP} = 0
\qqq
for every $p\in\widehat\xcP$ and $\beta \in \bB$.\ This identity follows immediately from identity \eqref{eq:B-conn} and,\ {\it vice versa},\ identity \eqref{eq:B-conn} is implied by
the one above since $\txT\xcR_\b$,\ when applied to any tangent space,\ is an isomorphism and thus preserves the dimension of subspaces.
\eroof

\bethe
For an $\Id$-reducibile Lie groupoid $\xcG$,\ there is a one-to-one correspondence between $\xcG$-fibred $\bB$-bundles with compatible connection and principaloid bundles with compatible connection.
\ethe
\beroof
In virtue of Thm.\,\ref{thm:Bequiv-as-princoid},\ the assumption about $\xcG$ ensures that every $\xcG$-fibred $\bB$-bundle is a principaloid bundle and {\it vice versa}.\ Moreover,\ the $\bB$-action on $\xcP$ is induced from the $\xcG$-action by Prop.\,\ref{prop:G-in-B_act},\ and so a principaloid bundle with compatible connection is automatically a $\xcG$-fibred $\bB$-bundle with compatible connection.

Let us,\ next,\ consider a model of a compatible connection on a $\xcG$-fibred $\bB$-bundle $\widehat\xcP$ in one of $\bB$-equivariant trivialisations $\widehat\xcP\t_i\colo\pi_{\widehat\xcP}^{-1}(O_i)\xrightarrow{\ \cong\ }O_i\x\xcG$ associated with a cover $\{O_i\}_{i\in I}$.\ In the light of Prop.\,\ref{prop:Bequiv-conns-as-forms},\ the local model of the horizontal distribution $\txH_p\widehat\xcP$ in $\widehat\xcP\t_i$ is given by the kernel of a 1-form
\qq\label{eq:loc-BTheta-gen}
\bigl(\widehat\xcP\t_i^{-1\,*}\widehat\Theta\bigr)(\si,g)=\id_{\txT\xcG}\rstr_{\txT_g\xcG}+\widehat\G{}_i(\si,g)\,,
\qqq
obtained by pullback of the connection 1-form $\widehat\Theta$ along $\widehat\xcP\t_i$,\ with a vector-valued 1-form $\widehat\G{}_i\colo\pr_1^*\txT O_i\to\pr_2^*\txT\xcG$.\ Here,\ we use the same conventions as in the proof of Thm.\,\ref{thm:loc-data-conn}.\ Using local coordinates $\{\si^\mu\}^{\mu\in\ovl{1,\dim\,\Si}}$ on a neighbourhood of $\si\in O_i$,\ we write
\qq\nn
\widehat\G{}_i(\si,\cdot)=\txd\si^\mu(\si)\ox\widehat\G{}^\si_{i\,\mu}(\cdot)\,,
\qqq
where $\widehat\G{}^\si_{i\,\mu}\in\G(\txT\xcG)$.\ Invoking the assumption of $\bB$-equivariance of $\widehat\Theta$,\ we note that while the identity operator $\id_{\txT\xcG}$ trivially commutes with $\txT R_\b$,\ the second term in \eqref{eq:loc-BTheta-gen} does so only if the $\widehat\G{}_{i\,\mu}^\si$ are $R(\bB)$-invariant,
\qq\nn
\widehat\G{}_{i\,\mu}^\si\in\G(\txT\xcG)^{R(\bB)}\,.
\qqq
Therefore,\ by Prop.\,\ref{prop:RinV-in-kerTs},\ we establish a restriction on the codomain of the $\widehat\theta{}_i$ in the following form:
\qq\nn
\widehat\G{}_i\colo\txT O_i\too\ker\,\txT s\,.
\qqq
Accordingly,\ horizontal vector fields are restricted to lie in the distribution $\pr_1^*\txT O_i\oplus\pr_2^*\ker\,\txT s$ within the model tangent bundle $\pr_1^*\txT O_i\oplus\pr_2^*\txT\xcG$,\ and so in the kernel of the local presentation of $\txT\mu$.

Invoking Prop.\,\ref{prop:rasoR},\ we conclude that $\txT\varrho_g$ is,\ for every $g\in t^{-1}(\{\mu(p)\})$,\ represented by $\txT R_{\beta_{g^{-1}}^{-1}}$ acting on $\pr_2^*\ker\,\txT s$.\ But the trivialisations are  $\bB$-equivariant,\ and therefore $\txT\varrho_g$ becomes $\txT\xcR_{\beta_{g^{-1}}^{-1}}$ globally over $\xcP$.\ The claim of the proposition now follows directly from the defining identity \eqref{eq:B-conn} and another application of Prop.\,\ref{prop:rasoR}.
\eroof

\section{The Atiyah sequence,\ gauge transformations,\ and the trident of $\xcP$}\label{sec:Atiyah}

\subsection{The Atiyah--Ehresmann groupoid and its short exact sequence}

The Atiyah--Ehresmann groupoid of an ordinary principal $\txG$-bundle is $(\sfP \x \sfP)/\txG$,\ where the quotient is relative to the diagonal action of $\txG$.\ This is a bundle associated to $\sfP$ by the defining right action $r$ on $\sfP$.\ Thus,\ it can be described as 
\qq\nn 
\bigsqcup_{i\in I}\,\bigl( O_i\x \sfP \bigr)/\sim_{r_{g_{\cdot\cdot}^{-1}}}
\qqq 
for the transition maps $g_{ij} \colo O_{ij} \to \txG$ of $\sfP$.\ We generalise this to  principaloid bundles as follows.
\bedef
We shall call the fibre bundle $\pi_{{\rm At}(\xcP)}\colo{\rm At}(\xcP)\to\Si$ with typical fibre $\xcP$ and model
\qq\nn
\bigsqcup_{i\in I}\,\bigl(\xcP\x O_i\bigr)/\sim_{\xcR_{\b_{\cdot\cdot}^{-1}}}\too\Si,\ [(p,\si)]\longmapsto\si\,,
\qqq
the {\bf Atiyah bundle of} $\xcP$.\ We shall also denote,\ in the above model,
\qq\nn
\pi_1\colo\bigsqcup_{i\in I}\,\bigl(\xcP\x O_i\bigr)/\sim_{\xcR_{\b_{\cdot\cdot}^{-1}}}\too\Si,\ [(p,\si)]\longmapsto\pi_\xcP(p)\,.
\qqq
\exdef

\brem
${\rm At}(\xcP)$ is canonically induced by the principaloid bundle $\xcP$ through its transition 1-cocycle $\{\b_{ij}\}_{i,j\in I^{\x 2}_\cO}$ realised by the defining action \eqref{eq:R-act-Poid}.
\erem

\bethe\label{thm:AtPoid}
The pair $({\rm At}(\xcP),\xcF)$ carries a canonical structure of a Lie groupoid,\ fitting into the following short exact sequence of Lie groupoids:
\qq\label{ses:Atiyah}\qquad
\alxydim{@C=2.5cm@R=1.5cm}{ {\rm Ad}(\xcP) \ar@{^{(}->}[r]^{j_{{\rm Ad}(\xcP)}} \ar@<.5ex>[d]_{\sfT\vert\;\; } \ar@<-.5ex>[d]^{\;\; \sfS\vert} & {\rm At}(\xcP) \ar@{->>}[r]^{\pi} \ar@<.5ex>[d]_{\sfT \;\;} \ar@<-.5ex>[d]^{\;\; \sfS} & \Si\x\Si{} \ar@<.5ex>[d]_{\pr_1\;\; } \ar@<-.5ex>[d]^{\;\; \:\pr_2}\\ \xcF \ar@{=}[r] & \xcF \ar[r]_{\pi_\xcF} & \Si }\,.
\qqq
Here,\ $\pi=(\pi_1,\pi_{{\rm At}(\xcP)})$ and $j_{{\rm Ad}(\xcP)}$ is the embedding of ${\rm Ad}(\xcP) = \pi^{-1}(\Id(\Sigma))$,\ where $\Id(\Sigma) \subset \Si \x \Si$ is the identity bisection.  
\ethe

\brem\label{rem:Atiyah-as-bndl} 
The $\pi$-fibres of ${\rm At}(\xcP)$  are isomorphic to the Lie groupoid $\xcG$. We can,\ therefore,\ read the exact sequence \eqref{ses:Atiyah} also as saying that $({\rm At}(\xcP), \mathrm{Pair}(\Sigma), \xcG, \pi)$ is a fibre-bundle object in the category of Lie groupoids: 
\qq\nn
\alxydim{@C=.5cm@R=1.cm}{\xcG \: \ar@{^{(}.>}[r] & {\rm At}(\xcP) \ar[d]^{\pi} \\ & {\rm Pair}(\Si)}\,.
\qqq
\erem

The proof of Thm.\,\ref{thm:AtPoid} bases upon a pair of propositions and a lemma stated below. The Remark above then follows from noting that the fibre of ${\rm At}(\xcP)$, which therefore also is the typical fibre of ${\rm Ad}(\xcP)$, is the Lie groupoid $\grpd{\xcG}{M}$.

\berop\label{prop:At-as-Foid} 
The pair $({\rm At}(\xcP),\xcF)$ composes,\ in a canonical way,\ a Lie groupoid,\ with the following structure maps,\ modelled on those of the Lie groupoid $\xcG$ in local trivialisations, 
\qq\label{diag:Grmatbndldiag}\qquad\qquad
\alxydim{@C=2cm@R=1.5cm}{ {\rm At}(\xcP)\x_{\xcF}{\rm At}(\xcP) \ar[r]^{\quad\sfM} & {\rm At}(\xcP) \ar[r]^\sfJ & {\rm At}(\xcP) \ar@<.25ex>[r]^{\sfS} \ar@<-.25ex>[r]_{\sfT} & \xcF \ar@/_2pc/[l]_{\sfI} }\,.
\qqq
\eerop
\beroof
We begin by noting that ${\rm At}(\xcP)$ admits further (model) resolution given by
\qq\label{eq:AtP-loc-res}
{\rm At}(\xcP)\cong\bigsqcup_{i,j\in I}\,\bigl(O_i^{(1)}\x\xcG\x O_j^{(2)}\bigr)/\sim_{L^{(1)}_{\b_{\cdot\cdot}}\circ R^{(2)}_{\b_{\cdot\cdot}^{-1}}}\,,
\qqq
with local charts glued by the identifications
\qq\nn
(\si_1,g,\si_2,k,l)\sim\bigl(\si_1,\b_{ik}(\si_1)\lact g\ract\b_{lj}(\si_2),\si_2,i,j\bigr)\,.
\qqq
Accordingly,\ we define the structure maps as follows 
\qq
&\sfS\colo {\rm At}(\xcP)\too\xcF,\ [(\si_1,g,\si_2,i,j)]\longmapsto[(\si_2, s( g),j)]\,,&\cr\cr
&\sfT\colo {\rm At}(\xcP)\too\xcF,\ [(\si_1,g,\si_2,i,j)]\longmapsto[(\si_1, t( g),i)]\,,&\nn\\ \label{eq:STIJ-At} \\
&\sfI\colo \xcF\too{\rm At}(\xcP),\ [(\si,m,i)]\longmapsto[(\si,\Id_m,\si,i,i)]\,,&\cr\cr
&\sfJ\colo {\rm At}(\xcP)\too{\rm At}(\xcP),\ [(\si_1,g,\si_2,i,j)]\longmapsto[(\si_2,g^{-1},\si_1,j,i)]&\nn
\qqq
and
\qq\nn
\sfM\colo {\rm At}(\xcP){}_{\sfS}\hspace{-3pt}\x_{\sfT}\hspace{-1pt}{\rm At}(\xcP)\too {\rm At}(\xcP),\ \bigl([(\si_1,g,\si_2,i,j)],[(\si_2,h,\si_3,k,l)]\bigr)\longmapsto[(\si_1,g.\bigl(\b_{jk}(\si_2)\lact h\bigr),\si_3,i,l)]\,.
\qqq
The well-definedness of the above maps is ensured by the identifications
\bit
\item for the source map,\ at $(\si_1,\si_2)\in O_{ik}\x O_{jl}$,
\qq\nn
\sfS\bigl([(\si_1,\b_{ki}(\si_1)\lact g\ract\b_{jl}(\si_2),\si_2,k,l)]\bigr)=[(\si_2,\b_{lj}(\si_2)\ulact s(g),l)]\equiv[(\si_2,s(g),j)]\,,
\qqq
see \Reqref{eq:BisAct-vs-str-i} in Appendix \ref{app:useful};
\item for the target map,\ at $(\si_1,\si_2)\in O_{ik}\x O_{jl}$,
\qq\nn
\sfT\bigl([(\si_1,\b_{ki}(\si_1)\lact g\ract\b_{jl}(\si_2),\si_2,k,l)]\bigr)=[(\si_1,\b_{ki}(\si_1)\ulact t(g),k)]\equiv[(\si_1,t(g),i)]\,,
\qqq
see \Reqref{eq:BisAct-vs-str-ii};
\item for the identity map,\ at $\si\in O_{ij}$,
\qq\nn
\sfI\bigl([(\si,\b_{ji}(\si)\ulact m,j)]\bigr)=[(\si,C_{\b_{ji}(\si)}(\Id_m),\si,j,j)]\equiv[(\si,\Id_m,\si,i,i)]\,,
\qqq
see \Reqref{eq:BisAct-vs-str-iii};
\item for the inverse map,\ at $(\si_1,\si_2)\in O_{ik}\x O_{jl}$,
\qq\nn
\sfJ[(\si_1,\b_{ki}(\si_1)\lact g\ract\b_{jl}(\si_2),\si_2,k,l)]=[(\si_2,\b_{lj}(\si_2)\lact g{}^{-1}\ract\b_{ik}(\si_2),\si_1,l,k)]\equiv[(\si_2, g{}^{-1},\si_1,j,i)]\,,
\qqq
see \Reqref{eq:BisAct-vs-str-iv};
\item for the multiplication map,\ at $(\si_1,\si_2,\si_3)\in O_{im}\x O_{jkno}\x O_{lp}$,
\qq\nn
&&\sfM\bigl([(\si_1,\b_{mi}(\si_1)\lact g\ract\b_{jn}(\si_2),\si_2,m,n)],[(\si_2,\b_{ok}(\si_2)\lact h\ract\b_{lp}(\si_3),\si_3,o,p)]\bigr)\cr\cr
&=&[(\si_1,(\b_{mi}(\si_1)\lact g\ract\b_{nj}(\si_2)^{-1}).(\b_{no}(\si_2)\lact(\b_{ok}(\si_2)\lact h)\ract\b_{lp}(\si_3)),\si_3,m,p)]\cr\cr
&=&[(\si_1,(\b_{mi}(\si_1)\lact g).((\b_{jn}(\si_2)\lact(\b_{no}(\si_2)\cdot\b_{ok}(\si_2))\lact h)\ract\b_{lp}(\si_3)),\si_3,m,p)]\cr\cr
&=&[(\si_1,(\b_{mi}(\si_1)\lact g).((\b_{jn}(\si_2)\cdot\b_{nk}(\si_2))\lact h\ract\b_{lp}(\si_3)),\si_3,m,p)]\cr\cr
&=&[(\si_1,\b_{mi}(\si_1)\lact(g.(\b_{jk}(\si_2)\lact h))\ract\b_{lp}(\si_3),\si_3,m,p)]\equiv[(\si_1,g.(\b_{jk}(\si_2)\lact h),\si_3,i,l)]\,,
\qqq
see \Reqref{eq:BisAct-vs-str-vi}.
\eit
The constitutive relations between the structure maps are implied by the same relations for their local models.\ The only seemingly non-obvious ones are those involving the multiplication map $\sfM$,\ but in the light of the above consistency check,\ we may rewrite the definition of $\sfM$ as
\qq\label{eq:M-At}\qquad\qquad
\sfM\colo {\rm At}(\xcP){}_{\sfS}\hspace{-1pt}\x_{\sfT}\hspace{-1pt}{\rm At}(\xcP)\too {\rm At}(\xcP),\ \bigl([(\si_1,g_1,\si_2,i,j)],[(\si_2,g_2,\si_3,j,l)]\bigr)\longmapsto[(\si_1,g_1.g_2,\si_3,i,l)]\,.
\qqq
\eroof

\bedef
We shall call $\grpd{{\rm At}(\xcP)}{\xcF}$ the {\bf Atiyah--Ehresmann groupoid} of $\xcP$.\ The short exact sequence \eqref{ses:Atiyah} of Lie groupoids shall be referred to as the {\bf Atiyah sequence} for $\xcP$.
\exdef

\belem\label{lem:At-Poid-over-PairSig}
The pair $(\pi,\pi_\xcF)$ is a Lie-groupoid epimorphism from ${\rm At}(\xcP)$ to ${\rm Pair}(\Si)$.
\elem
\beroof
The statement of the Lemma follows straightforwardly from the fact that the structure maps \eqref{eq:STIJ-At} and \eqref{eq:M-At} cover the corresponding structure maps of ${\rm Pair}(\Si)$.
\eroof

\berop\label{prop:trans-morphs}
Let $\grpd{\xcG_1}{M_1}$ and $\grpd{\xcG_2}{M_2}$ be Lie groupoids,\ and let
\qq\nn
\alxydim{@C=2cm@R=1.5cm}{ \xcG_1 \ar[r]^{\chi} \ar@<.5ex>[d]_{t_1\;\; } \ar@<-.5ex>[d]^{\;\; s_1} & \xcG_2 \ar@<.5ex>[d]_{t_2\;\; } \ar@<-.5ex>[d]^{\;\; s_2} \\ M_1 \ar[r]_{\chi_0}  & M_2 }
\qqq
be a morphism of Lie groupoids,\ whose component $\chi$ is transverse\footnote{Recall that a smooth manifold map $f\colo M\to N$ is said to be \emph{transverse} to a submanifold $S\subset N$ if at an arbitrary point $m\in f^{-1}(S)$,\ the following condition is satisfied:\ $\txT_m f(\txT_m M)+\txT_{f(m)}S=\txT_{f(m)}N$.} to the identity bisection $\Id(M_2)\subset\xcG_2$.\ The preimage $\chi^{-1}(\Id(M_2))$ of $\Id(M_2)$ is a Lie subgroupoid of $\,\grpd{\xcG_1}{M_1}$,\ {\it i.e.},\ it is,\ in particular,\ a Lie groupoid.
\eerop
\beroof
For $\chi^{-1}(\Id(M_2))$ to be a Lie subgroupoid in $\xcG_1$,\ the former has to be a submanifold in $\xcG_1$.\ In virtue of a classic variant of the Level-Set Theorem for submanifolds,\ this is ensured by the transversality of $\chi$.\ 

At this stage,\ it remains to check that the subset $\chi^{-1}(\Id(M_2))\subset\xcG_1$ is a subgroupoid.\ To this end,\ consider an arbitrary pair $(g_1,g'_1)\in\chi^{-1}(\Id(M_2))\,{}_{s_1}\hspace{-3pt}\x_{t_1}\hspace{-1pt}\chi^{-1}(\Id(M_2))$.\ There,\ then,\ exist points $m_2,m'_2\in M_2$ such that $\chi(g_1)=\Id_{m_2}$ and $\chi(g'_1)=\Id_{m'_2}$.\ As $\chi_0\circ s_1=s_2\circ\chi$ and $\chi_0\circ t_1=t_2\circ\chi$,\ we establish the identity
\qq\nn
m_2'=t_2(\Id_{m_2'})=t_2\bigl(\chi(g_1')\bigr)=\chi_0\bigl(t_1(g_1')\bigr)=\chi_0\bigl(s_1(g_1)\bigr)=s_2\bigl(\chi(g_1)\bigr)=s_2(\Id_{m_2})=m_2\,.
\qqq
Since $\chi$ is a morphism, this implies
\qq\nn
\chi(g_1.g_1')=\chi(g_1).\chi(g_1')=\Id_{m_2}.\Id_{m_2'}=\Id_{m_2}\,,
\qqq
which leads to the desired conclusion $g_1.g_1'\in\chi^{-1}(\Id(M_2))$.\ Similarly,\ for every $g_1\in\chi^{-1}(\Id(M_2))$,\ with $\chi(g_1)=\Id_{m_2}$,\ we have $\chi(g_1^{-1})=\chi(g_1)^{-1}=\Id_{m_2}^{-1}=\Id_{m_2}$,\ and so also $g_1^{-1}\in\chi^{-1}(\Id(M_2))$.
\eroof

\noindent\emph{Proof of Theorem \ref{thm:AtPoid}.}\ The existence of the structure of a Lie groupoid on $({\rm At}(\xcP),\xcF)$ is stated in Prop.\,\ref{prop:At-as-Foid}.\ The exactness of sequence \eqref{ses:Atiyah} at its node ${\rm Pair}(\Si)$ then follows from Lemma \ref{lem:At-Poid-over-PairSig}.\ As a set,\ ${\rm Ad}(\xcP)$ fits into the short exact sequence by definition,\ and the only thing that remains to be proven is the embedding of the pair $({\rm Ad}(\xcP),\xcF)$ in the Atiyah--Ehresmann groupoid $\grpd{{\rm At}(\xcP)}{\xcF}$ as a Lie subgroupoid.\ Since $(\pi,\pi_\xcF)$ is an epimorphism of Lie groupoids,\ its arrow component $\pi$ is automatically transverse to $\Id(\Si)\subset\Si\x\Si$ (as a submersion),\ and so we conclude the present proof by invoking Prop.\,\ref{prop:trans-morphs}. \qed 

\bedef
We shall call ${\rm Ad}(\xcP)$ the {\bf adjoint bundle},\ and  $\grpd{{\rm Ad}(\xcP)}{\xcF}$ the {\bf adjoint groupoid} of $\xcP$.
\exdef

\berop
The adjoint bundle ${\rm Ad}(\xcP)$ has a model
\qq\label{eq:clutch-xcA}
{\rm Ad}(\xcP)\cong\bigsqcup_{i\in I}\,\bigl( O_i\x \xcG\bigr)/\sim_{C_{\b_{\cdot\cdot}}}\,,
\qqq
written in terms of the transition 1-cocycle $\{\b_{ij}\}_{(i,j)\in I^{\x 2}_\cO}$ of $\xcP$,\ which we realise by conjugation of $\xcG$ by $\bB$. 
\eerop
\beroof
By definition,\ we have---in the notation of \eqref{eq:AtP-loc-res}--- 
\qq\nn
{\rm Ad}(\xcP)\equiv\pi^{-1}(\Id(\Sigma))=\bigsqcup_{i,j\in I}\,\bigl(O_i^{(1)}\x_\Si\bigl(\xcG\x O_j^{(2)}\bigr)\bigr)/\sim_{L^{(1)}_{\b_{\cdot\cdot}}\circ R^{(2)}_{\b_{\cdot\cdot}^{-1}}}\,.
\qqq
Thus,\ a point in ${\rm Ad}(\xcP)$ is an equivalence class $[(\si,g,\si,i,j)]\equiv[(\si,g\ract\b_{ji}(\si),\si,i,i)]\equiv[(\si,\b_{ji}(\si)\lact g,\si,j,j)]$.\ We define (smooth) maps
\qq\nn
\jmath_{{\rm Ad}(\xcP)}\colo\bigsqcup_{i\in I}\,\bigl( O_i\x \xcG\bigr)/\sim_{C_{\b_{\cdot\cdot}}}\too{\rm Ad}(\xcP),\ [(\si,g,i)]\longmapsto[(\si,g,\si,i,i)]
\qqq
and
\qq\nn
\iota_{{\rm Ad}(\xcP)}\colo{\rm Ad}(\xcP)\too\bigsqcup_{i\in I}\,\bigl( O_i\x \xcG\bigr)/\sim_{C_{\b_{\cdot\cdot}}},\ [(\si,g,\si,i,j)]\longmapsto[(\si,g\ract\b_{ji}(\si),i)]\,,
\qqq
which are readily checked to be each other's inverses.\ Their well-definedness is a consequence of the identities,\ written for $\si\in O_{ij}$ and $\si'\in O_{ijkl}$,
\qq\nn
\jmath_{{\rm Ad}(\xcP)}\bigl([(\si,C_{\b_{ji}(\si)}(g),j)]\bigr)\equiv[(\si,C_{\b_{ji}(\si)}(g),\si,j,j)]=[(\si,\b_{ji}(\si)\lact g\ract\b_{ji}(\si)^{-1},\si,j,j)]=[(\si,g,\si,i,i)]\,,
\qqq
and
\qq\nn
&&\iota_{{\rm Ad}(\xcP)}\bigl([(\si',\b_{ki}(\si')\lact g\ract\b_{lj}(\si')^{-1},\si',k,l)]\bigr)\equiv[(\si',\bigl(\b_{ki}(\si')\lact g\ract\b_{lj}(\si')^{-1}\bigr)\ract\b_{lk}(\si'),k)]\cr\cr
&=&[(\si',\b_{ki}(\si')\lact g\ract\b_{jl}(\si')\cdot\b_{lk}(\si'),k)]=[(\si',\b_{ki}(\si')\lact g\ract\b_{jk}(\si'),k)]\cr\cr
&=&[(\si',\b_{ki}(\si')\lact\bigl(g\ract\b_{jk}(\si')\cdot\b_{ki}(\si')\bigr)\ract\b_{ki}(\si')^{-1},k)]=[(\si',g\ract\b_{jk}(\si')\cdot\b_{ki}(\si'),i)]\cr\cr
&=&[(\si',g\ract\b_{ji}(\si'),i)]\,,
\qqq
where in the second,\ third and fourth lines,\ we have used commutativity of left--mul\-ti\-pli\-ca\-tions of $\xcG$ by $\bB$ with right-mul\-ti\-pli\-ca\-tions,\ alongside the 1-cocycle condition satisfied by the $\b_{ij}(\si')$.
\eroof

\brem
Applying the Lie functor to the short exact sequence \ref{ses:Atiyah},\
we obtain an exact sequence of Lie algebroids:
\qq\label{ses:Atiyah-alg}\qquad
\alxydim{@C=2.5cm@R=1.5cm}{ \xcE_0 \ar@{^{(}->}[r]^{j_{\xcE_0}} \ar[d] & \xcE \ar@{->>}[r]^{\pi_\xcE} \ar[d] & \txT\Si \ar[d] \\ \xcF \ar@{=}[r] & \xcF \ar[r]_{\pi_\xcF} & \Si}\,.
\qqq
We shall call $\xcE$ the {\bf Atiyah algebroid} of $\xcP$.\ The fibre of $\xcE_0=\ker\,\pi_\xcE$ is the Lie algebroid $E=\mathrm{Lie}(\xcG)$.\ This sequence appeared in the form of a $Q$-bundle in \cite{Kotov:2007nr},\ specifically for the case where the $Q$-manifolds correspond to Lie algebroids. 
\erem

\beg
If $\xcG=\txG$,\ then conjugation coincides with the adjoint action $\Ad\colo\txG\x\txG\to\txG,\ (g,h)\mapsto g h g^{-1}$ of the Lie group on itself,\ and we retrieve the adjoint bundle ${\rm Ad}(\xcP)\cong\Ad\,\txP$,\ associated with the pricipal $\txG$-bundle $\xcP\equiv\txP$ through that adjoint action,\ $\Ad\,\txP\equiv\txP\x_\Ad\txG$.\ Moreover,\ ${\rm At}(\xcP)$ is the standard Atiyah--Ehresmann groupoid of the principal $\txG$-bundle and $\xcE$ in Diag.\,\eqref{ses:Atiyah-alg} reduces to the standard Atiyah algebroid $\txT\txP/\txG$.
\eeg

\beg If $\xcG={\rm Pair}(M)$,\ then conjugation by $\b\equiv(f,\id_M)\in\bB$,\ with $f\in\Diff(M)$,\ is just the cartesian square of $f$,\ {\it i.e.},\ $C_\b\equiv f\x f$,\ and the adjoint bundle of the principaloid bundle $\xcP\cong\xcF\x M$ is the $\Si$-fibred square ${\rm Ad}(\xcP)\cong\xcF\x_\Si\xcF$ of $\xcF$.\ Here,\ ${\rm At}(\xcP)\equiv{\rm At}(\xcF\x M)$ is just the pair groupoid $\grpd{\xcF\x\xcF}{\xcF}$ with the groupoid morphism $\pi$ in Diag.\,\eqref{ses:Atiyah} given by  $\pi_\xcF \x \pi_\xcF$.\ Correspondingly,\ $\xcE = \txT\xcF$ and $\pi_\xcE = \txT\pi_\xcF$.
\eeg

The Atiyah sequence \eqref{ses:Atiyah} gives rise to a specific subgroup of bisections in $\grpd{\mathrm{At}(\xcP)}{\xcF}$, which will become of relevance in what follows.  

\bedef
A bisection $\b$ of the Atiyah--Ehresmann groupoid $\grpd{\mathrm{At}(\xcP)}{\xcF}$ for which there exists a bisection $\unl\b$ of ${\rm Pair}(\Si)$ such that the two  rectangles  at the right in the diagram
\qq\label{eq:Atiyah-bisec}\qquad
\alxydim{@C=2.5cm@R=1.5cm}{ {\rm Ad}(\xcP) \ar@{^{(}->}[r]^{j_{{\rm Ad}(\xcP)}} \ar@<.5ex>[d]_{\sfT\;\; } \ar@<-.5ex>[d]^{\;\; \sfS} & {\rm At}(\xcP) \ar@{->>}[r]^{\pi} \ar@<.5ex>[d]_{\sfT \;\;} \ar@<-.5ex>[d]^{\;\; \sfS} & \Si\x\Si{} \ar@<.5ex>[d]_{\pr_1\;\; } \ar@<-.5ex>[d]^{\;\; \:\pr_2} \\ \xcF \ar@{=}[r] & \xcF \ar@{->>}[r]_{\pi_\xcF} \ar@/^2.5pc/@<-.25ex>@{..>}[u]^{\b} & \Si \ar@/^-2.5pc/@<-.25ex>@{..>}[u]_{\unl{\b}} }
\qqq
are commutative shall be called {\bf $\pi$-projectable}.\ Such bisections compose a subgroup in ${\rm Bisec}({\rm At}(\xcP))$,\ which we shall denote henceforth as
\qq\nn
{\rm Bisec}_\pi\bigl({\rm At}(\xcP)\bigr)\,.
\qqq
\exdef

\begin{propanition}\label{prop:bisecAd-restr}
Every $\pi$-projectable bisection $\beta \in {\rm Bisec}_\pi({\rm At}(\xcP))$ such that $\unl\b=\Id\in{\rm Bisec}({\rm Pair}(\Si))$ induces uniquely a bisection $\b_{\rm v}\in{\rm Bisec}({\rm Ad}(\xcP))$.\ The subgroup of all such bisections shall be denoted as 
\qq\nn
{\rm Bisec}_\pi({\rm Ad}(\xcP))\,.
\qqq
\end{propanition}
\beroof
Obvious.
\eroof

\subsection{Automorphisms}\label{sub:autos}

\bedef\label{def:princ-auts}
Let $(\xcP,\Si,\xcG,\pi_\xcP)$ be a principaloid bundle with moment map $\mu$ and action $\varrho$.\ A {\bf principaloid-bundle automorphism of} $\xcP$ is a right-$\xcG$-equivariant bundle automorphism $\Phi\in\Diff_{\varrho(\xcG)}(\xcP)$ which covers a diffeomorphism $f\in\Diff(\Si)$ of the base $\Si$,\ as expressed by the commutative diagram
\qq\nn
\alxydim{@C=1.5cm@R=1.5cm}{ \xcP \ar[r]^{\Phi} \ar[d]_{\pi_\xcP} & \xcP \ar[d]^{\pi_\xcP} \\ \Si \ar[r]_{f} & \Si}\,.
\qqq
Whenever $f=\id_\Si$,\ we call $\Phi$ a {\bf vertical automorphism},\ or a {\bf gauge transformation}.

The {\bf group of} {\bf principaloid-bundle automorphisms} of $\xcP$ shall be denoted as 
\qq\nn
{\rm Aut}(\xcP)\equiv {\rm Aut}_{{\rm {\bf Bun}}(\Si)}(\xcP)\cap\Diff_{\varrho(\xcG)}(\xcP)\,,
\qqq
and its subgroup composed of vertical automorphisms,\ also to be referred to as the {\bf gauge group} of $\xcP$,\ as 
\qq\nn
{\rm Gauge}(\xcP) \equiv {\rm Aut}(\xcP)_{\rm vert}\,.
\qqq
\exdef

\berop\label{prop:P-Auts}
Let $\xcP$ be a principaloid ($\xcG$-)bundle,\ and let $\cO\equiv\{O_i\}_{i\in I}$ be a trivialising cover of its base $\Si$.\ Given $f\in\Diff(\Si)$,\ define a refined cover $\cO^f:=\{O^f_{(j,i)}\equiv f^{-1}(O_j)\cap O_i\}_{(j,i)\in I^{\x 2}_{f,\cO}}$ of $\Si$ with $I^{\x 2}_{f,\cO}:=\{\ (j,i)\in I^{\x 2} \:\vert\:O^f_{(j,i)}\neq\emptyset\ \}$.\ An automorphism $\Phi\in{\rm Aut}(\xcP)$ of $\xcP$ covering $f$ is locally presented by a family of smooth maps
\qq\nn
\g{}_{(j,i)}\colo O^f_{(j,i)}\too\bB\,,\qquad(j,i)\in I^{\x 2}_{f,\cO}
\qqq
as
\qq\label{eq:auts-local-pres}
\Phi\rstr\colo \pi_\xcP^{-1}\bigl(\cO^f_{(j,i)}\bigr)\xrightarrow{\ \cong\ }\pi_\xcP^{-1}\bigl(f\bigl(\cO^f_{(j,i)}\bigr)\bigr),\ \xcP\t_i^{-1}(\si,g)\longmapsto\xcP\t_j^{-1}\bigl(f(\si),L_{\g_{(j,i)}(\si)}(g)\bigr)\,.
\qqq
Over $\cO^f_{(j,i)(l,k)}\equiv O^f_{(j,i)}\cap O^f_{(l,k)}$,\ these maps are subject to the gluing relations 
\qq\label{eq:Auts-loc}
 \gamma{}_{(l,k)}\rstr_{\cO^f_{(j,i)(l,k)}}=\bigl(f^* \b{}_{lj}\cdot \g{}_{(j,i)}\cdot \b{}_{ik}\bigr)\rstr_{\cO^f_{(j,i)(l,k)}}\,.
\qqq

Conversely,\ every such family of maps determines an automorphism of $\xcP$.
\eerop
\beroof 
The local presentation \eqref{eq:auts-local-pres} of the $\xcG$-equivariant map $\Phi$ for the given choice of $\xcG$-equivariant trivialisations $\{\xcP\t_i\}_{i\in I}$ follows directly from Prop.\,\ref{prop:rGequiv-LB}.

For the local maps \eqref{eq:auts-local-pres} to come from a \emph{globally} smooth one $\Phi\in{\rm Aut}_{{\rm {\bf Bun}}(\Si)}(\xcP)$,\ the local presentations must glue over intersections $O^f_{(j,i)(l,k)}\ni\si$ as
\qq\nn
&&\xcP\t_j^{-1}\bigl(f(\si),\g_{(j,i)}(\si)\lact g\bigr)\equiv\Phi_{(j,i)}\bigl(\xcP\t_i^{-1}(\si, g)\bigr)\must\Phi_{(l,k)}\bigl(\xcP\t_i^{-1}(\si, g)\bigr)=\Phi_{(l,k)}\bigl(\xcP\t_k^{-1}\bigl(\si,\b_{ki}(\si)\lact g\bigr)\bigr)\cr\cr
&\equiv&\xcP\t_l^{-1}\bigl(f(\si),\g_{(l,k)}(\si)\cdot\b_{ki}(\si')\lact g\bigr)=\xcP\t_j^{-1}\bigl(f(\si),\b_{jl}\bigl(f(\si)\bigr)\cdot \g_{(l,k)}(\si)\cdot \b_{ki}(\si)\lact g\bigr)\,,
\qqq
whence the statement of the proposition follows upon setting $ g=\Id_m$ for an arbitrary $m\in M$.

The converse claim is proven by reversing the above reasoning in the usual manner.
\eroof

\berop\label{prop:Ups-ind}
There are two canonical group homomorphisms (made explicit in the proof below), 
\qq\nn
\xcF_*\colo {\rm Aut}(\xcP)\too{\rm Aut}(\xcF)\;, \qquad 
\xcA_*\colo {\rm Gauge}(\xcP)\too{\rm Aut}_{\rm vert}({\rm Ad}(\xcP))\,,
\qqq
satisfying,\ for every $\Phi\in{\rm Aut}(\xcP)$ and every $\Phi_{\rm v}\in{\rm Gauge}(\xcP)$,\ the equivariance identities 
\qq
&\xcD\circ\Phi=\xcF_*(\Phi)\circ\xcD\,,&\label{eq:duck-as-intertw}\\ \cr
&\sfS\circ\xcA_*(\Phi_{\rm v})=\xcF_*(\Phi_{\rm v})\circ\sfS\,,\qquad\qquad\xcA_*(\Phi_{\rm v})\circ\sfI=\sfI\circ\xcF_*(\Phi_{\rm v})\,,&\nn
\qqq
where $\sfS$ and $\sfI$ are the structure maps of Diag.\,\eqref{diag:Grmatbndldiag},\ implicitly restricted to ${\rm Ad}(\xcP)\subset{\rm At}(\xcP)$.
\eerop
\beroof
In the light of Prop.\,\ref{prop:P-Auts} (and in the notation introduced in it),\ we find,\ for $\,(\si,g)\in O^f_{(j,i)}\x\xcG$,
\qq\nn
(\xcD\circ\Phi)(\si,g)=\bigl(f(\si),t\bigl(L_{\g_{(j,i)}(\si)}(g)\bigr)\bigr)=\bigl(f(\si),t_*(\g_{(j,i)}(\si))(t(g))\bigr)\,,
\qqq
where the last equality uses Prop.\,\ref{prop:BisAct-vs-str}\,(ii).\ This provides us with the smooth candidate
\qq\nn
\xcF_*(\Phi)_{(j,i)}\colo \pi_\xcF^{-1}\bigl(O^f_{(j,i)}\bigr)\too\pi_\xcF^{-1}\bigl(f\bigl(O^f_{(j,i)}\bigr)\bigr),\ \xcF\t_i^{-1}(\si,m)\longmapsto\xcF\t_j^{-1}\bigl(f(\si),t_*\bigl(\g_{(j,i)}(\si)\bigr)(m)\bigr)
\qqq
for a local presentation of $\xcF_*(\Phi)$.\ We readily check the globality of the thus presented bundle automorphism through a direct calculation,\ carried out for $\si'\in O^f_{(j,i)(l,k)}$,
\qq\nn
&&\xcF_*(\Phi)_{(l,k)}\bigl(\xcF\t_i^{-1}(\si',m)\bigr)=\xcF_*(\Phi)_{(l,k)}\bigl(\xcF\t_k^{-1}\bigl(\si',\b_{ki}(\si')\ulact m\bigr)\bigr)=\xcF\t_l^{-1}\bigl(f(\si'),\g_{(l,k)}(\si')\cdot\b_{ki}(\si')\ulact m\bigr)\cr\cr
&=&\xcF\t_j^{-1}\bigl(\si',f^*\b_{jl}(\si')\cdot \g_{(l,k)}(\si')\cdot\b_{ki}(\si')\ulact m\bigr)=\xcF\t_j^{-1}\bigl(\si',\g_{(j,i)}(\si')\ulact m)\bigr)\bigr)\bigr)\equiv\xcF_*(\Phi)_{(j,i)}\bigl(\xcF\t_i^{-1}(\si',m)\bigr)\,,
\qqq
in which the homomorphic nature of $t_*$ has been taken into account,\ see Def.\,\ref{def:bisec}.\ The latter also ensures homomorphicity of $\xcF_*$.

In the next step,\ we define the obvious local presentation of the induced automorphism
\qq\nn
\xcA_*(\Phi)\colo {\rm Ad}(\xcP)\too{\rm Ad}(\xcP),\ \xcA\t_i^{-1}\bigl(\si, g\bigr)\longmapsto\xcA\t_j^{-1}\bigl(f(\si), C{}_{\gamma_{(j,i)}(\si)}\bigl( g\bigr)\bigr)\,,
\qqq
and verify that it comes from a global map by checking,\ at an arbitrary $\si'\in\cO^f_{(j,i)(l,k)}$,\ the desired equality
\qq\nn
&&\xcA_*(\Phi)\bigl(\xcA\t_k^{-1}\bigl(\si',C_{\b_{ki}(\si')}( g)\bigr)\bigr)\equiv\xcA\t_l^{-1}\bigl(f(\si'), C_{\g_{(l,k)}(\si')}\circ C{}_{\b_{ki}(\si')}(g)\bigr)\cr\cr
&=&\xcA\t_j^{-1}\bigl(f(\si'),C_{\b_{jl}(f(\si'))}\circ C_{\g_{(l,k)}(\si')}\circ C_{\b_{ki}(\si')}(g)\bigr)=\xcA\t_j^{-1}\bigl(f(\si'),C_{f^*\b_{jl}(\si')\cdot \g_{(l,k)}(\si')\cdot \b_{ki}(\si')}(g)\bigr)\cr\cr
&=&\xcA\t_j^{-1}\bigl(f(\si'),C_{\g_{(j,i)}(\si')}(g)\bigr)\equiv \xcA_*(\Phi)\bigl(\xcA\t_i^{-1}(\si', g)\bigr)\,.
\qqq

Finally,\ we verify the anticipated identities
\qq\nn
&&\sfS\circ \xcA_*(\Phi)\bigl(\xcA\t_i^{-1}(\si, g)\bigr)\equiv\xcF\t_j^{-1}\bigl(f(\si),s\circ C_{\g_{(j,i)}(\si)}(g)\bigr)=\xcF\t_j^{-1}\bigl(f(\si),t_*\bigl(\g_{(j,i)}(\si)\bigr)\bigl(s(g)\bigr)\bigr)\cr\cr
&\equiv&\xcF_*(\Phi)\bigl(\xcF\t_i^{-1}\bigl(\si,s(g)\bigr)\bigr)\equiv\xcF_*(\Phi)\circ\sfS\bigl(\xcA\t_i^{-1}\bigl(\si, g\bigr)\bigr)\,,
\qqq
see Cor.\,\ref{cor:BisCon-vs-str}\,(i),\ and 
\qq\nn
&&\xcA_*(\Phi)\circ\sfI\bigl(\xcF\t_i^{-1}(\si,m)\bigr)\equiv\xcA\t_j^{-1}\bigl(f(\si),C_{\g_{(j,i)}(\si)}(\Id_m)\bigr)=\xcA\t_j^{-1}\bigl(f(\si),\Id_{t_*(\g_{(j,i)}(\si))(m)}\bigr)\cr\cr
&\equiv&\sfI\circ\xcF\t_j^{-1}\bigl(f(\si),t_*\bigl(\g_{(j,i)}(\si)\bigr)(m)\bigr)=\sfI\circ\xcF_*(\Phi)\bigl(\xcF\t_i^{-1}(\si,m)\bigr)\,,
\qqq
see Prop.\,\ref{prop:BisAct-vs-str}\,(iii). 
\eroof

\bethe\label{thm:autP-from-BisAt}
There is a canonical group isomorphism (made explicit in the proof below)
\qq\nn
{\rm Bisec}_\pi\bigl({\rm At}(\xcP)\bigr)\cong{\rm Aut}(\xcP)\,,
\qqq
which restricts to a group isomorphism 
\qq\nn
{\rm Bisec}_\pi({\rm Ad}(\xcP))\cong{\rm Gauge}(\xcP)
\qqq
for the respective subgroups.
\ethe
\beroof
Consider an automorphism $\Phi\in{\rm Aut}(\xcP)$ which covers $f\in\Diff(\Si)$ and is presented by local data $\{\g_{(j,i)}\colo O^f_{(j,i)}\to\bB\}_{(j,i)\in I^{\x 2}_{f,\cO}}$.\ These can be used to define a map
\qq\nn
\b_\Phi\colo \xcF\too{\rm At}(\xcP),\ [(\si,m,i)]\longmapsto[(f(\si),\g_{(j,i)}(\si)(m),\si,j,i)]\,,
\qqq
written,\ in the `resolved' model \eqref{eq:AtP-loc-res} for ${\rm At}(\xcP)$,\ and for an \emph{arbitrary} index $j\in I$ with the property $O^f_{(j,i)}\ni\si$ (such an index exists,\ because the $O^f_{(j,i)}$ compose a cover $\Si$,\ which refines $\cO$).\ The map $\b_\Phi$ is readily checked to be globally smooth:\ Indeed,\ for every $\si'\in O^f_{(j,i)}\cap O^f_{(l,k)}$,\ we obtain the identity
\qq\nn
&&\b_\Phi\bigl([(\si',\b_{ki}(\si)\ulact m,k)]\bigr)\equiv[(f(\si),\g_{(l,k)}(\si)(\b_{ki}(\si)\ulact m),\si,l,k)]\cr\cr
&=&[(f(\si),(f^*\b_{lj}(\si)\cdot\g_{(j,i)}(\si))((t\cdot\b_{ik}(\si))(\b_{ki}(\si)\ulact m)).(\b_{ki}(\si)((t\circ\b_{ki}(\si))^{-1}(\b_{ki}(\si)\ulact m))^{-1},\si,l,k)]\cr\cr
&=&[(f(\si),f^*\b_{lj}(\si)(\g_{(j,i)}(\si)\ulact m).\g_{(j,i)}(\si)(m).(\b_{ki}(\si)(m))^{-1},\si,l,k)]\cr\cr
&\equiv&[(f(\si),f^*\b_{lj}(\si)\lact \g_{(j,i)}(\si)(m)\ract \b_{ik}(\si),\si,l,k)]=[(f(\si),\g_{(j,i)}(\si)(m),\si,j,i)]\,,
\qqq
in which the second equality is a simple rewrite of \Reqref{eq:Auts-loc}.

Next,\ we verify that $\b_\Phi$ is a bisection of the Atiyah--Ehresmann groupoid of $\xcP$.\ To this end,\ we compute
\qq\nn
\sfS\circ\b_\Phi\bigl([(\si,m,i)]\bigr)\equiv\sfS\bigl([(f(\si),\g_{(j,i)}(\si)(m),\si,j,i)]\bigr)\equiv[(\si,(s\circ(\g_{(j,i)}(\si)))(m),i)]=[(\si,m,i)]\,,
\qqq
using $\g_{(j,i)}(\si)\in\bB$ along the way,\ and
\qq\nn
&&\sfT\circ\b_\Phi\bigl([(\si,m,i)]\bigr)\equiv\sfT\bigl([(f(\si),\g_{(j,i)}(\si)(m),\si,j,i)]\bigr)\equiv[(f(\si),(t\circ(\g_{(j,i)}(\si)))(m),j)]\cr\cr
&\equiv&\xcF_*(\Phi)\bigl([(\si,m,i)]\bigr)\,,
\qqq
see the local presentation $\xcF_*(\Phi)_{(j,i)}$ of $\xcF_*(\Phi)$ in the proof of Prop.\,\ref{prop:Ups-ind}.\ Thus,\ altogether
\qq\label{eq:betaPh-as-bisec}
\sfS\circ\b_\Phi=\id_\xcF\,,\qquad\qquad\sfT\circ\b_\Phi=\xcF_*(\Phi)\in{\rm Aut}(\xcF)\,,
\qqq
which confirms our claim.\ The $\pi$-projectability of $\b_\Phi$ is manifest.

Conversely,\ given a $\pi$-projectable bisection $\b\in{\rm Bisec}_\pi({\rm At}(\xcP))$,\ we associate with it an automorphism $\Phi_\b\in{\rm Aut}(\xcP)$ as follows.\ Consider a local presentation of the bundle map $\b\colo\xcF\to{\rm At}(\xcP)$ relative to the `resolved' model \eqref{eq:AtP-loc-res} for ${\rm At}(\xcP)$,
\qq\label{eq:bisecAt-loc}
\b\colo\xcF\too{\rm At}(\xcP),\ [(\si,m,i)]\longmapsto[(f(\si),b_{(j,i)}(\si,m),\si,j,i)]\,.
\qqq
Here,\ we have already taken into account the assumption that $\b$ covers a bisection $(f,\id_\Si)$ of ${\rm Pair}(\Si)$,\ for some $f\in\Diff(\Si)$,\ and $j\in I$ is an (arbitrary) index such that $O^f_{(j,i)}\ni\si$.\ For the above presentation to be independent of the arbitrary choices (of indices) made,\ we need to demand that the identity 
\qq\label{eq:bisecAt-glue}
b_{(l,k)}\bigl(\si',\b_{ki}(\si')\ulact m\bigr)=f^*\b_{lj}(\si')\lact b_{(j,i)}(\si',m)\ract\b_{ik}(\si')
\qqq
hold true for every $m\in M$ and $\si'\in O^f_{(j,i)}\cap O^f_{(l,k)}$.\ Define
\qq\nn
\g_{(j,i)}(\si):=b_{(j,i)}(\si,\cdot)\colo M\too\xcG\,,
\qqq
and invoke the local presentation of the structure maps $\sfS$ and $\sfT$ to read off from \eqref{eq:bisecAt-loc} the conditions
\qq\nn
s\circ\g_{(j,i)}(\si)=\id_M\,,\qquad\qquad t\circ\g_{(j,i)}(\si)\in\Diff(M)\,.
\qqq
These jointly imply
\qq\nn
\g_{(j,i)}\colo O^f_{(j,i)}\too\bB\,.
\qqq
Upon rewriting \Reqref{eq:bisecAt-glue} in terms of the $\g_{(j,i)}$ and moving the transition maps of $\xcP$ to the other side of the equality sign,\ we obtain the gluing law
\qq\label{eq:Auts-loc-atm}
\g_{(j,i)}(\si')(m)=f^*\b_{jl}(\si')\lact\g_{(l,k)}(\si')\bigl(\b_{ki}(\si')\ulact m\bigr)\ract\b_{ki}(\si')\,,
\qqq
in which we recognise the familiar formula \eqref{eq:Auts-loc},\ evaluated at $m$.\ Hence,\ by Prop.\,\ref{prop:P-Auts},\ $\b$ does,\ indeed,\ determine an automorphism $\Phi_\b\in{\rm Aut}(\xcP)$ with local data $\g_{(j,i)}$ which covers $f\in\Diff(\Si)$.\ We have the obvious identity $\b\equiv\b_{\Phi_\b}$.

The last statement of the theorem follows straightforwardly from Prop.\,\ref{prop:bisecAd-restr}.
\eroof

We conclude the present section with a natural definition of automorphisms of the $\xcG$-fibred $\bB$-bundles of Def.\,\ref{def:Bequiv-Gbndl},\ which provides us with an alternative characterisation of ${\rm Aut}(\xcP)$ in the bisection-complete case.
\bedef\label{def:Bequiv-Gbndl-auts}
An {\bf automorphism} of a $\xcG$-fibred $\bB$-bundle $\widehat\xcP$ is an arbitrary $\bB$-equivariant bundle automorphism of $\widehat\xcP$.
\exdef
\berop
Let $\xcG$ be a bisection-complete Lie groupoid,\ and let $\widehat\xcP$ be a $\xcG$-fibred $\bB$-bundle.\ Every automorphism of $\widehat\xcP$ is locally realised by elements of $L(\bB)$ as in Prop.\,\ref{prop:P-Auts}.\ Conversely,\ every automorphism of $\widehat\xcP$,\ viewed as a pricipaloid bundle,\ is $\bB$-equivariant with respect to the induced $\bB$-action \eqref{eq:R-act-Poid}. 
\eerop
\beroof
The first part of the proposition follows directly from Prop.\,\ref{prop:RB-comm-LB},\ and the second one is a consequence of the fact that left-multiplication commutes with right-multiplication. 
\eroof

\brem
Our discussion for bisection-complete structure groupoids is summarised by the following three-floor commutative diagram
\qq\nn\hspace{-3.75cm}
\alxydim{@C=1.6cm@R=1.5cm}{ & & & \xcP\x_\Si\xcP{}^{\circlearrowleft\bB} \ar@{^{(}->}[rr] \ar@<.5ex>[dl]^{\pr_1} \ar@<-.5ex>[dl]_{\pr_2} \ar@{->>}'[d]^{\xcD\x\xcD}[dd] & & \xcP\x\xcP{}^{\circlearrowleft\bB} \ar@{->>}[rr]^{\pi_\xcP\x\pi_\xcP} \ar@<.5ex>[dl]^{\pr_1} \ar@<-.5ex>[dl]_{\pr_2} \ar@{->>}'[d]^{\xcD\x\xcD}[dd] & & \Si\x\Si{}^{\circlearrowleft\bB} \ar@<.5ex>[dl]^{\pr_1} \ar@<-.5ex>[dl]_{\pr_2} \ar@{=}[dd] & & \\ & & \xcP{}^{\circlearrowleft\bB} \ar@{=}[rr] \ar@{->>}[dd]^(.3){\xcD} & & \xcP{}^{\circlearrowleft\bB} \ar@{->>}[rr]^(.7){\pi_\xcP} \ar@{->>}[dd]^(.3){\xcD}|(0.62)\hole \ar@/^2pc/@<-.25ex>[ur]^(.35){(\Phi,\id_\xcP)} & & \Si{}^{\circlearrowleft\bB} \ar@{=}[dd] \ar@/^-2pc/@<.25ex>[ur]_(.35){\unl\b{}_\Phi\equiv(f,\id_\Si)} & & \\ & & & \xcF\x_\Si\xcF \ar@{^{(}->}'[r][rr] \ar@<.5ex>[dl]^{\pr_1} \ar@<-.5ex>[dl]_{\pr_2} & & \xcF\x\xcF \ar@{->>}'[r]^{\qquad\pi_\xcF\x\pi_\xcF}[rr] \ar@<.5ex>[dl]^{\pr_1} \ar@<-.5ex>[dl]_{\pr_2} & & \Si\x\Si \ar@<.5ex>[dl]^{\pr_1} \ar@<-.5ex>[dl]_{\pr_2} \ar@{=}[dd] & & \\ & & \xcF \ar@{=}[rr] \ar@{=}[dd] & & \xcF \ar@{->>}[rr]^(.75){\pi_\xcF} \ar@{=}[dd] \ar@/^2pc/@<-.25ex>[ur]^(.35){(\xcF_*(\Phi),\id_\xcF)} & & \Si \ar@{=}[dd] \ar@/^-2pc/@<.25ex>[ur]_(.35){\unl\b{}_\Phi} & & \\ & & & {\rm Ad}(\xcP) \ar@{^{(}->}'[r][rr] \ar@<.5ex>[dl]^{\sfT} \ar@<-.5ex>[dl]_{\sfS} \ar@{->>}'[u]_{(\txT,\txS)}[uu] & & {\rm At}(\xcP) \ar@{->>}'[r]^{\quad\pi}[rr] \ar@<.5ex>[dl]^{\sfT} \ar@<-.5ex>[dl]_{\sfS} \ar@{->>}'[u]_{(\sfT,\sfS)}[uu] & & \Si\x\Si \ar@<.5ex>[dl]^{\pr_1} \ar@<-.5ex>[dl]_{\pr_2} & & \\ & & \xcF \ar@{=}[rr] & & \xcF \ar[rr]^{\pi_\xcF} \ar@/^2pc/@<-.25ex>[ur]^(.35){\b_\Phi} & & \Si \ar@/^-2pc/@<.25ex>[ur]_(.35){\unl\b{}_\Phi} & & }\,,
\qqq
which we may call the {\bf extended Atiyah sequence}.\ It shows the interplay between $\bB$-equivariant automorphisms of the bundle $\xcP$ on the top floor,\ and the corresponding bisections of the groupoids on the bottom floor,\ where $\bB$-equivariance is encoded by their construction.\ Note also that while the middle floor is redundant for ordinary principal $\txG$-bundles,\ here,\ it carries information on the induction of automorphisms of the shadow bundle from those of the principaloid bundle.

Whenever the underlying bisection $\unl\b{}_\Phi$ of ${\rm Pair}(\Si)$ is the identity bisection $\Id(\Si)$,\ with $f=\id_\Si$,\ then $\b_\Phi$ induces a bisection of the adjoint groupoid ${\rm Ad}(\xcP)$,\ which then corresponds to a gauge transformation of $(\xcP,\xcF)$.
\erem

\subsection{The trident of the principaloid bundle}

\bethe\label{thm:princ-as-At-mod} 
The Atiyah--Ehresmann groupoid acts on the principaloid bundle and its shadow in the following way:
\bit
\item[$\bullet$] On every principaloid bundle $\xcP$,\ there exists a canonical structure of a left ${\rm At}(\xcP)$-module,\ with momentum $\mu_\xcP=\xcD\colo\xcP\to\xcF$ and action
\qq\label{eq:lamPa}\qquad\qquad
\la_\xcP\colo {\rm At}(\xcP){}_{\sfS}\hspace{-2pt}\x_{\xcD}\hspace{-2pt}\xcP\too\xcP,\ \bigl([(\si_1,g,\si_2,i,j)],[(\si_2,h,k)]\bigr)\longmapsto[(\si_1,g.(\b_{jk}(\si_2)\lact h),i)]\,,
\qqq
written in the local models:\ \eqref{eq:AtP-loc-res} for ${\rm At}(\xcP)$ and \eqref{eq:clutch-xcP} for $\xcP$,\ for $\si_1\in O_i$ and $\si_2\in O_{jk}$.\
\item[$\bullet$] This structure covers the canonical structure of a left ${\rm Pair}(\Si)$-module on $\Si$,\ with momentum $\mu_\Si\equiv\id_\Si$ and action
\qq\nn
\la_\Si\colo(\Si\x\Si)\,{}_{\pr_2}\hspace{-3pt}\x_{\id_\Si}\hspace{-2pt}\Si\too\Si,\ \bigl((\si_1,\si_2),\si_2\bigr)\longmapsto\si_1\,.
\qqq
\item[$\bullet$] On the corresponding shadow bundle $\xcF$,\ there exists a canonical structure of a left ${\rm At}(\xcP)$-module,\ with momentum $\mu_\xcF=\id_\xcF$ and action 
\qq\label{eq:lamFa}
\la_\xcF\colo{\rm At}(\xcP){}_{\sfS}\hspace{-2pt}\x_{\id_\xcF}\hspace{-2pt}\xcF\too\xcF,\ \bigl([(\si_1,g,\si_2,i,j)],[(\si_2,m,k)]\bigr)\longmapsto[(\si_1,t(g),i)]\,.
\qqq
\item[$\bullet$] The action $\la_\xcP$ is intertwined with $\la_\xcF$ by the sitting-duck map $\xcD$,\ as reflected in the identities
\qq\label{eq:lamPa-int-lamFa}
\mu_\xcP=\mu_\xcF\circ\xcD\,,\qquad\qquad\la_\xcF\circ\bigl(\id_{{\rm At}(\xcP)}\x\xcD\bigr)=\xcD\circ\la_\xcP\,.
\qqq
\eit
\ethe
\beroof
Let us first analyse the triple $(\xcP,\mu_\xcP,\la_\xcP)$.\ Composability of the two arrows:\ $g$ and $\b_{jk}(\si_2)\lact h$ in \eqref{eq:lamPa} is ensured by the identity
\qq\nn
[(\si_2,s(g),j)]\equiv\sfS\bigl([(\si_1,g,\si_2,i,j)]\bigr)=\xcD\bigl([(\si_2,h,k)]\bigr)\equiv[(\si_2,t(h),k)]=[(\si_2,t_*(\b_{jk}(\si_2))(t(h)),j)]\,.
\qqq
Using the 1-cocycle condition for the $\b_{ij}$ alongside identities \eqref{eq:BisAct-vs-str-vi},\ we readily establish,\ for any $\si_1\in O_{il}$ and $\si_2\in O_{jkmn}$,\ the identity 
\qq\nn
&&\la_\xcP\bigl([(\si_1,\b_{li}(\si_1)\lact g\ract\b_{jm}(\si_2),\si_2,l,m)],[(\si_2,\b_{nk}(\si_2)\lact h,n)]\bigr)\cr\cr
&\equiv&[(\si_1,(\b_{li}(\si_1)\lact g\ract\b_{jm}(\si_2)).(\b_{mn}(\si_2)\lact(\b_{nk}(\si_2)\lact h)),l)]\cr\cr
&\equiv&[(\si_1,(\b_{li}(\si_1)\lact g\ract\b_{jm}(\si_2)).(\b_{mk}(\si_2)\lact h),l)]\cr\cr
&=&[(\si_1,(\b_{li}(\si_1)\lact g).(\b_{jm}(\si_2)\lact(\b_{mk}(\si_2)\lact h)),l)]=[(\si_1,(\b_{li}(\si_1)\lact g).(\b_{jk}(\si_2)\lact h),l)]\cr\cr
&=&[(\si_1,\b_{li}(\si_1)\lact(g.(\b_{jk}(\si_2)\lact h)),l)]=[(\si_1,g.(\b_{jk}(\si_2)\lact h),i)]\,,
\qqq
which demonstrates that $\la_\xcP$ is a well-defined globally smooth action map sought after.\ Its constitutive properties are readily verified in the following calculations:\ Axiom (GlM1) of Def.\,\ref{def:gr-mod} is checked in
\qq\nn
&&(\mu_\xcP\circ\la_\xcP)\bigl([(\si_1,g,\si_2,i,j)],[(\si_2,h,k)]\bigr)\equiv\xcD\bigl([(\si_1,g.(\b_{jk}(\si_2)\lact h),i)]\bigr)=[(\si_1,t(g.(\b_{jk}(\si_2)\lact h)),i)]\cr\cr
&=&[(\si_1,t(g),i)]\equiv \sfT\bigl([(\si_1,g,\si_2,i,j)]\bigr)\,.
\qqq
Axiom (GlM2) is seen to hold true in virtue of
\qq\nn
&&\la_\xcP\bigl(\sfI\circ\mu_\xcP\bigl([(\si,g,i)]\bigr),[(\si,g,i)]\bigr)\equiv\la_\xcP\bigl([(\si,\Id_{t(g)},\si,i,i)],[(\si,g,i)]\bigr)\equiv[(\si,\Id_{t(g)}.(\b_{ii}(\si)\lact g),i)]\cr\cr
&=&[(\si,g,i)]\,.
\qqq
Axiom (GlM3) is established,\ with the help of the first of identities \eqref{eq:BisAct-vs-str-v},\ in
\qq\nn
&&\la_\xcP\bigl([(\si_1,g_1,\si_2,i,j)],\la_\xcP\bigl([(\si_2,g_2,\si_3,k,l)],[(\si_3,h,m)]\bigr)\bigr)\cr\cr
&=&\la_\xcP\bigl([(\si_1,g_1,\si_2,i,j)],[(\si_2,g_2.(\b_{lm}(\si_3)\lact h),k)]\bigr)=[(\si_1,g_1.(\b_{jk}(\si_2)\lact(g_2.(\b_{lm}(\si_3)\lact h))),i)]\cr\cr
&=&[(\si_1,(g_1.(\b_{jk}(\si_2)\lact g_2)).(\b_{lm}(\si_3)\lact h))),i)]\equiv\la_\xcP\bigl([(\si_1,g_1.(\b_{jk}(\si_2)\lact g_2),\si_3,i,l)],[(\si_3,h,m)]\bigr)\cr\cr
&\equiv&\la_\xcP\bigl(\sfM\bigl([(\si_1,g_1,\si_2,i,j)],[(\si_2,g_2,\si_3,k,l)]\bigr),[(\si_3,h,m)]\bigr)\,.
\qqq
Projectability of the above structure of a left ${\rm At}(\xcP)$-module on $\xcP$ onto the canonical structure of a left ${\rm Pair}(\Si)$-module on $\Si$ (along $(\pi,\pi_\xcF)$) is self-evident. 

Passing to the triple $(\xcF,\mu_\xcF,\la_\xcF)$,\ we note well-definedness of $\la_\xcF$ and manifest projectability of the entire structure onto $(\Si,\id_\Si,\la_\Si)$.\ The first of identities \eqref{eq:lamPa-int-lamFa} is trivially satisfied,\ and so we are left with the second one to prove.\ First of all,\ note the obvious identity
\qq\nn
\bigl(\id_{{\rm At}(\xcP)}\x\xcD\bigr)\bigl({\rm At}(\xcP){}_{\sfS}\hspace{-2pt}\x_{\xcD}\hspace{-2pt}\xcP\bigr)\equiv{\rm At}(\xcP){}_{\sfS}\hspace{-2pt}\x_{\id_\xcF}\hspace{-2pt}\xcF\,,
\qqq
which ensures meaningfulness of that ${\rm At}(\xcP)$-equivariance condition.\ The latter is checked in a direct computation:
\qq\nn
&&\bigl(\la_\xcF\circ\bigl(\id_{{\rm At}(\xcP)}\x\xcD\bigr)\bigr)\bigl([(\si_1,g,\si_2,i,j)],[(\si_2,h,k)]\bigr)\equiv\la_\xcF\bigl([(\si_1,g,\si_2,i,j)],[(\si_2,t(h),k)]\bigr)\equiv[(\si_1,t(g),i)]\cr\cr
&=&[(\si_1,t(g.(\b_{jk}(\si_2)\lact h)),i)]\equiv\xcD\bigl([(\si_1,g.(\b_{jk}(\si_2)\lact h),i)]\bigr)\equiv\bigl(\xcD\circ\la_\xcP\bigr)\bigl([(\si_1,g,\si_2,i,j)],[(\si_2,h,k)]\bigr)\,.
\qqq
\eroof

\brem
Note that the independence of $\la_\xcP$ of the choice of local charts permits us to simplify \eqref{eq:lamPa} as
\qq\label{eq:lamPa-simpl}
\la_\xcP\colo {\rm At}(\xcP){}_{\sfS}\hspace{-2pt}\x_{\xcD}\hspace{-2pt}\xcP\too\xcP,\ \bigl([(\si_1,g,\si_2,i,j)],[(\si_2,h,j)]\bigr)\longmapsto[(\si_1,g.h,i)]\,.
\qqq
\erem
 
In analogy with the construction of Remark \ref{rem:Bisec-geom-act},\ and consistently with Thm.\,\ref{thm:autP-from-BisAt},\ we establish 
\berop
The structure of a left ${\rm At}(\xcP)$-module on $\xcP$ gives rise to a geometric implementation of automorphisms ${\rm Aut}(\xcP)\ni\Phi$ on $\xcP$ (resp.\ on $\xcF$) through restriction of $\la_\xcP$ (resp.\ of $\la_\xcF$) to the corresponding submanifolds $\b_\Phi(\xcF)\subset{\rm At}(\xcP)$,\ embedded by the $\pi$-projectable bisections $\b_\Phi\in{\rm Bisec}_\pi({\rm At}(\xcP))$,
\qq\nn
\alxydim{@C=1.cm@R=1.cm}{ & & & \xcP \ar[ddl]_{\mu_\xcP\equiv\xcD} \ar@{->>}[dr]^{\xcD} & \\ {\rm At}(\xcP) \ar@{=>}[drr] \ar@{->>}[dd]_{\pi} & \b_\Phi(\xcF) \ar@{_{(}->}[l] & & & \xcF \ar[dll]^{\mu_\xcF\equiv\id_\xcF} \ar@{->>}[dd]^{\pi_\xcF} \\ & & \xcF \ar@{->>}[dd]^{\pi_\xcF} \ar[ul]_{\b_\Phi} & & \\ \Si\x\Si \ar@{=>}[drr] & (f,\id_\Si)(\Si) \ar@{_{(}->}[l] & & & \Si \ar[dll]^{\mu_\Si} \\ & & \Si \ar[ul]_{(f,\id_\Si)} & & }\,.
\qqq
This is captured by the following commutative diagrams:
\qq\nn
\alxydim{@C=3.cm@R=1.5cm}{ \xcF{}_{\id_\xcF}\hspace{-3pt}\x_{\mu_\xcP}\hspace{-1pt}\xcP \ar[r]^{\quad\b_\Phi\x\id_\xcP} & {\rm At}(\xcP){}_{\sfS}\hspace{-2pt}\x_{\xcD}\hspace{-2pt}\xcP \ar[d]^{\la_\xcP} \\ \xcP \ar[u]^{(\mu_\xcP,\id_\xcP)} \ar[r]_{\Phi} & \xcP }\,,
\qqq
and
\qq\nn
\alxydim{@C=3.cm@R=1.5cm}{ & {\rm At}(\xcP){}_{\sfS}\hspace{-2pt}\x_{\id_\xcF}\hspace{-2pt}\xcF \ar[d]^{\la_\xcF} \\ \xcF \ar[ur]^{\quad\b_\Phi\x\id_\xcF} \ar[r]_{\xcF_*(\Phi)} & \xcF }\,.
\qqq
\eerop
\beroof
Given an automorphism $\Phi$ covering $f\in\Diff(\Si)$,\ with local data $\{\g_{(j,i)}\colo O^f_{(j,i)}\to\bB\}_{(j,i)\in I^{\x 2}_{f,\cO}}$,\ we obtain
\qq\nn
&&\la_\xcP\bigl(\b_\Phi\bigl([(\si,t(\b_{ik}(\si)\lact g),i)]\bigr),[(\si,g,k)]\bigr)=\la_\xcP\bigl([(f(\si),\g_{(j,i)}(\si)(\b_{ik}(\si)\ulact t(g)),\si,j,i)],[(\si,g,k)]\bigr)\cr\cr
&=&[(f(\si),\g_{(j,i)}(\si)(\b_{ik}(\si)\ulact t(g)).(\b_{ik}(\si)\lact g),j)]=[(f(\si),(\g_{(j,i)}(\si)(\b_{ik}(\si)\ulact t(g))\ract\b_{ik}(\si)).g,j)]\cr\cr
&=&[(f(\si),\g_{(j,k)}(\si)(t(g)).g,j)]\equiv[(f(\si),\g_{(j,k)}(\si)\lact g,j)]\equiv\Phi\bigl([(\si,g,k)]\bigr)\,.
\qqq
The well-definedness of the above calculation is ensured by the identity $\sfS\circ\b_\Phi\circ\mu_\xcP=\mu_\xcP$ (see \eqref{eq:betaPh-as-bisec}),\ which implies composability of the arrows in the third expression.\ The third equality is implied by the second identity in \Reqref{eq:BisAct-vs-str-vi},\ and the next one follows upon invoking the explicit form \eqref{eq:Auts-loc-atm} of the gluing law \eqref{eq:Auts-loc}.\ Thus,\ altogether,\ we arrive at the anticipated identity
\qq\nn
\la_\xcP\circ\bigl(\b_\Phi\x\id_\xcP\bigr)\circ\bigl(\mu_\xcP,\id_\xcP\bigr)=\Phi\,.
\qqq
The corresponding identity for $\la_\xcF$ is now readily derived with the help of identities \eqref{eq:lamPa-int-lamFa} and the first one of the equivariance identities of Prop.\,\ref{prop:Ups-ind}:
\qq\nn
&&\la_\xcF\circ\bigl(\b_\Phi\x\id_\xcF\bigr)\circ\xcD\equiv\la_\xcF\circ\bigl(\b_\Phi\x\id_\xcF\bigr)\circ\bigl(\mu_\xcF,\id_\xcF\bigr)\circ\xcD=\la_\xcF\circ\bigl(\b_\Phi\x\id_\xcF\bigr)\circ\bigl(\mu_\xcP,\xcD\bigr)\cr\cr
&\equiv&\la_\xcF\circ\bigl(\id_{{\rm At}(\xcP)}\x\xcD\bigr)\circ\bigl(\b_\Phi\x\id_\xcP\bigr)\circ\bigl(\mu_\xcP,\id_\xcP\bigr)=\xcD\circ\la_\xcP\circ\bigl(\b_\Phi\x\id_\xcP\bigr)\circ\bigl(\mu_\xcP,\id_\xcP\bigr)=\xcD\circ\Phi\cr\cr
&=&\xcF_*(\Phi)\circ\xcD\,,
\qqq
and taking into account the surjectivity of $\xcD$.
\eroof

\bedef
Let $\grpd{\xcG_A}{M_A},\ A\in\{1,2,3\}$ be Lie groupoids,\ let $\widehat P$ be a principal $(\xcG_1,\xcG_2)$-bibundle,\ as in Def.\,\ref{def:bibndl},\ and let $\Si$ be a manifold.\ We shall call the quintuple $(\xcG_1,\widehat P,\xcG_2;\Si,\xcG_3)$ a ({\bf left}) {\bf trident} with {\bf base} $\Si$ and fibre $\xcG_3$ if the following conditions are satisfied:
\bit
\item $\widehat P$ is the total space of a fibre bundle $\pi_{\widehat P}\colo\widehat P\to\Si$ with base $\Si$ and typical fibre $\xcG_3$;
\item the right $\xcG_2$-action $\rho_2\equiv\mact$ preserves $\pi_{\widehat P}$-fibres,\ {\it i.e.},\ $\pi_{\widehat P}(p\mact g_2)=\pi_{\widehat P}(p)$ for all $(p,g_2)\in\widehat P{}_{\mu_2}\hspace{-3pt}\x_{t_2}\hspace{-1pt}\xcG_2$ ;
\item the Lie groupoid $\grpd{\xcG_1}{M_1}$ is the total space of a fibre-bundle object in the category of Lie groupoids,\ with base ${\rm Pair}(\Si)$ and typical fibre $\grpd{\xcG_3}{M_3}$
\item the $\xcG_1$-module structure on $\widehat P$ covers the canonical left ${\rm Pair}(\Si)$-module structure on $\Si$,\ and is modelled on the canonical left $\xcG_3$-module structure on $\xcG_3$ (in a local trivialisation),\ as captured by the diagram:
\qq\nn
\alxydim{@C=.75cm@R=.5cm}{ \xcG_3 \ar@{=>}[rd] \ar@{^{(}.>}[dd] & & \xcG_3 \ar[dl]^{t_3} \ar@{^{(}.>}[dd] \\ & M_3 \ar@{^{(}.>}[dd] & & \\ \xcG_1 \ar@{=>}[rd] \ar[dd]_{\pi_{\xcG_1}} & & \widehat P \ar[dl]^{\mu_1} \ar[dd]^{\pi_{\widehat P}} \\ & M_1 \ar[dd]_{\pi_{M_1}} & & \\ \Si\x\Si \ar@{=>}[rd] & & \Si \ar[dl]^{\id_\Si} \\ & \Si & }\,.
\qqq
\eit
The trident shall be represented by the following diagram:
\qq\label{diag:trident}
\alxydim{@C=.75cm@R=1cm}{ & & \xcG_3 \ar@{^{(}.>}[d] & & \\ \xcG_1 \ar@{=>}[rd] & & \widehat P \ar[dl]_{\mu_1} \ar[rd]^{\mu_2} \ar[dd]^{\pi_{\widehat P}} & & \xcG_2 \ar@{=>}[ld] \\ & M_1 & & M_2 & \\ & & \Si & &}\,.
\qqq 
\exdef

\bethe\label{thm:W-str-on-P}
For every principaloid bundle $\xcP$,\ the quintuple $({\rm At}(\xcP),\xcP,\xcG;\Si,\xcG)$ is a trident,\ captured by the following diagram
\qq\nn
\alxydim{@C=.75cm@R=1cm}{ & & \xcG \ar@{^{(}.>}[d] & & \\ {\rm At}(\xcP) \ar@{=>}[rd] & & \xcP \ar[dl]_{\xcD} \ar[rd]^{\mu} \ar[dd]^{\pi_\xcP} & & \xcG \ar@{=>}[ld] \\ & \xcF & & M & \\ & & \Si & &}\,.
\qqq
\ethe
\beroof
The structure of a right principal $\xcG$-module on $(\xcP,\xcF,\xcD)$ was established in Thm.\,\ref{thm:duck-as-prince}.\ The underlying right $\xcG$-action on $\xcP$ preserves $\pi_\xcP$-fibres by definition.\ Furthermore,\ the identification of ${\rm At}(\xcP)$ as a fibre-bundle object in the category of Lie groupoids,\ with base ${\rm Pair}(\Si)$ and typical fibre $\xcG$,\ was made in Remark \ref{rem:Atiyah-as-bndl}.\ Therefore,\ it remains to investigate the left action of ${\rm At}(\xcP)$ on $\xcP$,\ as identified in Thm.\,\ref{thm:princ-as-At-mod}. 

We begin by noting that $\la_\xcP\colo {\rm At}(\xcP){}_{\sfS}\hspace{-2pt}\x_{\xcD}\hspace{-2pt}\xcP\too\xcP$ preserves $\mu$-fibres.\ Indeed,\ the action restricts to fibres of $\xcP$ as the \emph{left}-multiplication,\ and so leaves the moment map $\mu$,\ locally modelled on $s$,\ invariant. 

The action gives rise to a smooth map
\qq\nn
(\la_\xcP,\pr_2)\colo {\rm At}(\xcP){}_{\sfS}\hspace{-2pt}\x_{\xcD}\hspace{-2pt}\xcP\too\xcP{}_{\mu}\hspace{-2pt}\x_{\mu}\hspace{-2pt}\xcP\,,
\qqq
whose well-definedness is ensured by the following equality (see \Reqref{eq:lamPa-simpl}):
\qq\nn
\mu\bigl([(\si_1,g.h,i)]\bigr)\equiv s(g.h)=s(h)\equiv\mu\bigl([(\si_2,h,j)]\bigr)\,.
\qqq
Consider the division map 
\qq\nn
\psi_\xcP\colo\xcP{}_{\mu}\hspace{-2pt}\x_{\mu}\hspace{-2pt}\xcP\too{\rm At}(\xcP),\ \bigl([(\si_1,g,i)],[(\si_2,h,j)]\bigr)\longmapsto[(\si_1,g.h^{-1},\si_2,i,j)]\,.
\qqq
It is well-defined since $s(g)\equiv\mu([(\si_1,g,i)])=\mu([(\si_2,h,j)])\equiv s(h)=t(h^{-1})$.\ It now remains to verify the identity
\qq\label{eq:psiP-as-laPinv}
(\psi_\xcP,\pr_2)=(\la_\xcP,\pr_2)^{-1}\,.
\qqq
We have
\qq\nn
\sfS\bigl([(\si_1,g.h^{-1},\si_2,i,j)]\bigr)\equiv [\si_2,s(g.h^{-1}),j]=[\si_2,t(h),j]\equiv\xcD\bigl([(\si_2,h,j)]\bigr)\,,
\qqq
so that $(\psi_\xcP,\pr_2)$ maps $\xcP{}_{\mu}\hspace{-2pt}\x_{\mu}\hspace{-2pt}\xcP$ to ${\rm At}(\xcP){}_{\sfS}\hspace{-2pt}\x_{\xcD}\hspace{-2pt}\xcP$.\ Identity \eqref{eq:psiP-as-laPinv} can now be checked in a direct calculation.\ This establishes $\la_\xcP$ as a principal action.

That $\la_\xcP$ covers the canonical left action of ${\rm Pair}(\Si)$ on $\Si$ is part of Thm.\,\ref{thm:princ-as-At-mod}. 
\eroof

\brem
In the approach to the theory of principal $\txG$-bundles due to Ehresmann \cite{Ehresmann:1950},\ we encounter triples $(\sfP,\txG,{\rm At}(\sfP))$ consisting of a principal $\txG$-bundle $\sfP$,\ its structure \emph{group} $\txG$,\ and its structure \emph{groupoid} ${\rm At}(\sfP)$.\ This is captured by Diag.\,\eqref{diag:W-diagram} for $M_2$ a point and $M_1 = \Si$ the base of $\sfP$.\ In the above,\ instead,\ the right wing of the $W$-diagram is a general groupoid.\ This aligns with an idea of Pradines\footnote{The idea later resurfaced in various incarnations in the study of foliations and generalised (Morita) morphisms of groupoids,\ in particular in the works of Hilsum and Skandalis \cite{Hilsum:1983,Hilsum:1987} (leading to the related notion of Hilsum--Skandalis maps),\ and Haefliger \cite{Haefliger:1984}.} \cite{Pradines:1977} (see also \cite{Pradines:2007}),\ whose goal was to symmetrise the structure as in Def.\,\ref{def:bibndl},\ except that here $M_1$ cannot remain to be $\Si$,\ but it is to be replaced,\ in an essential way,\ by the bundle $\xcF$,\ whose fibre is $M_2=M$.\ Thus,\ in distinction to Pradines,\ there emerges the trident diagram \eqref{diag:trident}.
\erem

\beg
In the setting of Example \ref{eg:duck-for-action},\ with $\xcP=\sfP\x M$,\ the action of $\xcG=\txG\x M$ on $\xcP$ coincides with the diagonal action of $\,\txG$ on $\sfP\x M$,\ and the action of ${\rm At}(\xcP)$ is by automorphisms of the factor $\sfP$ in $\sfP\x M$.\ That the two actions commute can be used to see that the automorphisms of $\sfP$ descend to automorphisms of the associated bundle $\xcF=\sfP\x_\la M\equiv(\sfP\x M)/\txG$.\ This includes the important vertical automorphisms,\ which correspond to gauge transformations. 
\eeg

\subsection{Gauge transformations of connections}

\bedef
Let $\xcP$ be a principaloid bundle with connection 1-form $\Theta\in\Om^1(\xcP,\txV\xcP)$.\ Given an automorphism $\Phi\in{\rm Aut}(\xcP)\cong{\rm Bisec}_\pi({\rm At}(\xcP))$,\ we denote 
\qq\label{eq:gauge-transform}
\Theta^\Phi:=\Phi^{-1\,*}\Theta\equiv\txT\Phi\circ\Theta\circ\txT\Phi^{-1}\in\Om^1(\xcP,\txV\xcP)\,.
\qqq
Whenever $\Phi\in{\rm Gauge}(\xcP)\equiv{\rm Aut}(\xcP)_{\rm vert}$,\ the 1-form $\Theta^\Phi$ shall be called the {\bf gauge transform} of $\Theta$.
\exdef

\berop
Let $\xcP$ be a principaloid bundle with connection 1-form $\Theta$ with local connection data $\{\txA_i\in\Om^1(O_i\x M,\pr_2^*E)\}_{i\in I}$ associated with a trivialising cover $\cO\equiv\{\cO_i\}_{i\in I}$ of the base $\Si$.\ For every automorphism $\Phi\in{\rm Aut}(\xcP)$ covering a diffeomorphism $f\in\Diff(\Si)$ and presented by local data $\g_{(j,i)}\colo O^f_{(j,i)}\to\bB,\ (j,i)\in I^{\x 2}_{f,\cO}$,\ local connection data $\{\txA_i^\Phi\}_{i\in I}$ of $\Theta^\Phi$ obey,\ for any $\si\in O^f_{(j,i)}$ and $m\in M$,\ the following transformation law:
\qq\label{eq:connoid-gaulaw}
\txA_j^\Phi\bigl(f(\si),t_*\bigl(\g_{(j,i)}(\si)\bigr)(m)\bigr)\circ\txT_\si f=\txT_{\Id_m}C_{\g_{(j,i)}(\si)}\circ\txA_i(\si,m)-\bigl(\ev_m\circ \g_{(j,i)}\bigr)^*\theta_{\rm R}(\si)\,,
\qqq
For $\Phi\in{\rm Gauge}(\xcP)$,\ the above formula specialises to
\qq\label{gaugetrafo}
\txA_i^\Phi\bigl(\si,t_*\bigl(\g_i(\si)\bigr)(m)\bigr)=\txT_{\Id_m}C_{\g_i(\si)}\circ\txA_i(\si,m)-\bigl(\ev_m\circ \g_i\bigr)^*\theta_{\rm R}(\si)\,,
\qqq
where $\g_i\equiv\g_{(i,i)}$.
\eerop
\beroof
The statement of the proposition is verified through a calculation fully analogous to the derivation of the gluing law for local connection data,\ carried out in the proof of Thm.\,\ref{thm:loc-data-conn}.\ More specifically,\ upon fixing local coordinates:\ $\{\z^\a\}^{\a\in\ovl{1,\dim\,\xcG}}$ on a neighbourhood $\xcU$ of $g\in\xcG$ and $\{\widehat\z{}^\a\}^{\a\in\ovl{1,\dim\,\xcG}}$ on a neighbourhood $\widehat\xcU$ of $\g_{(j,i)}(\si)\lact g$,\ formula \eqref{eq:gauge-transform} yields an equality
\qq\nn
&&\id_{\txT\xcG}\rstr_{\txT_g\xcG}+\txT_{\Id_{t(g)}}r_g\circ\txA_i\bigl(\si,t(g)\bigr)\cr\cr
&=&\txT_{\check\eta{}_{(j,i)}(\si,g)}L_{\g_{(j,i)}(\si)^{-1}}\circ\bigl(\check\eta{}_{(j,i)}^*\txd\widehat\z{}^\a(\si,g)\ox\widehat\p{}_\a\rstr_{\check\eta{}_{(j,i)}(\si,g)}+\txT_{\Id_{t(\check\eta{}_{(j,i)}(\si,g))}}r_{\check\eta{}_{(j,i)}(\si,g)}\circ\txA_j^\Phi\bigl(f(\si),t\bigl(\check\eta{}_{(j,i)}(\si,g)\bigr)\bigr)\circ\txT_\si f\bigr)\,,
\qqq
expressed in terms of the smooth map
\qq\nn
\check\eta{}_{(j,i)}\colo O^f_{(j,i)}\x\xcG\too\xcG,\ (\si,g)\longmapsto \g_{(j,i)}(\si)\lact g\,.
\qqq
Computing the right-hand side of the above condition as in the proof quoted,\ we obtain 
\qq\nn
&&\id_{\txT\xcG}\rstr_{\txT_g\xcG}+\txT_{\Id_{t(g)}}r_g\circ\txA_i\bigl(\si,t(g)\bigr)\cr\cr
&=&\id_{\txT\xcG}\rstr_{\txT_g\xcG}+\txT_{\check\eta{}_{(j,i)}(\si,g)}L_{\g_{(j,i)}(\si)^{-1}}\circ\bigl(\txT_\si\check\eta{}_{(j,i)}(\cdot,g)+\txT_{\Id_{t(\check\eta{}_{(j,i)}(\si,g))}}r_{\check\eta{}_{ij}(\si,g)}\circ\txA_j^\Phi\bigl(f(\si),t\bigl(\check\eta{}_{(j,i)}(\si,g)\bigr)\bigr)\circ\txT_\si f\bigr)\,.
\qqq
The anticipated transformation law now follows readily,\ again,\ as in the said proof.
\eroof

\berop\label{prop:covariance-of-shadow}
For every \emph{vertical} automorphism $\Phi\in{\rm Gauge}(\xcP)$ of a principaloid bundle $\xcP$ with connection 1-form $\Theta$,\ the shadow connection 1-form $\Theta_\xcF\in\Om^1(\xcF,\txV\xcF)$ satisfies the identity
\qq\nn
\xcF_*(\Phi)^*\:\Theta_\xcF^\Phi=\Theta_\xcF\,,
\qqq
in which $\Theta_\xcF^\Phi$ is the shadow connection induced by the gauge transform $\Theta^\Phi$,\ and in which the pullback along $\xcF_*(\Phi)$ is explicitly given by
\qq\nn
\xcF_*(\Phi)^*\:\Theta_\xcF^\Phi\equiv\txT\bigl(\xcF_*(\Phi)\bigr)^{-1}\circ\Theta_\xcF^\Phi\circ\txT\bigl(\xcF_*(\Phi)\bigr)\,.
\qqq
\eerop
\beroof
In the light of Prop.\,\ref{prop:shadconn-idef},\ the gauge-transformed connection 1-form $\Theta_\xcF^\Phi$ satisfies 
\qq\nn
\Theta_\xcF^\Phi\circ\txT\xcD=\txT\xcD\circ\Theta^\Phi\,.
\qqq
Upon sandwiching it between $\txT(\xcF_*(\Phi))^{-1}$ and $\txT\Phi$ and using the tangent variant of the intertwiner property \eqref{eq:duck-as-intertw} of $\xcD$:
\qq\nn
\txT\bigl(\xcF_*(\Phi)\bigr)\circ\txT\xcD=\txT\xcD\circ\txT\Phi\,,
\qqq
alongside definition \eqref{eq:gauge-transform},\ we obtain
\qq\nn
\xcF_*(\Phi)^*\:\Theta_\xcF^\Phi\circ\txT\xcD&=&
\txT\bigl(\xcF_*(\Phi)\bigr)^{-1}\circ\Theta_\xcF^\Phi\circ\txT\xcD\circ\txT\Phi=\txT\bigl(\xcF_*(\Phi)\bigr)^{-1}\circ\txT\xcD\circ\Theta^\Phi\circ\txT\Phi=\txT\xcD\circ\txT\Phi^{-1}\circ\Theta^\Phi\circ\txT\Phi\cr\cr
&=&\txT\xcD\circ\Theta=\Theta_\xcF\circ\txT\xcD\,,
\qqq
where the last equality uses \Reqref{eq:shadconn-idef} again.\ The statement of the proposition now follows by the surjectivity of $\txT\xcD$.
\eroof

\subsection{Covariant derivatives}

\bedef
The {\bf covariant derivative} of a section $\varphi\in\G(\xcF)$ of the shadow bundle $\xcF$ of a principaloid bundle $\xcP$ {\bf relative to} a connection 1-form $\Theta$ on $\xcP$ is the $\bR$-linear mapping
\qq\nn
\nabla^\Theta_\cdot\varphi\colo \G(\txT\Si)\too\G(\txV\xcF),\cV\longmapsto\Theta_\xcF\circ\txT\varphi(\cV)\equiv\nabla^\Theta_\cV\varphi\,.
\qqq
\exdef

\berop
For a section $\varphi\in\G(\xcF)$ of the shadow bundle $\xcF$ of a principaloid bundle $\xcP$,\ and a gauge transformation $\Phi\in{\rm Gauge}(\xcP)$ of $\xcP$,\ let $\varphi^\Phi$ denote the corresponding gauge transform of $\varphi$: 
\qq\nn
\varphi^\Phi\equiv\xcF_*(\Phi)\circ\varphi\,.
\qqq
The covariant derivative of $\varphi$ relative to a connection 1-form $\Theta$ on $\xcP$ transforms \emph{covariantly} as 
\qq\nn
\nabla^{\Theta^\Phi}_\cdot\bigl(\varphi^\Phi\bigr)=\txT\bigl(\xcF_*(\Phi)\bigr)\circ\nabla^\Theta_\cdot\varphi\,.
\qqq
\eerop
\beroof
Upon invoking Prop.\,\ref{prop:covariance-of-shadow},\ we calculate directly:
\qq\nn
&&\nabla^{\Theta^\Phi}_\cdot\bigl(\varphi^\Phi\bigr)\equiv\Theta_\xcF^\Phi\circ\txT(\xcF_*(\Phi)\circ\varphi)=\txT\bigl(\xcF_*(\Phi)\bigr)\circ\bigl(\txT\bigl(\xcF_*(\Phi)\bigr)^{-1}\circ{}^\Phi\hspace{-2pt}\Theta_\xcF\circ\txT\bigl(\xcF_*(\Phi)\bigr)\bigr)\circ\txT\varphi\cr\cr
&=&\txT\bigl(\xcF_*(\Phi)\bigr)\circ\Theta_\xcF\circ\txT\varphi\equiv\txT\bigl(\xcF_*(\Phi)\bigr)\circ\nabla^\Theta_\cdot\varphi\,.
\qqq
\eroof

\appendix

\section{Conventions on maps and differential forms}\label{app:maps}

\becon
Let $M_A$ and $N_A$,\ with $A\in\{1,2\}$,\ be smooth manifolds and $f_A\colo M_A\to N_A$  smooth maps.\ We write
\qq\nn
f_1\x f_2\colo M_1\x M_2\too N_1\x N_2,\ (m_1,m_2)\longmapsto\bigl(f_1(m_1),f_2(m_2)\bigr)\,.
\qqq
\econ

\becon
Let $M$ and $N_A$,\ with $A\in\{1,2\}$,\ be manifolds and $f_A\colo M\to N_A$ smooth maps.\ We write 
\qq\nn
(f_1,f_2)\colo M\too N_1\x N_2,\ m\longmapsto\bigl(f_1(m),f_2(m)\bigr)\,.
\qqq
\econ

\bedef\label{def:bas-fib-pullbck}
Let $\pi_\bV\colo \bV \to B$ be a (real) vector bundle.\ A $\bV$-valued $p$-form on $B$ is an arbitrary global section $\om\in\G(\bigwedge{}^p\,\txT^*B\ox\bV)\equiv\Om^p(B,\bV)$,\ or---equivalently---a vertical vector-bundle morphism
\qq\nn
\alxydim{@C=1.75cm@R=1.5cm}{\bigwedge{}^p\,\txT B \ar[r]^{\quad \om} \ar[d]_{\pi_{\bigwedge{}^p\txT B}} & \bV \ar[d]^{\pi_\bV}\\ B \ar@{=}[r]_{\id_B} & B}\,.
\qqq
\exdef

\bedef\label{def:folia}
Let $M$ be a smooth manifold.\ A ({\bf smooth}) {\bf foliation} on $M$ is an integrable subbundle $\xcF\emb\txT M$,\ {\it i.e.},\ a subbundle with the {\bf Frobenius property}
\qq\nn
[\G(\xcF),\G(\xcF)]\subset\G(\xcF)\,.
\qqq
\exdef

\bedef\label{def:foliaform}
Let $X$ and $M$ be smooth manifolds,\ $\xcF$ a foliation on $X$,\ and $\pi_\bV \colo \bV \to M$ a vector bundle.\ An $\xcF${\bf -foliated differential 1-form on} $X$ {\bf with values in} $\bV$ is a vector-bundle morphism 
\qq\nn
\alxydim{@C=1.5cm@R=1.5cm}{ \xcF \ar[r]^{\eta} \ar[d]_{\pi_{\txT X}\rstr_\xcF} & \bV \ar[d]^{\pi_\bV}\\ X \ar[r]_{h} & M}\,.
\qqq
\exdef

\section{Conventions on Lie groupoids}\label{app:Liegrpd}

\bedef\label{def:grpd}
A \textbf{groupoid} is a small category with all morphisms invertible.\ Thus,\ it is a septuple $\Gr=(\obj\Gr,\morf\Gr, s,  t,\Id,\Inv,\txm\equiv.)$ composed of a pair of sets: 
\bit
\item the \textbf{object set} $\obj\Gr$; 
\item the \textbf{arrow set} $\morf\Gr$,
\eit
and a quintuple of \textbf{structure maps}: 
\bit
\item the \textbf{source map} $ s\colo \morf\Gr\to\obj\Gr$;
\item the \textbf{target map} $  t\colo \morf\Gr\to\obj\Gr$; 
\item the \textbf{unit map} $\Id\colo \obj\Gr\to\morf\Gr,\ m\mapsto\Id_m$;
\item the \textbf{inverse map} $\Inv\colo \morf\Gr\to\morf\Gr,\  g\mapsto \Inv(
 g)\equiv g^{-1}$;
\item the \textbf{multiplication map} $\txm\colo \morf\Gr\,{}_ s\hspace{-3pt}\x_t\hspace{-1pt}\morf\Gr\to\morf\Gr,\ ( g, h)\mapsto \txm( g, h)\equiv g.h$,
\eit 
defined in terms of the subset $\morf\Gr\,{}_ s\hspace{-3pt}\x_  t\hspace{-1pt}\morf\Gr$ of composable morphisms,
\qq\nn
\morf\Gr\,{}_s\hspace{-3pt}\x_t\hspace{-1pt}\morf\Gr=\{\
( g,  h)\in\morf\Gr\x\morf\Gr
\quad\vert\quad  s( g)=  t( h) \ \}\equiv\morf\Gr\x_{\obj\Gr}\morf\Gr\,,
\qqq
and subject to the conditions (in force whenever the expressions are well-defined):
\bit
\item[(i)] $ s( g. h)= s(
 h),\   t( g. h)=  t(
 g)$;
\item[(ii)] $( g. h). k= g.( h.k)$;
\item[(iii)] $\Id_{  t( g)}. g=
 g= g.\Id_{ s( g)}$
;
\item[(iv)] $ s( g^{-1})=  t( g),\   t(
 g^{-1})= s( g),\  g
. g^{-1}=\Id_{  t( g)},\
 g^{-1}. g=\Id_{ s(g)}$.
\eit

A \textbf{morphism} between two groupoids $\Gr_A,\ A\in\{1,2\}$ is a functor $\Phi\colo \Gr_1\to\Gr_2$.

A \textbf{Lie groupoid} is a groupoid whose object and arrow sets are smooth manifolds,\ whose structure maps are smooth,\ and whose source and target maps are surjective submersions.\ A morphism between two Lie groupoids is a functor between them with smooth object and morphism components. 

We shall represent a Lie groupoid $\Gr$ by a (non-commutative) diagram 
\qq\label{diag:Grdiag}\qquad\qquad
\alxydim{@C=1.5cm@R=1.5cm}{ \morf\Gr\x_{\obj\Gr}\morf\Gr \ar[r]^{\qquad\quad\txm} & \morf\Gr \ar[r]^{\Inv} & \morf\Gr \ar@<.25ex>[r]^{ s} \ar@<-.25ex>[r]_{  t} & \obj\Gr \ar@/_1.75pc/[l]_{{\rm Id}} }\,.
\qqq
\exdef

\brem
From now onwards,\ we shall mostly use the shorthand notation $\obj\Gr\equiv M$ and $\morf\Gr\equiv\xcG$,\ and,\ furthermore,\ write $\xcG\sim\Gr$.
\erem

\beg\label{eg:grpasgrpd}
A Lie group $\txG$ can be viewed as a Lie groupoid $(\{\bullet\},\txG,\bullet,\bullet,\bullet\mapsto e,\Inv_\txG,\cdot)$ with the singleton $\{\bullet\}$ as the object manifold.
\eeg

\beg\label{eg:pairgrpd}
An important example of a Lie groupoid is provided by the {\bf pair groupoid} ${\rm Pair}(M)$ of a smooth manifold $\obj\,({\rm Pair}(M))=M$,\ with the object manifold $M$ and the arrow manifold $\morf\,({\rm Pair}(M))=M\x M$,\ the source map $ s=\pr_2$ and the target map $  t=\pr_1$ given by the canonical projections,\ the composition of morphisms
\qq\nn
\txm\colo (M\x M)\,{}_{\pr_2}\hspace{-3pt}\x_{\pr_1}\hspace{-1pt}(M\x M)\too M\x M,\ \bigl((m_3,m_2),(m_2,m_1)\bigr)\longmapsto(m_2,m_1)\,,
\qqq
the unit map
\qq\nn
\Id\colo M\too M\x M,\ m\longmapsto(m,m)\,,
\qqq
and the inversion map
\qq\nn
\Inv\colo M\x M\too M\x M,\ (m_2,m_1)\longmapsto(m_1,m_2)\,.
\qqq
The pair groupoid contains,\ as a proper Lie subgroupoid,\ the {\bf Lie groupoid of} $M$ obtained through restriction of the arrow manifold to the diagonal $M\x_M M\equiv M$.

Whenever $M$ is the total space of a surjective submersion $\pi_M\colo M\to \Si$,\ the pair groupoid ${\rm Pair}(M)$ contains,\ as a proper Lie subgroupoid,\ the {\bf $\Si$-fibred pair groupoid} ${\rm Pair}_\Si(M)$ of $M$,\ with the arrow manifold
\qq\nn
\morf\,({\rm Pair}_\Si(M))=M\x_\Si M\equiv\{\ (m_1,m_2)\in M\x M\quad\vert\quad \pi_M(m_1)=\pi_M(m_2)\ \}\,.
\qqq
\eeg

\beg\label{eg:actgrpd}
Another Lie groupoid of relevance is the {\bf action groupoid} $\txG\,\lx_\la M$ associated with a smooth action 
\qq\nn
\la\colo \txG\x M\too M,\ m\longmapsto\la_g(m)\equiv g\lact m\equiv g.m
\qqq 
of a Lie group $\txG$ on a smooth manifold $M$,\ with the object manifold $\obj\,(\txG\,\lx_\la M)=M$ and the arrow manifold $\morf\,(\txG\,\lx_\la M)=\txG\x M$,\ the source map $ s=\pr_2$ given by the canonical projection,\ the target map $  t=\la$ given by the smooth (left) action $\la$,\ the composition of morphisms
\qq\nn
\txm\colo (\txG\x M)\,{}_{\pr_2}\hspace{-3pt}\x_\la\hspace{-1pt}(\txG\x
M)\too\txG\x M,\ \bigl((h,g.m),(g,m)\bigr)\longmapsto(h\cdot
g,m)=:(h,g.m).(g,m)\,,
\qqq
the unit map ($e\in\txG$ is the group unit)
\qq\nn
\Id\colo M\too\txG\x M,\ m\longmapsto\Id_m=(e,m)\,,
\qqq
and -- finally -- the inversion map
\qq\nn
\Inv\colo \txG\x M\too\txG\x M,\ (g,m)\longmapsto\bigl(g^{-1},g.m\bigr)
=:\Inv(g,m)\equiv(g,m)^{-1}\,.
\qqq

The action groupoid can be endowed with the structure of a $\txG\x\txG$-space in several ways,\ amidst which there is a distinguished one:\ The (left) action of $\txG\x\txG$ intertwined with $\la\circ(\pr_2\x\id_M)$ by $s$ and with $\la\circ(\pr_1\x\id_M)$ by $  t$ is given by the formula
\qq\label{eq:ell1}\qquad\qquad
\xcL\colo (\txG\x\txG)\x(\txG\x M)\too\txG\x M,\ \bigl((h_1,h_2),(g,m)\bigr)
\longmapsto\bigl(h_1\cdot g\cdot h_2^{-1},h_2.m\bigr)\,.
\qqq
This can equivalently be seen as an `adjoint' action of $\txG\,\lx_\la M$ on itself.
\eeg

\beg\label{eg:symplgrpd}
Let $(M,\xcG, s,  t,\Id, \Inv,\txm)$ be a Lie groupoid.\ We call it a {\bf symplectic groupoid},\ after Refs.\,\cite{Weinstein:1987,Coste:1987},\ whenever there exists a closed non-degenerate 2-form $\om\in\Om^2(\xcG)$ which makes $(\xcG,\om)$ a symplectic manifold,\ and such that the graph ${\rm graph}(\txm)$ of $\txm$ is a lagrangean submanifold in $(\xcG^{\x 3},(\pr_1^*+\pr_2^*-\pr_3^*)\om)$.\ As shown in \Rxcite{Thm.\,1.1,\ Chap.\,II}{Coste:1987},\ the object manifold $M$ of a symplectic groupoid carries a canonical Poisson bivector $\Pi\in\G(\bigwedge^2\txT^*M)$ relative to which the source map is Poisson and the target map is anti-Poisson.
\eeg

\beg\label{eg:cotgrp}
The geometry of the cotangent bundle $\txT^*\txG$ of a Lie group $\txG$ combines structures from Examples \ref{eg:actgrpd} and \ref{eg:symplgrpd}.\ Indeed,\ its global trivialisation,\ $  \txT^*\txG\cong\txG\x\ggt^*$,\ turns it into the arrow manifold of the action groupoid $\txG\lx_{\txT_e\Ad^*}\ggt^*$ (relative to the coadjoint action),\ and the choice of the Poincar\'e 2-form on it ({\it i.e.},\ the exterior derivative of the tautological 1-form) as the symplectic form gives rise to a symplectic groupoid,\ see \cite{Coste:1987}.
\eeg

\beg\label{eg:GLV}
Let $(\bV,B,\bR^N,\pi_\bV)$ be a (real) vector bundle of rank $N\in\bN^\x$ with base $B$.\ The {\bf general linear groupoid of} $\bV$ is the Lie groupoid ${\rm GL}(\bV)$ with the object manifold $\obj\,{\rm GL}(\bV)=B$,\ the arrow manifold
\qq\nn
\morf\,{\rm GL}(\bV)=\bigsqcup_{x,y\in B}\,{\rm Iso}_\bR(\bV_x,\bV_y)\,,
\qqq
the source and arrow maps 
\qq\nn
 s\ &:&\ \bigsqcup_{x,y\in B}\,{\rm Iso}_\bR(\bV_x,\bV_y)\too B,\ (\chi_{u,v},u,v)\longmapsto u\,,\cr\cr
  t\ &:&\ \bigsqcup_{x,y\in B}\,{\rm Iso}_\bR(\bV_x,\bV_y)\too B,\ (\chi_{u,v},u,v)\longmapsto v\,,
\qqq
the composition of morphisms
\qq\nn
\txm\colo \morf\,{\rm GL}(\bV)\,{}_s\hspace{-3pt}\x_t\hspace{-1pt}\morf\,{\rm GL}(\bV)\too\morf\,{\rm GL}(\bV),\ \bigl((\chi_{y,z},y,z),(\chi_{x,y},x,y)\bigr)\longmapsto(\chi_{y,z}\circ\chi_{x,y},x,z)\,,
\qqq
the unit map
\qq\nn
\Id\colo \obj\,{\rm GL}(\bV)\too\morf\,{\rm GL}(\bV),\ x\longmapsto(\id_{\bV_x},x,x)\,,
\qqq
and the inversion map
\qq\nn
\Inv\colo \bigsqcup_{x,y\in B}\,{\rm Iso}_\bR(\bV_x,\bV_y)\circlearrowleft,\ (\chi_{u,v},u,v)\longmapsto(\chi_{u,v}^{-1},v,u)\,.
\qqq
\eeg

\bedef\label{def:SES-Lie}
Let $\Gr_a,\ a\in\{1,2,3\}$ be Lie groupoids and let $\jmath\colo \Gr_1\to\Gr_2$ and $\pi\colo \Gr_2\to\Gr_3$ be Lie-groupoid morphisms.\ We say that the quintuple $(\Gr_1,\Gr_2,\Gr_3,\jmath,\pi)$ composes a {\bf short exact sequence of Lie groupoids}
\qq\nn
\alxydim{@C=1cm@R=1.cm}{ \bd1 \ar@{-->}[r] & \Gr_1 \ar[r]^{\jmath} & \Gr_2 \ar[r]^{\pi} & \Gr_3 \ar@{-->}[r] & \bd1 }
\qqq
if $\jmath$ is a monomorphism,\ $\pi$ is an epimorphism,\ and the $\pi$-preimage of the identity bisection $\Id(\obj\Gr_3)\subset\morf\Gr_3$ is canonically isomorphic to the (faithful) $\jmath$-image of $\Gr_1$ in $\Gr_2$.\ If,\ moreover,\ there exists a Lie-groupoid morphism $\si\colo \Gr_3\to\Gr_2$ such that the identity $\pi\circ\si=\id_{\Gr_3}$ obtains,\ then we say that the short exact sequence {\bf splits},\ and call it a {\bf split short exact sequence} ({\bf of Lie groupoids}).
\exdef

\bedef[\cite{Moerdijk:2003mm}]\label{def:bisec}
Let $\Gr=(M,\xcG, s,  t,\Id,\Inv,.)$ be a Lie groupoid.\ A ({\bf global}) {\bf bisection of} $\Gr$ is a section $\b\colo M\to\xcG$ of the surjective submersion $ s\colo \xcG\to M$ such that the induced map 
\qq\nn
  t_*\b\equiv  t\circ\b\colo M\too M
\qqq
is a diffeomorphism.\ Equivalently,\ it is a submanifold $S\subset\xcG$ with the property that both restrictions:\ $ s\rstr_S$ and $t\rstr_S$ are diffeomorphisms.\ We shall denote the set of bisections as 
\qq\nn
\BisGr\,.
\qqq

A {\bf local bisection of} $\Gr$ is a local section $\b\colo O\to\xcG$ of $ s$ over an open subset $O\subset M$,\ such that the induced map 
\qq\nn
  t_*\b\equiv  t\circ\b\colo O\too  t\circ\b(O)
\qqq
is a diffeomorphism.\ We shall denote the set of local bisections as ${\rm Bisec}_{\rm loc}(\Gr)$.
\exdef

\bedef
The {\bf group of bisections of} $\Gr$ is the canonical structure of a group on $\BisGr$.\ Its binary operation is defined as
\qq\nn
\cdot\colo \BisGr \x\BisGr \too\BisGr,\ (\b_2,\b_1)\longmapsto\b_2\bigl(t\circ\b_1(\cdot)\bigr).\b_1(\cdot)\equiv\b_2\cdot\b_1\,.
\qqq
The neutral element is $\Id$,\ termed the {\bf unit bisection} in the present context,\ and the corresponding inverse is 
\qq\nn
\Inv\colo \BisGr\too\BisGr,\ \b\longmapsto\Inv\circ\b\circ\bigl(  t_*\b\bigr)^{-1}\equiv\b^{-1}\,,
\qqq
where the inverse on the right-hand side is the one in $\Diff(M)$.
\exdef

\begin{figure}[hbt!]

$$
 \raisebox{-50pt}{\begin{picture}(0,195)
  \put(-160,-5){\scalebox{0.5}{\includegraphics{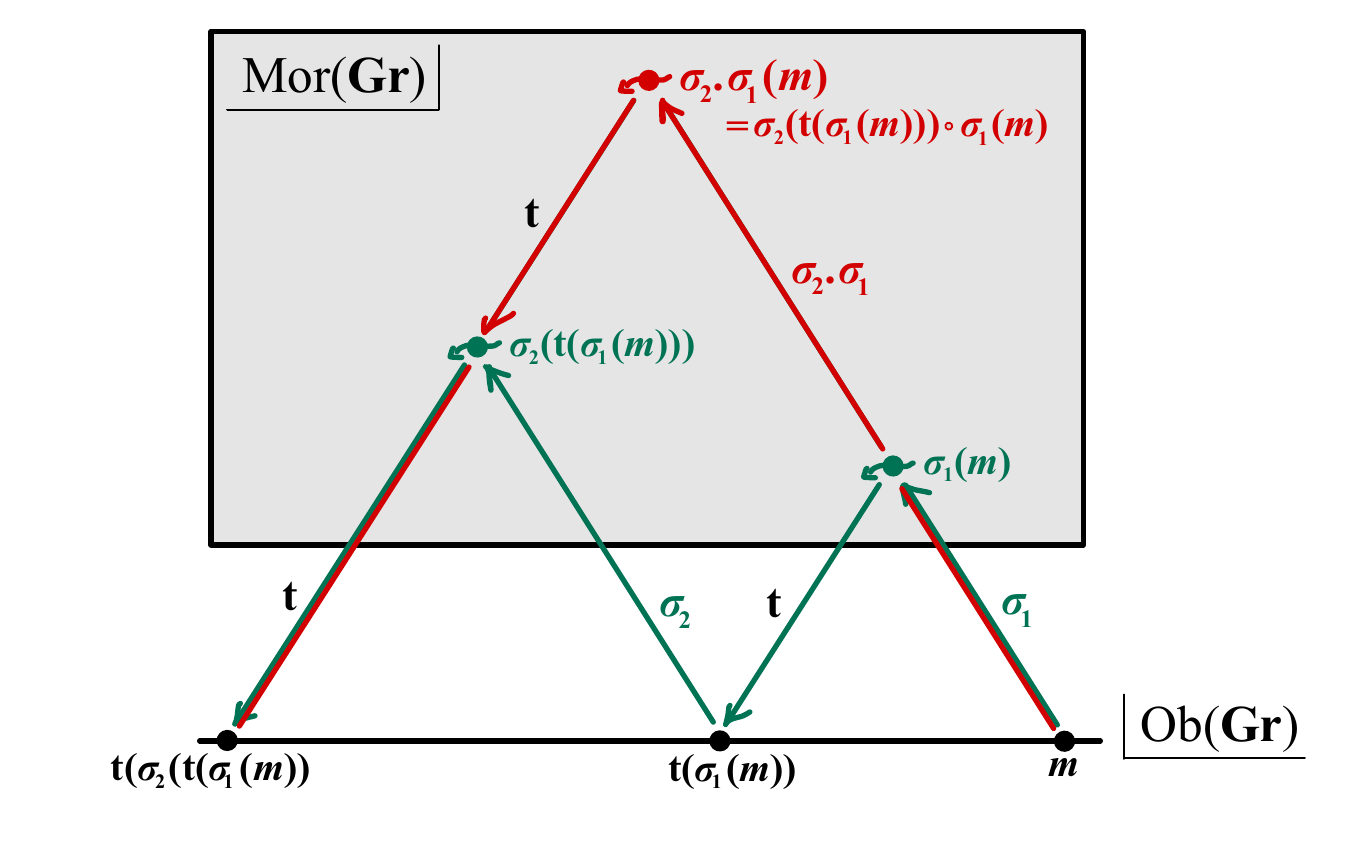}}}
  \end{picture}}
$$

\caption{The binary operation on $\BisGr$.} \label{fig:bisec_bin}
\end{figure}

\brem
The binary operation of $\BisGr$ admits a simple pictorial representation in Fig.\,\ref{fig:bisec_bin}
\erem

\brem
For a fixed Lie groupoid,\ we shall often use the shorthand notation $\BisGr\equiv\bB$.
\erem

\bedef
The {\bf left regular action of} $\bB$ ({\bf on itself}) is
\qq\nn
\ell\colo \bB\x\bB\too\bB,\ (\g,\b)\longmapsto\g\cdot\b\equiv\ell_\g(\b)\,,
\qqq
and the {\bf right regular action of} $\bB$ ({\bf on itself}) is
\qq\nn
\wp\colo \bB\x\bB\too\bB,\ (\b,\g)\longmapsto\b\cdot\g\equiv\wp_\g(\b)\,.
\qqq
The {\bf adjoint action of} $\bB$ {\bf on itself} is
\qq\nn
c\colo \bB\x\bB\too\bB,\ (\g,\b)\longmapsto\g\cdot\b\cdot\g^{-1}\equiv c_\g(\b)\,.
\qqq
\exdef

\bedef
The {\bf shadow action of} $\bB$ {\bf on} $M$ is 
\qq\nn
t_*\colo \bB\x M\too M,\ (\b,m)\longmapsto t\bigl(\b(m)\bigr)\,.
\qqq
By the usual abuse of the notation,\ we shall refer by the same name to and use the same symbol for the group homomorphism 
\qq\label{eq:BisGr-act-obj}
  t_*\colo \bB \too{\rm Diff}(M)\,.
\qqq 
\exdef

\bedef\label{def:bisec-act}
The {\bf left-multiplication of} $\xcG$ {\bf by} $\bB$ is the left action
\qq\nn
 L\colo \bB \x\xcG\too\xcG,\ \bigl(\b, g\bigr)\longmapsto\b\bigl(  t\bigl( g\bigr)\bigr). g\equiv L_\b\bigl( g\bigr)\equiv \b\lact g\,.
\qqq
The {\bf right-multiplication of} $\xcG$ {\bf by} $\bB$ is the right action
\qq\nn
 R\colo  \xcG\x\bB\too\xcG,\ \bigl( g,\b\bigr)\longmapsto g.\bigl(\b^{-1}\bigl(s\bigl( g\bigr)\bigr)\bigr)^{-1}\equiv R{}_\b\bigl( g\bigr)\equiv g\ract\b\,.
\qqq
The {\bf conjugation of} $\xcG$ {\bf by} $\bB$ is the \emph{left} action
\qq\nn
 C\colo \bB \x\xcG\too\xcG,\ (\b, g)\longmapsto\b\bigl(  t( g)\bigr). g.\b\bigl( s( g)\bigr)^{-1}\equiv C_\b( g)\equiv \b\lact g\ract\b^{-1}\,.
\qqq
\exdef
\brem
Taking into account the significance of the right-multiplication in the present paper,\ we simplify its definition for the sake of transparency.\ Thus,\ we find
\qq\nn
g\ract\b\equiv g.\b\bigl((t_*\b)^{-1}\bigl(s(g)\bigr)\bigr)\,.
\qqq
\erem

\section{Groupoid modules and principal groupoid bundles}\label{app:princ-grpd-bndl}

\bedef\label{def:gr-mod}
Given a Lie groupoid $\Gr=(M,\xcG,s,t,\Id,\Inv,.)$,\ a \textbf{right-$\xcG$-module space} is a triple $(X,\mu,\varrho)$ composed of 
\bit
\item a smooth manifold $X$;
\item a smooth map $\mu\colo X\to M$,\ called the (\textbf{right}) \textbf{moment map};
\item a smooth map
\qq\nn
\varrho\colo  X{}_{\mu}\hspace{-3pt}\x_t\hspace{-1pt}\xcG\too
 X,\ (x, g)\longmapsto\varrho(x, g)\equiv\varrho_{ g}(x)\equiv x\mact g\,,
\qqq
termed the (\textbf{right}) \textbf{action} ({\bf map}),
\eit 
subject to the relations (in force whenever the expressions are well-defined):
\bit
\item[(GrM1)] $\mu(x\mact g)= s( g)$;
\item[(GrM2)] $x\mact\Id_{\mu(x)}=x$;
\item[(GrM3)] $(x\mact g)\mact h=x\mact( g. h)$.
\eit

Similarly,\ a {\bf left-$\xcG$-module space} is a triple $(X,\mu,\la)$ composed of 
\bit
\item a smooth manifold $X$;
\item a smooth map $\mu\colo X\to M$,\ called the (\textbf{left}) \textbf{moment map};
\item a smooth map
\qq\nn
\la\colo \xcG{}_{s}\hspace{-2pt}\x_\mu\hspace{-1pt}X\too
 X,\ (g,x)\longmapsto\la(g,x)\equiv\la_g(x)\equiv g\mlact x\,, 
\qqq
termed the (\textbf{left}) \textbf{action} ({\bf map}),
\eit 
subject to the relations (in force whenever the expressions are well-defined):
\bit
\item[(GlM1)] $\mu(g\mlact x)= t(g)$;
\item[(GlM2)] $\Id_{\mu(x)}\mlact x=x$;
\item[(GlM3)] $h\mlact(g\mlact x)=(h.g)\mlact x$.
\eit

A right action $\varrho$ is termed \textbf{free} iff the following implication obtains:
\qq\nn
x\mact g=x\qquad\Longrightarrow\qquad g=
\Id_{\mu(x)}\,,
\qqq
so that,\ in particular,\ the \textbf{isotropy group} $\xcG_m= s^{-1}(\{m\})\cap  t^{-1}(\{m\})$ of $m\in M$ acts freely (in the usual sense) on the fibre $\mu^{-1}(\{m\})$.\ A free left action is defined analogously.

A right action $\varrho$ is termed \textbf{transitive} iff for any two points $x,x'\in X$ there exists an arrow $ g\in\xcG$ such that $x'=x\mact g$.\ A transitive left action is defined analogously.

Let $\Gr_A,\ A\in\{1,2\}$ be a pair of Lie groupoids and let $( X_A,\mu_A,\varrho^A)$ be the respective
right-$\xcG_A$-module spaces.\ A \textbf{morphism} between the latter is a pair $(\Theta,\Phi)$ consisting of a smooth manifold map $\Theta\colo  X_1\to X_2$ together with a functor $\Phi\colo \Gr_1\to\Gr_2$ for which the following diagrams commute
\qq
&\alxydim{@C=1.5cm@R=1.25cm}{ X_1 \ar[r]^{\Theta} \ar[d]_{\mu_1} &
 X_2 \ar[d]^{\mu_2}\cr M_1 \ar[r]_{\Phi} &
M_2}\,,& \label{diag:mu-Th-mu}\\\cr\cr
&\alxydim{@C=1.75cm@R=1.25cm}{ X_1\fibx{\mu_1}{  t_1}\,\xcG_1
\ar[r]^{\Theta\x\Phi} \ar[d]_{\varrho^1} &  X_2
\fibx{\mu_2}{  t_2}\,\xcG_2 \ar[d]^{\varrho^2}\cr  X_1
\ar[r]_{\Theta} &  X_2}\,.& \label{diag:Th-intertw}
\qqq
\exdef

\noindent There are two canonical left-module structures associated with the groupoid itself,\ which we stumble upon in our considerations.\ The first of them is indicated in  
\beg\label{eq:GrpdonOb}
The object manifold $M$ of a Lie groupoid $\Gr=(M,\xcG, s,  t,\Id,\Inv,.)$ carries a natural structure of a left-$\xcG$-module space given by
\qq\nn
\bigl(M,\id_M,{\rm Aim}^\xcG\bigr)\,,
\qqq
where the action reads
\qq\nn
{\rm Aim}^\xcG:=  t\circ\pr_1\colo \xcG{}_{ s}\hspace{-3pt}\x_{\id_M}\hspace{-1pt}M\too M,\ \bigl( g, s( g)\bigr)\longmapsto  t( g)\equiv{\rm Aim}^\xcG_g\bigl( s(g)\bigr)\,.
\qqq
\eeg
\noindent Another one is described in the following
\beg\label{eg:lr-fibreg-act}
We may regard $\xcG$ as a left-$\xcG^{\x 2}$-module space for the \textbf{product groupoid}
\qq\nn
\alxydim{@C=2.5cm}{\Gr^{\x 2}\qquad : \hspace{-2cm} & {\rm Mor}\bigl(\Gr^{\x 2}\bigr)=\xcG^{\x 2} \ar@<.5ex>[r]^{ s^{\x 2}:= s\x s} \ar@<-.5ex>[r]_{  t^{\x 2}:=  t\x  t} & M^{\x 2}={\rm Ob}\bigl(\Gr^{\x 2}\bigr) }\,,
\qqq
with the composition map defined as
\qq\nn
\bigl( g{}_1, g{}_2\bigr)\bullet\bigl( h{}_1, h{}_2\bigr):=\bigl( g{}_1. h{}_1, g{}_2. h{}_2\bigr)\,,
\qqq
the identity map
\qq\nn
\Id_{(m_1,m_2)}=(\Id_{m_1},\Id_{m_2})\,,
\qqq
and the inversion map
\qq\nn
\Inv(g_1,g_2)=\bigl(\Inv(g_1),\Inv(g_2)\bigr)\,.
\qqq
The left-$\xcG^{\x 2}$-module space 
\qq\nn
(\xcG,\mu_{\rm c},{\rm c})
\qqq 
is the triple composed of the smooth manifold $\xcG$,\ the moment map
\qq\nn
\mu_{\rm c}\equiv( t,s)\colo \xcG\too M^{\x 2}\,,
\qqq
defining the set
\qq\nn
\xcG^{\x 2}{}_{s^{\x 2}}\hspace{-3pt}\x_{\mu_{\rm c}}\hspace{-1pt}\xcG:=\bigl\{\ \bigl( ( h_1, h_2), g\bigr)\in\xcG^{\x 3} \quad\vert\quad  s^{\x 2}( h_1,  h_2)=\mu_{\rm c}( g) \ \bigr\}\,,
\qqq
and the action
\qq\nn
{\rm c}\colo \xcG^{\x 2}{}_{ s^{\x 2}}\hspace{-3pt}\x_{\mu_{\rm c}}\hspace{-1pt}\xcG\too\xcG,\ \bigl(( h_1, h_2),  g\bigr)\longmapsto h_1. g. h_2^{-1}\,.
\qqq
Given these data for $\Gr=\txG\,\lx_\la M$,\ we may rewrite \Reqref{eq:ell1} as 
\qq
\xcL\bigl((h_1,h_2),(g,m)\bigr)=(h_1,g.m).(g,m).\bigl(h_2^{-1},h_2.m\bigr)={\rm c}\bigl(\bigl(\bigl(h_1,  t(g,m)\bigr),\bigl(h_2, s(g,m)\bigr)\bigr),(g,m)\bigr)\,.\nn \\ \label{eq:Ad-symm-grpd}
\qqq

Note also that the above left-$\xcG^{\x 2}$-module structure on $\xcG$ is,\ in fact,\ a combination of two `chiral' $\xcG$-module structures on the same manifold:\ the left-$\xcG$-module structure 
\qq\nn
(\xcG,  t,l\equiv.)\,,
\qqq
with the {\bf left-fibred action}
\qq\label{eq:Lgrpdact}
l\bigl( h, g\bigr)\equiv l_{ h}\bigl( g\bigr)= h. g\,,
\qqq
and the right-$\xcG$-module structure 
\qq\label{eq:can-R-Gr-mod}
(\xcG, s,r\equiv.)\,,
\qqq
with the {\bf right-fibred action}
\qq\label{eq:Rgrpdact}
r\bigl( g, h\bigr)\equiv r_{ h}\bigl( g\bigr)= g. h\,.
\qqq
\eeg

\beg\label{eg:GLalgbd} 
Let $\bE=(\cE,M,\bR^N,\pi_\cE,\a_\cE,[\cdot,\cdot]_\cE)$ be a Lie algebroid over the object manifold of a Lie groupoid $\Gr=(M,\xcG, s,  t,\Id,\Inv,.)$.\ A {\bf linear representation of} $\Gr$ in $\bE$ is a Lie-groupoid morphism $\varrho\colo \Gr\to{\rm GL}(\bV)$,\ see Example \ref{eg:GLV}.\ It endows the vector bundle $\bV$ with the structure of a left-$\xcG$-module space,\ with moment map $\mu_\bV=\pi_\bV$ and action
\qq\nn
\varrho^\bV\equiv\ev\colo \bigsqcup_{x,y\in B}\,{\rm Iso}_\bR(\bV_x,\bV_y){}_{ s}\hspace{-1pt}\x_{\mu_\bV}\hspace{-3pt}\bV\too\bV,\ \bigl((\chi_{x,y},x,y),v\bigr)\longmapsto\chi_{x,y}(v)\,.
\qqq
\eeg

\bedef\cite[Sec.\,1.2]{Moerdijk:1991}\label{def:princ-gr-bun}
A \textbf{right principal $\xcG$-bundle} is a quintuple $(\breP,\Si,\pi_{\breP},\mu,\varrho)$ composed of a pair of smooth manifolds: 
\bit
\item the \textbf{total space} $\breP$;
\item the \textbf{base} $\Si$,
\eit 
and a triple of smooth maps: 
\bit
\item a surjective submersion $\pi_{\breP}\colo \breP\to\Si$,\
termed the \textbf{base projection}; 
\item the \textbf{moment map} $\mu\colo \breP\to M$;
\item the \textbf{action} ({\bf map}) $\varrho\colo \breP{}_{\mu}\hspace{-3pt}\x_t\hspace{-1pt}\xcG\to\breP$,
\eit
with the following properties:
\bit
\item[(PGr1)] $(\breP,\mu,\varrho)$ is a right $\xcG$-module space;
\item[(PGr2)] $\pi_{\breP}$ is $\xcG$-invariant in the sense made precise by the following commutative diagram (in which $\pr_1$ is the canonical projection)
\qq\nn
\alxydim{@C=1.5cm@R=1.cm}{\breP{}_{\mu}\hspace{-3pt}\x_t\hspace{-1pt}
\xcG \ar[r]^{\quad\varrho} \ar[d]_{\pr_1} & \breP
\ar[d]^{\pi_{\breP}}\cr \breP \ar[r]_{\pi_{\breP}} & \Si}\,;
\qqq
\item[(PGr3)] the map
\qq\nn
(\pr_1,\varrho)\colo \breP{}_{\mu}\hspace{-3pt}\x_t\hspace{-1pt}\xcG\too\breP\fibx{\pi_{\breP}}{\pi_{\breP}}\,\,\breP,\ (p,g)\longmapsto(p,p\mact g)
\qqq
is a diffeomorphism,\ so that $\xcG$ acts freely and transitively on $\pi_{\breP}$-fibres.\ The smooth inverse of $(\pr_1,\varrho)$ takes the form
\qq\nn
(\pr_1,\varrho)^{-1}=:(\pr_1,\phi_{\breP})\,,\qquad\qquad\phi_{\breP}\colo \breP\fibx{\pi_{\breP}}{\pi_{\breP}}\,\,\breP\too\xcG
\qqq
and $\phi_{\breP}$ is called the \textbf{division map}.
\eit
\noindent The Lie groupoid $\xcG$ is termed the {\bf structure groupoid} of $\breP$.

We shall represent a right principal $\xcG$-bundle by the simplified diagram 
\qq\label{diag:princGrdiag}\qquad\qquad
\alxydim{@C=.75cm@R=1cm}{ \breP \ar@{->>}[d]_{\pi_\breP} \ar[rd]^{\mu} & & \morf\Gr \ar@{=>}[ld] \\ \Si & \obj\Gr & }\,,
\qqq
in which the remaining structure is implicit.

Let $(\breP{}_A,\Si,\pi_{\breP{}_A},\mu_A,\varrho^A),\ A\in\{1,2\}$ be a pair of right principal $\xcG$-bundles over a common base $\Si$.\ A \textbf{morphism}\footnote{In \Rxcite{Sec.\,5.7}{Moerdijk:2003mm},\ these morphisms were termed ``equivariant maps''.} between the two bundles is a morphism $(\Theta,\Id_\Gr)$ between the corresponding right $\xcG$-modules $(\breP{}_A,\mu_A,\varrho^A)$ which maps $\pi_{\breP{}_1}$-fibres to $\pi_{\breP{}_2}$-fibres. 

{\bf Left principal $\xcG$-bundles} $(\breP,\Si,\pi_\breP,\mu,\la)$ (and morphisms between them) are defined analogously.\ The corresponding diagrams take the self-explanatory form
\qq\nn
\alxydim{@C=.75cm@R=1cm}{ \morf\Gr \ar@{=>}[rd] & & \breP \ar@{->>}[d]^{\pi_\breP} \ar[ld]_{\mu} \\ & \obj\Gr & \Si}\,,
\qqq

Customarily,\ principal $\xcG$-bundles are taken to be right $\xcG$-modules,\ and so whenever the term is used without a qualifier,\ it is to be understood that we are dealing with a right principal $\xcG$-bundle.
\exdef

\beg\cite[Remark 5.34(1)]{Moerdijk:2003mm}\label{def:triv-grpd-bndle}
There is a canonical structure of a principal $\xcG$-bundle on $\Gr$,\ given by 
\qq\nn
U_\Gr:=(\xcG,M,t,s,r)\,,
\qqq 
where $r$ is the right-fibred action \eqref{eq:Rgrpdact}.\ We call this structure the \textbf{unit bundle of} $\Gr$.

Axioms (PGr1) and (PGr2) from Def\,\ref{def:princ-gr-bun} are self-evident,\ and so it remains to verify property (PGr3):\ The map
\qq\nn
(\pr_1,r)\colo \xcG{}_s\hspace{-3pt}\x_t\hspace{-1pt}\xcG\too\xcG{}_t\hspace{-3pt}\x_t\hspace{-1pt}\xcG,\ (g,h)\longmapsto(g,g.h)
\qqq
admits the smooth inverse
\qq\nn
(\pr_1,\phi_{U_\Gr})\colo \xcG{}_t\hspace{-3pt}\x_t\hspace{-1pt}\xcG\too\xcG{}_s\hspace{-3pt}\x_t
\hspace{-1pt}\xcG,\ (g,h)\longmapsto\bigl(g,g^{-1}.h\bigr)\,.
\qqq
\eeg

\beg\cite[Remark 5.34(2),(3)]{Moerdijk:2003mm}
A \textbf{trivial principal $\xcG$-bundle} over $\Si$ is the pullback $f^*U_\Gr$ of the unit bundle $U_\Gr$ along an arbitrary smooth map $f\colo \Si\to M$,\ that is the principal $\xcG$-bundle 
\qq\nn
f^*U_\Gr:=\bigl(f^*\xcG,\Si,\pr_1,s\circ\pr_2,f^*\varrho\bigr)
\qqq
with $f^*\xcG=\Si\,{}_f\hspace{-3pt}\x_t\hspace{-1pt}\xcG$ and
\qq\nn
f^*\varrho\colo f^*\xcG\fibx{s\circ\pr_2}{t}\hspace{2pt}\xcG\too f^*\xcG,\ \bigl((\si,g),h\bigr)\longmapsto(\si,g.h)\,.
\qqq
\eeg

\begin{propanition}\cite[Remark 5.34(4)]{Moerdijk:2003mm}\label{eq:3M-loc-triv}
Let $(\breP,\Si,\pi_{\breP},\mu,\varrho)$ be a principal $\xcG$-bundle.\ Given any open subset $O\subset\Si$ endowed with a local section $s\colo O\to\breP$,\ define the corresponding {\bf local moment map} as $\mu_O:=\mu\circ s$.\ The map $\t_O\colo\pi_{\breP}^{-1}(O)\to\mu_O^*U_\Gr,\ p\mapsto\phi_{\breP}(s(\pi_{\breP}(p)),p)$ is a diffeomorphism with inverse $\t_O^{-1}\colo \mu_O^*U_\Gr\to\pi_{\breP}^{-1}(O)),\ (\si,g)\mapsto s(\si)\mact g$.\ We call $\t_O$ a {\bf local trivialisation} of $\breP$ over $O$.\ Thus,\ every principal $\xcG$-bundle is locally trivial.
\end{propanition}

\begin{coronition}\label{cor:3M-model}
Let $(\breP,\Si,\pi_{\breP},\mu,\varrho)$ be a principal $\xcG$-bundle,\ and let $\cO\equiv\{O_i\}_{i\in I}$ be an open cover of its base $\Si$ endowed with local sections $s_i\colo O_i\to\breP$.\ The bundle admits a model 
\qq\label{eq:3M-model}
\breP\cong\bigsqcup_{i\in I}\,\mu_i^*U_\Gr/\sim_{l_{\g_{\cdot\cdot}}}\,,
\qqq
in which the $\mu_i\equiv\mu\circ s_i\colo O_i\to\breP$ are the local moment maps induced by the $s_i$,\ and the quotient is taken with respect to the equivalence relation 
\qq\nn
(\si,g,i)\sim(\si',g',j)\qquad\Longleftrightarrow\qquad\bigl(\ \si'=\si\in O_{ij}\quad\land\quad g'=\g_{ji}(\si).g\ \bigr)\,,
\qqq
written in terms of the {\bf transition 1-cocycle}
\qq\nn
\g_{ij}:=\phi_\breP\circ(s_i,s_j)\colo O_{ij}\too\xcG\,.
\qqq
The collection $(\{\g_{ij}\}_{(i,j)\in I^{\x 2}_\cO},\{\mu_i\}_{i\in I})$ is termed the {\bf trivialising data} of $\breP$.
\end{coronition}
\beroof
A constructive proof relies on the identities 
\qq\label{eq:def-props-trans-arrows}
s\circ\g_{ij}=\mu_j\,,\qquad t\circ\g_{ij}=\mu_i\,,\qquad\g_{ij}.\g_{jk}\rstr_{O_{ijk}}=\g_{ik}\rstr_{O_{ijk}}
\qqq
(following directly from Def.\,\ref{def:princ-gr-bun}),\ and is conceptually the same as that of the classic Clutching Theorem for fibre bundles.\ Its technical details can be found,\ {\it e.g.},\ in \Rxcite{Sec.\,3.2}{Rossi:2004b}.
\eroof

\bethe[The Godement Criterion]\cite{Serre:1964,Fernandes:2024}\label{thm:Godement}
Let $M$ be a smooth manifold and let $\sim$ be an equivalence relation on M,\ with graph
\qq\nn
\alxydim{@C=1.cm@R=1.5cm}{ \cR_\sim \ar@{^{(}->}[r] & M\x M \ar[dl]_{\pr_1} \ar[dr]^{\pr_2} \\ M & & M}\,.
\qqq
There exists a smooth structure on the quotient 
\qq\nn
M//\sim=\{\ [m]_\sim \quad\vert\quad m\in M \ \}\,, 
\qqq
compatible with the quotient topology,\ and such that $\pi\colo M\too M//\sim$ is a submersion,\ iff the graph $\cR_\sim$ is a proper submanifold of $M\x M$ and the restriction of the projection $\pr_1\colo M\x M\too M$ to $\cR_\sim$ is a submersion.
\ethe

\brem\label{rem:Godement}
In virtue of the above Theorem,\ there exists a canonical diffeomorphism
\qq\nn
\breP/\xcG\cong\Si
\qqq
for every principal $\xcG$-bundle $\,\breP$,\ with the left-hand side given by the smooth quotient of the total space $\,\breP\ni p_1,p_2\,$ by the equivalence relation 
\qq\nn
p_1\sim p_2\qquad\Longleftrightarrow\qquad\exists\ g\in\xcG\colo p_2=p_1\mact g\,.
\qqq
\erem

\brem
The definition of a principal $\xcG$-bundle can be viewed as a structural relation between three Lie groupoids.\ Indeed,\ we may rephrase it as a statement of existence,\ for a given surjective submersion $\pi_\breP\colo \breP\to\Si$ and a (smooth) map $\mu\colo \breP\to M$,\ of 
\bit
\item an action Lie groupoid $\breP{}_\mu\hspace{-3pt}\rx_t\hspace{-1pt}\xcG$ with $\breP$ as the object manifold and $\breP{}_\mu\hspace{-3pt}\x_t\hspace{-1pt}\xcG$ as the arrow manifold,\ with $\pr_1$ as the source map,\ a smooth map $\varrho\colo \breP{}_\mu\hspace{-3pt}\x_t\hspace{-1pt}\xcG\to\breP$ as the target map,\ $(\id_\breP,\Id\circ\mu)$ as the unit map,\ $(\varrho,\Inv\circ\pr_2)$ as the inverse map,\ and $(\varrho(p,g),h).(p,g)=(p,g.h)$ as the multiplication map,
\item a Lie-groupoid morphism $\Phi_\breP\colo{\rm Pair}_\Si(\breP)\to\Gr$ with $\mu$ as the object component and a smooth map $\phi_\breP\colo \breP\x_\Si\breP\to\xcG$ as the morphism component,
\eit
such that the Lie-groupoid morphism $\widetilde\varrho\colo\breP{}_\mu\hspace{-3pt}\rx_t\hspace{-1pt}\xcG\to{\rm Pair}_\Si(\breP)$,\ with the object component $\id_\breP$ and the morphism component $(\pr_1,\varrho)$,\ is invertible,\ with the inverse $\widetilde\phi{}_\breP$ given by ($\id_\breP$ on objects,\ and) $(\pr_1,\phi_\breP)$ on morphisms.\ Thus,\ the entire information on $\breP$ is neatly encoded in the following commutative diagram in the category of Lie groupoids: 
\qq\nn
\alxydim{@C=.75cm@R=1.5cm}{ & \breP{}_\mu\hspace{-3pt}\rx_t\hspace{-1pt}\xcG \ar@<.25ex>[ld]_{\widetilde\varrho} \ar[rd]^{\widehat\mu} & \\ {\rm Pair}_\Si(\breP) \ar[rr]_{\Phi_\breP} \ar@<.25ex>[ru]_{\widetilde\phi{}_\breP\equiv\widetilde\varrho{}^{-1}} & & \Gr }\,,
\qqq
in which the Lie-groupoid morphism $\widehat\mu$ has $\mu$ as the object component and $\pr_2$ as the morphism component.
\erem

\bedef\label{def:bibndl}
Let $\grpd{\xcG_A}{M_A},\ A\in\{1,2\}$ be Lie groupoids.\ A {\bf $(\xcG_1,\xcG_2)$-bibundle} is a manifold $\widehat P$ which carries the structure of a left $\xcG_1$-module $(\widehat P,\mu_1,\la_1\equiv\mlact)$ and that of a right $\xcG_2$-module $(\widehat P,\mu_2,\rho_2\equiv\mact)$,\ such that the two actions commute and each moment map is invariant with respect to the other action,\ {\it i.e.},\ we have,\ for all $(g_1,p,g_2)\in\xcG_1{}_{s_1}\hspace{-3pt}\x_{\mu_1}\hspace{-1pt}\widehat P{}_{\mu_2}\hspace{-3pt}\x_{t_2}\hspace{-1pt}\xcG_2$,
\bit
\item $(g_1\mlact p)\mact g_2=g_1\mlact(p\mact g_2)$;
\item $\mu_1(p\mact g_2)=\mu_1(p)$;
\item $\mu_2(g_1\mlact p)=\mu_2(p)$.
\eit
Whenever $(\widehat P,M_2,\mu_2,\mu_1,\la_1)$ is a (left) principal $\xcG_1$-bundle (with base $M_2\cong\widehat P/\xcG_1$),\ and $(\widehat P,M_1,\mu_1,\mu_2,\rho_2)$ is a (right) principal $\xcG_2$-bundle (with base $M_1\cong\widehat P/\xcG_2$),\ we call $\widehat P$ a ({\bf bi}){\bf principal $(\xcG_1,\xcG_2)$-bibundle},\ and depict it by the following $W$-shaped diagram:
\qq\label{diag:W-diagram}
\alxydim{@C=.5cm@R=1.cm}{ \xcG_1 \ar@{=>}[rd] & & \widehat P \ar[dl]_{\mu_1} \ar[rd]^{\mu_2} & & \xcG_2 \ar@{=>}[ld] \\ & M_1 & & M_2 & }\,.
\qqq
\exdef

\section{Conventions on Lie algebroids}

\bedef\label{def:Liealgbrd}
Let $M$ be a smooth manifold.\ A ({\bf real}) {\bf Lie algebroid over} $M$ ({\bf of rank} $N\in\bN^\x$) is a quintuple $(()\cE,M,\bR^N,\pi_\cE),\rho_\cE,[\cdot,\cdot]_\cE)$ composed of
\bit
\item a vector bundle $(\cE,M,\bR^N,\pi_\cE)$;
\item a vector-bundle morphism 
\qq\nn
\alxydim{@C=1.cm@R=1.cm}{ \cE \ar[r]^{\rho_\cE\quad} \ar[d]_{\pi_\cE} & \txT M \ar[d]^{\pi_{\txT M}}\\ M \ar@{=}[r]_{\id_M} & M}\,,
\qqq
termed the {\bf anchor} ({\bf map});
\item a binary operation $[\cdot,\cdot]_\cE\colo \G(\cE)\x\G(\cE)\to\G(\cE)$,
\eit
satisfying the following conditions:
\bit
\item $[\cdot,\cdot]_\cE$ is a Lie bracket;
\item $\forall_{\vep_1,\vep_2\in\G(\cE)}\ \forall_{f\in C^\infty(M;\bR)}\colo [\vep_1,f\lact\vep_2]_\cE=f\lact[\vep_1,\vep_2]+\rho_\cE(\vep_1)(f)\lact\vep_2$ (the {\bf Leibniz rule}).
\eit
\exdef

\beg\label{eg:tanalgbrd}
The {\bf tangent Lie algebroid of} $M$ is the canonical structure of a Lie algebroid on the tangent bundle $\pi_{\txT M}\colo \txT M\to M$ with anchor $\rho_{\txT M}=\id_{\txT M}$,\ and the standard Lie bracket $[\cdot,\cdot]_{\txT M}$ of vector fields on $M$.
\eeg

\beg\label{eg:tanLiealgbrd}
Fix a Lie groupoid $\xcG$.\ Let 
\qq\nn
\G(\txT\xcG)_{\rm R}=\bigl\{\  \xcV\in\G(\ker\,\txT s) \quad\vert\quad \forall_{g\in\xcG}\ \forall_{h\in s^{-1}(\{t(g)\})}\colo \txT_h r_{g}\bigl(\xcV(h)\bigr)=\xcV(h.g) \ \bigr\}
\qqq
be the set of {\bf right-invariant vector fields on} $\xcG$,\ defined in terms of the right-fibred action of Example \ref{eg:lr-fibreg-act},\ and let 
\qq\nn
\G(\txT\xcG)_{\rm L}=\bigl\{\  \xcV\in\G(\ker\,\txT t) \quad\vert\quad \forall_{g\in\xcG}\ \forall_{h\in t^{-1}(\{s(g)\})}\colo \txT_h l_{g}\bigl(\cV( h)\bigr)=\cV(g.h) \ \bigr\}
\qqq
be the set of {\bf left-invariant vector fields on} $\xcG$,\ defined in terms of the left-fibred action of Example \ref{eg:lr-fibreg-act}.\ The {\bf right tangent algebroid of} $\xcG$ is the vector bundle 
\qq\nn
\pr_1\colo \gt{gr}_{\rm R}:=\Id^*\ker\,\txT s=M\,{}_{\Id}\hspace{-3pt}\x_{\pi_{\txT\xcG}}\hspace{-1pt}\ker\,\txT s\too M
\qqq 
with anchor 
\qq\nn
\rho_{\gt{gr}_{\rm R}}=\txT t\circ\pr_2\,,
\qqq 
and Lie bracket 
\qq\nn
[\cdot,\cdot]_{\gt{gr}_{\rm R}}=\iota_{\rm R}^{-1}\circ[\cdot,\cdot]_{\txT\xcG}\circ(\iota_{\rm R}\x\iota_{\rm R})
\qqq
induced by the canonical ($\bR$-linear) isomorphism (note that -- by definition -- $\pi_{\txT\xcG}\circ\xcS=\Id$)
\qq\nn
\iota_{\rm R}\colo \G(\gt{gr}_{\rm R})\xrightarrow{\ \cong\ }\G(\txT\xcG)_{\rm R},\ \bigl(\id_M,\xcS\bigr)\longmapsto\txT_{\Id_{t(\cdot)}} r_\cdot\bigl(\xcS\bigl(t(\cdot)\bigr)\bigr)\,,
\qqq
with the obvious inverse
\qq\nn
\iota_{\rm R}^{-1}=\ev_{\Id}\colo \G(\txT\xcG)_{\rm R}\too\G(\gt{gr}_{\rm R}),\ \xcV\longmapsto\bigl(\cdot,\xcV\bigl(\Id(\cdot)\bigr)\bigr)\,.
\qqq

The {\bf left tangent algebroid of} $\xcG$ is defined analogously as the Lie algebroid
\qq\nn
\bigl(\gt{gr}_{\rm L}:=\Id^*\ker\,\txT t,M,\bR^{\rk(\ker\,\txT t)},\pr_1,\txT s\circ\pr_2,\iota_{\rm L}^{-1}\circ[\cdot,\cdot]_{\txT\xcG}\circ(\iota_{\rm L}\x\iota_{\rm L})\bigr)\,,
\qqq
where
\qq\nn
\iota_{\rm L}\colo \G(\gt{gr}_{\rm L})\xrightarrow{\ \cong\ }\G(\txT\xcG)_{\rm L},\ \bigl(\id_M,\xcS\bigr)\longmapsto\txT_{\Id_{s(\cdot)}} l_\cdot\bigl(\xcS\bigl(s(\cdot)\bigr)\bigr)
\qqq
is the left-handed counterpart of $\iota_{\rm R}$.
\eeg

\beg
The tangent Lie algebroid of the Lie groupoid $\txG$ of Example \ref{eg:grpasgrpd} is the tangent Lie algebra $\ggt\equiv{\rm Lie}(\txG)$ of the Lie group $\txG$,\ with the Lie-algebra bracket as the Lie bracket.
\eeg

\beg
The tangent Lie algebroid of the pair groupoid ${\rm Pair}(M)$ of Example \ref{eg:pairgrpd} is the tangent bundle $\txT M$ over the object manifold $M$,\ with the identity map as the anchor and the commutator of vector fields as the Lie bracket.
\eeg

\beg
An important example of a Lie algebroid is the (right) tangent algebroid of the action groupoid from Example \ref{eg:actgrpd}.\ It goes by the name of the {\bf action algebroid} and is denoted as $\ggt\,\lx_\la M$ (here,\ $\ggt\equiv{\rm Lie}\txG$ is the (tangent) Lie algebra of the Lie group $\txG$).\ One readily finds the convenient description
\qq\nn
\G\bigl(\txT(\txG\x M)\bigr)_{\rm R}=\bigl\{\  f^A\circ\la(\cdot)\lact R_A\bigl(\pr_1(\cdot)\bigr) \quad\vert\quad f^A\in C^\infty(M;\bR)\,,\quad A\in\ovl{1,\dim_\bR\,\txT_e\txG} \ \bigr\}\,,
\qqq
written out in terms of the basis right-invariant vector fields $R_A$ on $\txG$ associated with an arbitrary basis $\t=\{t_A\}_{A\in\ovl{1,\dim_\bR\,\txT_e\txG}}$ of $\txT_e\txG$.\ Upon applying $\iota_{\rm R}^{-1}$,\ we obtain
\qq\nn
\G\bigl(\txT(\txG\x M)\bigr)_{\rm R}\ni f^A\circ\la(\cdot)\lact R_A\bigl(\pr_1(\cdot)\bigr)\longmapsto f^A(\cdot)\lact t_A\in\G\bigl((e,\cdot)^*\ker\,\txT\pr_2\bigr)\equiv\G\bigl(\Id^*\ker\,\txT s\bigr)\,,
\qqq
which illustrates the bundle isomorphism
\qq\nn
\Id^*\ker\,\txT s\cong M\x\ggt\,.
\qqq
On global sections,\ the anchor evaluates as
\qq\nn
\G\a_{\ggt\,\lx_\la M}\colo \G\bigl(\Id^*\ker\,\txT s\bigr)\too\G(\txT M)\colo f^A(\cdot)\lact t_A\longmapsto-f^A(\cdot)\lact\xcK_A(\cdot)\,,
\qqq
with,\ on the right-hand side,\ the fundamental vector fields of $\la$, 
\qq\nn
\xcK_A\equiv\xcK_{t_A}=\txT_{(e,\cdot)}\la\bigl(-t_A,\brd0_{\txT M}(\cdot)\bigr)\,.
\qqq
Their bracket can be expressed in terms of the structure constants $f_{BC}^{\ \ \ A}$ of $\ggt$,\ associated with $\t$ through the structure equations
\qq\nn
[t_A,t_B]_{\ggt}=f_{AB}^{\ \ \ C}\,t_C\,,
\qqq
as
\qq\nn
[f_1^A\lact t_A,f_2^B\lact t_B]_{\ggt\,\lx_\la M} =\bigl(f_2^B\,\cK_B(f_1^A)-f_1^B\,\cK_B(f_2^A)-f_1^B\,f_2^C\,f_{BC}^{\ \ \ A}\bigr)\lact t_A\,.
\qqq
\eeg

\beg
The tangent Lie algebroid of the symplectic groupoid of Example \ref{eg:symplgrpd} is canonically isomorphic to the Lie algebroid $((\txT^*M,M,\bR^{\dim\,M},\pi_{\txT^*M}),\Pi^\#,\{\cdot,\cdot\})$,\ with the anchor given by the vector-bundle morphism $\Pi^\#\colo\txT^*M\to\txT M$  canonically induced by the Poisson bivector $\Pi$,\ and the Lie bracket 
\qq\nn
\{\cdot,\cdot\}\colo\G(\txT^*M)\x\G(\txT^*M)\too\G(\txT^*M),\ (\om_1,\om_2)\longmapsto\imath_{\Pi^\#(\om_1)}\sfd\om_2-\imath_{\Pi^\#(\om_2)}\sfd\om_1+\sfd\bigl(\om_1\wedge\om_2(\Pi)\bigr)\,,
\qqq
see \Rxcite{Chap.\,II, Prop.\,2.1}{Coste:1987}.
\eeg

Lie algebroids provide us with canonical examples of foliated vector-valued differential 1-forms.

\beg[\cite{Fernandes:2014}]\label{eg:MC}
Fix a Lie groupoid $\xcG$ with the tangent Lie algebroids $\ggt_{\rm H},\ H\in\{L,R\}$.\ The {\bf left-invariant Maurer--Cartan form on} $\xcG$ is the $\ker\,\txT t$-foliated 1-form with values in $\gt{gr}_{\rm L}$ given in 
\qq\nn
\alxydim{@C=1.5cm@R=1.5cm}{ \ker\,\txT t \ar[r]^{\ \MCL} \ar[d]_{\pi_{\txT\xcG}\rstr_{\ker\,\txT t}} & \gt{gr}_{\rm L} \ar[d]^{\pi_{\gt{gr}_{\rm L}}\equiv\pr_1}\\ \xcG \ar[r]_{s\ } & M}\,,\qquad\quad \MCL(g)(v)=\bigl(s(g),\txT_g l_{g^{-1}}(v)\bigr)\,,\quad v\in(\ker\,\txT t)_g\,.
\qqq
Similarly,\ the {\bf right-invariant Maurer--Cartan form on} $\xcG$ is the $\ker\,\txT s$-foliated 1-form with values in $\gt{gr}_{\rm R}$ given in  
\qq\nn
\alxydim{@C=1.5cm@R=1.5cm}{ \ker\,\txT s \ar[r]^{\ \MCR} \ar[d]_{\pi_{\txT\xcG}\rstr_{\ker\,\txT s}} & \gt{gr}_{\rm R} \ar[d]^{\pi_{\gt{gr}_{\rm R}}\equiv\pr_1}\\ \xcG \ar[r]_{t\ } & M}\,,\qquad\quad\MCR(g)(w)=\bigl(t(g),\txT_g r_{g^{-1}}(w)\bigr)\,,\quad w\in(\ker\,\txT s)_g\,.
\qqq
Thus,\ $\MCL(g)\in(\ker\,\txT t)^*_g\ox_\bR(\ker\,\txT t)_{\Id_{s( g )}}$ and $\MCR(g)\in(\ker\,\txT s)^*_g\ox_\bR(\ker\,\txT s)_{\Id_{t(g)}}$.
\eeg

\section{Useful properties of Lie groupoids,\ bisections,\ and Lie algebroids}\label{app:useful}

\berop\label{prop:BisAct-vs-str}
The left- and right-multiplications of $\xcG$ by $\bB$ from Def.\,\ref{def:bisec-act} have the following properties relative to the structure maps of $\xcG$ (written for arbitrary $(\b,g)\in\bB\x\xcG$):
\bit
\item[(i)] $s\circ L_\b=s\,,\quad s\circ R_\b=t_*(\b^{-1})\circ s$ ($s$ intertwines $L$ with the trivial action and $R$ with $t_*\circ\Inv$);
\item[(ii)] $t\circ L_\b=t_*\b\circ t\,,\quad t\circ R_\b=t$ ($t$ intertwines $L$ with $t_*$ and $R$ with the trivial action);
\item[(iii)] $L_\b\circ\Id=\b\,,\quad R_\b\circ\Id=\b\circ t_*\bigl(\b^{-1}\bigr)$;
\item[(iv)] $\Inv\circ L_\b=R_{\b^{-1}}\circ\Inv\,,\quad \Inv\circ R_\b=L_{\b^{-1}}\circ\Inv$ ($\Inv$ intertwines $L$ with $R\circ\Inv$);
\item[(v)] $L_\b\circ r_g=r_g\circ L_\b\,,\quad R_\b\circ l_g=l_g\circ R_\b$;
\item[(vi)] $r_g\circ R_\b=r_{L_\b(g)}\,,\quad l_g\circ L_\b=l_{R_\b(g)}$.
\eit
\eerop
\beroof
\bit
\item[Ad (i)] The first identity is trivial.\ For the second one,\ we compute explicitly,\ for any $h\in\xcG$,
\qq\nn
s\circ R_\b(h)&\equiv& s\bigl(h.\bigl(\b^{-1}\bigl(s(h)\bigr)\bigr)^{-1}\bigr)=t\bigl(\b^{-1}\bigl(s(h)\bigr)\bigr)\equiv t\bigl(\bigl(\b\circ(t_*\b)^{-1}\bigl(s(h)\bigr)\bigr)^{-1}\bigr)=s\circ\b\circ(t_*\b)^{-1}\bigl(s(h)\bigr)\cr\cr
&=&(t_*\b)^{-1}\circ s(h)=t_*\bigl(\b^{-1}\bigr)\circ s(h)\,.
\qqq
where the last equality follows from the homomorphicity of $t_*$.
\item[Ad (ii)] The second identity is trivial.\ For the first one,\ we compute explicitly,\ for any $h\in\xcG$,
\qq\nn
t\circ L_\b(h)\equiv t\bigl(\b\bigl(t(h)\bigr).h\bigr)=t\circ\b\bigl(t(h)\bigr)\equiv t_*\b\circ t(h)\,.
\qqq
\item[Ad (iii)] We find,\ for an arbitrary $m\in M$,
\qq\nn
L_\b\circ\Id(m)\equiv\b\bigl(t(\Id_m)\bigr).\Id_m=\b(m)
\qqq
and
\qq\nn
R_\b\circ\Id(m)\equiv\Id_m.\bigl(\b^{-1}\bigl(s(\Id_m)\bigr)\bigr)^{-1}=\bigl(\b^{-1}(m)\bigr)^{-1}\equiv\b\circ(t_*\b)^{-1}(m)=\b\circ t_*\bigl(\b^{-1}\bigr)(m)\,.
\qqq
\item[Ad (iv)] For the first identity,\ we compute explicitly,\ for any $h\in\xcG$,
\qq\nn
\Inv\circ L_\b(h)\equiv\bigl(\b\bigl(t(h)\bigr).h\bigr)^{-1}=h^{-1}.\b\bigl(t(h)\bigr)^{-1}\equiv h^{-1}.\bigl(\b^{-1}\bigr)^{-1}\bigl(s\bigl(h^{-1}\bigr)\bigr)^{-1}\equiv R_{\b^{-1}}\circ\Inv(h)\,.
\qqq
The second identity now follows by replacing $\b$ with $\b^{-1}$ in the one just proved,\ and subsequently sandwiching both sides of it between two copies of $\Inv$.
\item[Ad (v)] Take an arbitrary arrow $h\in s^{-1}(\{t(g)\})$ and calculate directly:
\qq\nn
L_\b\circ r_g(h)\equiv\b\bigl(t(h.g)\bigr).(h.g)=\b\bigl(t(h)\bigr).h.g\equiv r_g\bigl(\b\bigl(t(h)\bigr).h\bigr)\equiv r_g\circ L_\b(h)\,.
\qqq
The proof of the second identity is fully analogous. 
\item[Ad (vi)] For the first identity,\ take an arbitrary arrow $h\in s^{-1}(t_*\b(t(g)))$ and calculate directly
\qq\nn
r_g\circ R_\b(h)\equiv\bigl(h.\bigl(\b^{-1}\bigl(s(h)\bigr)\bigr)^{-1}\bigr).g=h.\bigl(\b^{-1}\bigl(t_*\b\bigl(t(g)\bigr)\bigr)\bigr)^{-1}.g\equiv h.\b\bigl(t(g)\bigr).g\equiv h.\bigl(L_\b(g)\bigr)\equiv r_{L_\b(g)}(h)\,.
\qqq
The second identity follows by replacing $\b$ with $\b^{-1}$ and $g$ with $g^{-1}$ in the one just proved,\ and subsequently using (iv).
\eit
\eroof

\brem\label{rem:BisAct-vs-str-simpl}
Upon evaluation on the respective arguments $h\in\xcG,\ m\in M,\ u\in s^{-1}(\{t(g)\}),\ v\in t^{-1}(\{s(g)\}),\ w\in s^{-1}(t_*\b(t(g))),\ y\in t^{-1}((t_*\b)^{-1}(s(g)))$ and in the shorthand notation $L_\b(h)\equiv\b\lact h,\ R_\b(h)\equiv h\ract\b$ and $t_*\b(m)\equiv\b\ulact m$,\ the functional identities of Prop.\,\ref{prop:BisAct-vs-str} take the following form:
\qq
& s(\b\lact h)=s(h)\,,\qquad\qquad s(h\ract\b)=\b^{-1}\ulact(s(h)\bigr)\,,&\label{eq:BisAct-vs-str-i}\\ \nn\\
&t(\b\lact h)=\b\ulact t(h)\,,\qquad\qquad t(h\ract\b)=t(h)\,,&\label{eq:BisAct-vs-str-ii}\\ \nn\\
&\b\lact\Id_m=\b(m)\,,\qquad\qquad\Id_m\ract\b=\b\bigl(\b^{-1}\ulact m\bigr)\,,&\label{eq:BisAct-vs-str-iii}\\ \nn\\
&(\b\lact h)^{-1}=h^{-1}\ract\b^{-1}\,,\qquad\qquad(h\ract\b)^{-1}=\b^{-1}\lact h^{-1}\,,&\label{eq:BisAct-vs-str-iv}\\ \nn\\
&\b\lact(u.g)=(\b\lact u).g\,,\qquad\qquad(g.v)\ract\b=g.(v\ract\b)\,,&\label{eq:BisAct-vs-str-v}\\ \nn\\
&(w\ract\b).g=w.(\b\lact g)\,,\qquad\qquad g.(\b\lact y)=(g\ract\b).y\,.&\label{eq:BisAct-vs-str-vi}
\qqq
\erem

\becor\label{cor:BisCon-vs-str}
The conjugation of $\xcG$ by $\bB$ from Def.\,\ref{def:bisec-act} has the following properties relative to the structure maps of $\xcG$ (written for arbitrary $(\b,g)\in\bB\x\xcG$):
\bit
\item[(i)] $s\circ C_\b=t_*\b\circ s$ ($s$ intertwines $C$ with $t_*$);
\item[(ii)] $t\circ C_\b=t_*\b\circ t$ ($t$ intertwines $C$ with $t_*$);
\item[(iii)] $C_\b\circ\Id=\Id\circ t_*\b$ ($\Id$ intertwines $C$ with $t_*$);
\item[(iv)] $\Inv\circ C_\b=C_\b\circ\Inv$ ($\Inv$ intertwines $C$ with itself);
\item[(v)] $C_\b\circ\txm=\txm\circ(C_\b\x C_\b)$ ($C$ distributes over $\txm$).
\eit
\ecor
\begin{corproof}
\bit
\item[Ad (i)] Point (i) of Prop.\,\ref{prop:BisAct-vs-str} implies $s\circ C_\b=s\circ R_{\b^{-1}}=t_*\b\circ s$.
\item[Ad (ii)] Point (ii) of Prop.\,\ref{prop:BisAct-vs-str} implies $t\circ C_\b=t\circ L_\b=t_*\b\circ t$.
\item[Ad (iii)] Point (iii) of Prop.\,\ref{prop:BisAct-vs-str} implies $C_\b\circ\Id\equiv R_{\b^{-1}}\circ(L_\b\circ\Id)=R_{\b^{-1}}\circ\b=R_{\b^{-1}}\circ(R_\b\circ\Id\circ t_*\b)=R_{\b\cdot\b^{-1}}\circ\Id\circ t_*\b=R_\Id\circ\Id\circ t_*\b=\Id\circ t_*\b$.
\item[Ad (iv)] Point (iv) of Prop.\,\ref{prop:BisAct-vs-str} implies $\Inv\circ C_\b\equiv(\Inv\circ L_\b)\circ R_{\b^{-1}}=R_{\b^{-1}}\circ(\Inv\circ R_{\b^{-1}})=(R_{\b^{-1}}\circ L_\b)\circ\Inv\equiv C_\b\circ\Inv$.
\item[Ad (v)] We calculate directly,\ in the notation of Remark \ref{rem:BisAct-vs-str-simpl} and using identities \eqref{eq:BisAct-vs-str-v} and \eqref{eq:BisAct-vs-str-vi} along the way,
\qq\nn
&&C_\b\circ\txm(g_2,g_1)\equiv\bigl(\b\lact(g_2.g_1)\bigr)\ract\b^{-1}=\bigl((\b\lact g_2).g_1\bigr)\ract\b^{-1}=\bigl(\bigl((\b\lact g_2)\ract\b^{-1}\bigr).(\b\lact g_1)\bigr)\ract\b^{-1}\cr\cr
&\equiv & \bigl(\bigl((\b\lact g_2)\ract\b^{-1}\bigr).\bigl(\bigl((\b\lact g_1)\ract\b^{-1}\bigr)\equiv C_\b(g_2).C_\b(g_1)\equiv\txm\circ(C_\b\x C_\b)(g_2,g_1)\,.
\qqq
\eit
\end{corproof}

\berop\label{prop:rasoR}
For every $g\in\xcG$ and every global bisection $\beta_g$ such that
\qq\label{eq:gBg}
g=\beta_g\bigl(s(g)\bigr)\,,
\qqq
the following identities hold true:
\bit
\item[(i)] $\beta_g((t_*\beta_g)^{-1}(t(g)))=g$;
\item[(ii)] the \emph{functional} identity:
\qq\nn
r_g=R_{\beta_g}\rstr_{s^{-1}(t(g))}\,,
\qqq 
and so also,\ for every $h\in s^{-1}(t(g))$,
\qq\label{eq:Trahg-TRahb}
\txT_h r_g=\txT_h R_{\beta_g}\rstr_{\ker\,\txT_h s}\,.
\qqq
\eit
\eerop
\beroof
For (i),\ evaluate the target map on both sides of the defining equality \eqref{eq:gBg},
\qq\nn
t(g)=(t_*\beta_g)\bigl(s(g)\bigr)\,,
\qqq
and subsequently invoke the invertibility of $t_*\beta_g\in\Diff(\xcG)$ to obtain the identity
\qq\nn
s(g)=(t_*\beta_g)^{-1}\bigl(t(g)\bigr)\,.
\qqq
Upon substituting the latter back in \eqref{eq:gBg},\ we arrive at the desired identity,
\qq\label{eq:BgtBgt}
g=\beta_g\bigl(s(g)\bigr)=\beta_g\bigl((t_*\beta_g)^{-1}\bigl(t(g)\bigr)\bigr)\,.
\qqq

For (ii),\ invoke \Reqref{eq:BisAct-vs-str-ii} to find the equality,\ for every $h\in s^{-1}(t(g))$,
\qq\nn
r_g(h)\equiv h.g=h.\bigl(\Id_{t(g)}\ract \beta_g\bigr)=\bigl(h.\Id_{t(g)}\bigr)\ract \beta_g=R_{\beta_g}(h)\,.
\qqq
Its differentiation yields formula \eqref{eq:Trahg-TRahb}.
\eroof

\bethe\cite{Rybicki:2002ALG}\label{thm:LieBis}
The tangent Lie algebra of the (Fr\'echet--)Lie group $\bB$ of bisections of a Lie groupoid $\xcG$ is isomorphic to the Lie algebra $\G_{\rm c}(\ggt_{\rm R})$ (resp.\ $\G_{\rm c}(\ggt_{\rm L})$) of smooth compactly supported sections of the tangent Lie algebroid $\ggt_{\rm R}$ (resp.\ $\ggt_{\rm L}$) of $\xcG$.
\ethe

\end{document}